\documentclass[11pt]
{article}

\usepackage{amsthm}
\usepackage[margin=30mm]{geometry}
\usepackage{amsmath}%
\usepackage{amsfonts}%
\usepackage{amssymb}%
\usepackage{amscd}
\usepackage{graphicx}
\usepackage{color}
\usepackage{xcolor}
\usepackage{bbm}
\usepackage[all,cmtip]{xy}
\usepackage{rotating}
\usepackage{url}
\usepackage{currfile}
\usepackage{enumerate}
\usepackage{subcaption}
\usepackage{enumitem}
\usepackage{wrapfig}
\usepackage{eurosym}
\usepackage{epigraph}
\usepackage{courier}
\definecolor{grau}{rgb}{0.1,0.1,0.1}
\usepackage[colorlinks=true, urlbordercolor=grau, linkcolor=grau, citecolor=grau, urlcolor=grau]{hyperref}
\usepackage{bibspacing}
\usepackage{titling}
\usepackage[bottom]{footmisc}

\numberwithin{equation}{section}

\newtheorem{theorem}[equation]{Theorem}

\newtheorem{lemma}[equation]{Lemma}
\newtheorem{corollary}[equation]{Corollary}

\theoremstyle{definition}
\newtheorem{definition}[equation]{Definition}

\newtheorem{remark}[equation]{Remark}
\newtheorem{remarks}[equation]{Remarks}

\newcommand{\bezout}{B\'{e}zout}

\newcommand{\kollar}{Koll\'ar}
\newcommand{\lojasiewicz}{\L{}ojasiewicz}

\newcommand{\RR}{\mathbbm{R}} 
\newcommand{\QQ}{\mathbbm{Q}} 
\newcommand{\ZZ}{\mathbbm{Z}} 
\newcommand{\CC}{\mathbbm{C}} 
\newcommand{\FF}{\mathbbm{F}} 
\renewcommand{\AA}{\mathbbm{A}} 
\newcommand{\PP}{\mathbbm{P}} 
\DeclareMathOperator*{\wo}{\backslash}
\DeclareMathOperator{\im}{im} 

\newcommand{\textdef}[1]{\textnormal{\textit{#1}}}

\newcommand{\incl}{\hookrightarrow} 
\newcommand{\st}{\ |\ } 
\newcommand{\Impl}{\Longrightarrow} 
\DeclareMathOperator{\codim}{codim} 

\newcommand{\mathsmall}[1]{{\begingroup\everymath{\scriptstyle}\small{#1}\endgroup}}

\newcommand{\eps}{\varepsilon}

\newcommand{\wt}[1]{\widetilde{#1}}

\newcommand{\blank}{\underline{\ \ }}

\newcommand{\rad}{\textnormal{rad}}

\newcommand{\dist}{\textnormal{dist}}

\newcommand{\Gal}{\textnormal{Gal}}
\newcommand{\Qbar}{\overline\QQ}
\newcommand{\OK}{\mathcal{O}_K}

\newcommand{\OL}{\mathcal{O}_L}
\newcommand{\divides}{\,|\,}
\newcommand{\ndivides}{\nmid}
\DeclareMathOperator{\Char}{char}

\DeclareMathOperator{\adj}{adj}

\renewcommand{\phi}{\varphi}
\renewcommand{\angle}{\measuredangle}

\newcommand{\Radical}[1]{{{#1}^{\textnormal{rad}}}}

\begin{document}

\title{Proofs by example}

\author{%
Benjamin Matschke%
\setcounter{footnote}{-1}%
\thanks{%
Boston University, matschke@bu.edu. 
Partially supported by Simons Foundation grant \#{}550023.%
}
\setcounter{footnote}{-1}%
}
\date{Sept.\ 1, 2019}




\maketitle


\begin{abstract}

We study the proof scheme ``proof by example'' in which a general statement can be proved by verifying it for a single example.
This strategy can indeed work if the statement in question is an algebraic identity and the example is ``generic''. 
This article addresses the problem of constructing a \emph{practical} example, which is sufficiently generic, for which the statement can be verified efficiently, and which even allows for a numerical margin of error.

Our method is based on diophantine geometry, in particular an arithmetic B\'ezout theorem, an arithmetic Nullstellensatz, and a new effective Liouville--\L{}ojasiewicz type inequality for algebraic varieties.
As an application we discuss theorems from plane geometry and how to prove them by example.

\end{abstract}


\makeatletter
\renewcommand*\l@section[2]{%
  \ifnum \c@tocdepth >\z@
    \addpenalty\@secpenalty
    \addvspace{0.5em \@plus\p@}%
    \setlength\@tempdima{1.5em}%
    \begingroup
      \parindent \z@ \rightskip \@pnumwidth
      \parfillskip -\@pnumwidth
      \leavevmode 
      \advance\leftskip\@tempdima
      \hskip -\leftskip
      {\normalsize 
      #1}\nobreak\hfil 
      \nobreak\hb@xt@\@pnumwidth{\hss #2}\par\vspace{0.5em}
    \endgroup
  \fi}
\makeatother


\renewcommand\contentsname{Contents}
\setcounter{tocdepth}{2}
\noindent
{\small
\tableofcontents}

\section{Introduction}

\subsection{Motivation}

Proof by example often refers to the 
venturesome
illusive
idea of a proving scheme, in which 
it is attempted
to prove a general statement of the form ``for all $x\in X$, $G(x)$ holds'' by verifying it for a single example $p\in X$: ``$G(p)$ holds, therefore $G(x)$ holds for all $x\in X$.''
Obviously this does not work in a large generality, as otherwise we could prove the statement ``all primes are even'' by verifying it for the single example ``$2$ is even''.
This explains why proof by example is usually considered as an inappropriate generalization, or as a logical fallacy~\cite{wikipediaProofByExample}.

On the other hand, there are situations in which sufficiently generic examples, random examples, or examples without any apparent particularity towards the statement, will lead to at least a heuristic that the statement in question is true in general.
For instance we may think of Thales' theorem from euclidean geometry (about $90^\circ$-angles in semi-circles).
If one makes a sufficiently generic sketch for Thales' theorem, which is as well sufficiently precise, then the angle in question will be measured and observed as quite close to $90^\circ$ in this sketch.
One will feel heuristically convinced that this angle should measure \emph{exactly} $90^\circ$, and not only for this sketch, but \emph{in general}.
Why is this?

Are there situations, in which a suitable example can be sufficient to prove the general statement?
Of course, we aim at positive answers that are useful in practice.

In this article we study general statements that can be 
phrased \emph{algebraically}.
For instance for Thales' theorem this can be done: If $\overline{AB}$ is a diameter of a circle that also contains the point $C$, then the set of all such points C 
is described by a circle equation of the form $f(A,B,C)=0$, 
and the statement ``$\angle BCA=90^\circ$'' translates also into a polynomial equation of the form $g(A,B,C)=0$ in the coordinates of the three points.

\paragraph{Algebraic proof by example scheme \textnormal{(exact version, first attempt)}.}
Suppose our general statement can be reformulated into the following algebraic one, namely that a polynomial~$g$ in $n$ variables vanishes on a given affine variety $X=V(f_1,\ldots,f_m)\subseteq\AA^n$.
A point $P\in X(\CC)=\{P\in\CC^n\st f_i(P)=0\ \forall i\}$ should be called \emph{sufficiently generic for~$g$}, if the following implication holds:
If $g(P)=0$ then $g|_{X(\CC)}=0$.

We would call $P$ the \emph{example} in this \emph{proof by example scheme}. (Later will will call this proof by example scheme \emph{exact}, as it deals with the exact equations $P\in X(\CC)$ and $g(P)=0$.)
As stated, $P$ may depend on $g$, and hence in principle one could simply take any non-root of $g|_{X(\CC)}$ if such exist.
Thus this first attempt for a proof by example scheme may look trivial.
However we aim at constructing sufficiently generic points~$P$ that only \emph{mildly} depend on~$g$, for instance on a certain ``arithmetic complexity'' of~$g$.

First ideas would be the following.

1. In the scheme-theoretic sense, $X$ has a generic point $*$.
Now, $g(*)$ equals $g$ in the affine coordinate ring of $X$, so $g(*)$ vanishes if and only if $g|_{X(\CC)}=0$.
This is a trivial equivalence, and hence of no practical use for us, since the equation $g(*)=0$ is not easier to verify than the identity $g|_{X(\CC)}=0$.

2. Case $X=\AA^n$: The Schwartz-Zippel lemma (see Section~\ref{secHistory}) states in a precise way that for a polynomial $g$ in $n$ variables over a field, the probability for $g(P)$ to vanish on a random point (using a suitable probability distribution) is small if $g$ is non-zero.
In other words, random points will with high probability serve as an example in the proof by example scheme.
As we want to be certain, a deterministic version is preferable.
This exists in form of Alon's Combinatorial Nullstellensatz; see Section~\ref{secHistory}.
It can be regarded as an exact proof by example scheme for $X=\AA^n$ that requires only a finite number of examples.

3. Case $X=\AA^1$: For univariate polynomials there exist classical root bounds, e.g.\ by Lagrange and by Cauchy, such that if $g(p)=0$ for suitably large $|p|$, $p\in\CC$, then $g=0$.
Here, the expression ``suitably large'' depends only on a degree bound for~$g$ and simple inequalities in the coefficients of~$g$.
This can be interpreted as a one-dimensional proof by example scheme; see Section~\ref{secHistory}.

\paragraph{Demands on a practical proof by example scheme.}
In practice, we want a proof by example scheme that solves the following issues:
\begin{enumerate}
\item
The example needs to be \emph{sufficiently generic}.
In the exact scheme as above this means: If $g(P)=0$ then $g|_{X(\CC)} = 0$.
\item
The example $P$ needs to be easy to construct.
This refers to flexibility in choosing the example as well as to computational efficiency. 
\item
The computation of $g(P)$ should be easy.
In the exact scheme as above, this means that the equality $g(P)=0$ can be efficiently verified.
\item
We want to allow a numerical margin of error, as 
in practice it may not be suitable to compute a point $P$ that lies exactly on~$X(\CC)$.
\end{enumerate}

\subsection{Main theorem}

The logarithmic Weil height (or simply the height) of a rational number $x=a/b$ with $\gcd(a,b)=1$ is defined as $h(x) = \log\max(|a|,|b|)$.
The height $h(f)$ of a polynomial $f\in\ZZ[x_1,\ldots,x_n]$ is defined as the maximal height of its coefficients.

Given some parameter $H$, the notation $a\ll_H b$ 
 ($a,b\in\RR$) will mean that there exists an explicit  
order-preserving bijection $\psi_H:\RR\to\RR$ depending only on $H$, such that $\psi_H(a)\leq b$.

For polynomials $f_1,\ldots,f_m\in\ZZ[x_1,\ldots,x_n]$, the variety $X=V(f_1,\ldots,f_m)\subseteq\AA^n$ can be thought of the set of common zeros of these polynomials, so for example $X(\CC) = \{P\in\CC^n\st f_i(P)=0~\forall i\}$.

We can now state a simplified version of our main theorem, which is the basis of this paper's proof by example scheme.
In Section~\ref{secRobustMainTheorem} we state it more generally over arbitrary number fields, for arbitrary valuations, and with explicit bounds. 

\begin{theorem}[Robust non-effective main theorem over $\QQ$ using the standard norm]
\label{thmRobustMainTheorem_IntroVersion}
Let $f_1,\ldots,f_m,g\in\ZZ[x_1,\ldots,x_n]$ be polynomials such that the affine variety $X = V(f_1,\ldots,f_m)\subseteq \AA^n$ is irreducible over~$\QQ$ and of dimension~$d$. 
Let 
\[
H = h(f_1)+\ldots+h(f_m)+h(g)+\deg f_1+\ldots+\deg f_m+\deg g+n.
\]
Let $P=(p_1,\ldots,p_n)\in \QQ^n$ be a point such that
\begin{equation}
\label{eqRobustMainTheoremIntro_boundsFor_hpi}
0\ \ll_{H}\ h(p_1)\ \ll_{H}\ \ldots\ \ll_{H}\ h(p_d).
\end{equation}
There exists an $\eps = \eps(H,h(p_d))>0$ such that the following holds.
If 
\begin{equation}
\label{eqRobustMainTheoremIntro_boundsFor_fiP_and_giP}
|f_i(P)|_\infty < \eps \ \ \textnormal{ for each } i=1,\ldots,m,\quad \textnormal{ and }\quad |g(P)|_\infty < \eps,
\end{equation}
then $g$ vanishes on $X(\CC)$.
\end{theorem}

In the theorem the only role of $H$ is to measure the ``arithmetic complexity'' of the given polynomials, on which the subsequent inequalities depend.
The point $P$ plays the role of the \emph{example} in the proof by example scheme.
The iterated height bounds~\eqref{eqRobustMainTheoremIntro_boundsFor_hpi} for the first~$d$ coordinates of~$P$ make the example \emph{sufficiently generic}. 
The absolute value bounds~\eqref{eqRobustMainTheoremIntro_boundsFor_fiP_and_giP} for $f_i(P)$ and $g(P)$ make $P$ an almost common zero of these polynomials, which can be interpreted as follows:
The example~$P$ has \emph{sufficient precision}; it is \emph{sufficiently close to} $X(\CC)$.

\paragraph{Interpretation as a numerical Nullstellensatz.}
Under these assumptions, the theorem asserts that~$g$ has to vanish on all common roots of $f_1,\ldots,f_m$.
Equivalently (via Hilbert's Nullstellensatz), $g$ has to lie in the radical ideal of $\langle f_1,\ldots,f_m\rangle$, $g\in \Radical{\langle f_1,\ldots,f_m\rangle}$, which means that there exist an integer $N\geq 1$ and polynomials $\lambda_1,\ldots,\lambda_m\in \QQ[x_1,\ldots,x_n]$ such that
\begin{equation}
\label{eqNullstellensatz_gN_as_polynomialCombination_of_fis}
g^N = \lambda_1 f_1 + \ldots + \lambda_m f_m.
\end{equation}
In this sense, Theorem~\ref{thmRobustMainTheorem_IntroVersion} can be regarded as a ``robust one-point Nullstellensatz'':
\begin{align*}
g\in \langle f_1,\ldots,f_m\rangle^\rad & \quad \Longleftrightarrow \quad 
g|_{X(\CC)} = 0 \\
& \quad \Longleftrightarrow \quad
``g(P)\approx 0\textnormal{ for some suitably generic $P$ close to $X$''}.
\end{align*}
Here, the first equivalence is Hilbert's Nullstellensatz, and the second one is Theorem~\ref{thmRobustMainTheorem_IntroVersion} or Theorem~\ref{thmRobustMainTheorem}.

The irreducibility assumption on $X$ is important, without it one could easily construct counter-examples where $f$ vanishes only on one of the components of $X$ and where $P$ is chosen $|\ .\ |_\infty$-close to this component.
Note that $X$ only needs to be irreducible over~$\QQ$, not over~$\Qbar$. 
A version of the main theorem for reducible varieties~$X$ is given in Section~\ref{secReducibleMainTheorem}.

\paragraph{New witness for $g|_{X(\CC)}=0$.}
Suppose we need a proof of the equality $g|_{X(\CC)}=0$ that is efficiently verifiable.
We call such a proof a \emph{witness} for $g|_{X(\CC)}=0$.
As in most cases, $|X(\CC)|=\infty$, checking $g(Q)=0$ for each point $Q\in X(\CC)$ separately is out of question.

Assume that $f_1,\ldots,f_m,g\in K[x_1,\ldots,x_n]$, where $K$ is a number field.
One way to obtain such a witness is by Hilbert's Nullstellensatz:
One may find polynomials $\lambda_1,\ldots,\lambda_m\in K[x_1,\ldots,x_n]$ such that~\eqref{eqNullstellensatz_gN_as_polynomialCombination_of_fis} holds.
The polynomials $\lambda_i$ can be encoded in finite memory, and to bound the size of the memory in terms of $f_1,\ldots,f_m$ and $g$ one can make use of an arithmetic Nullstellensatz, see Theorem~\ref{thmArithmNullstellensatzGeneral}.
Moreover, given $\lambda_1,\ldots,\lambda_m$,~\eqref{eqNullstellensatz_gN_as_polynomialCombination_of_fis} is efficiently verifiable.
Thus $(\lambda_1,\ldots,\lambda_m)$ is a witness for $g|_{X(\CC)}=0$.

A side application of 
Theorems~\ref{thmRobustMainTheorem_IntroVersion} and~\ref{thmRobustMainTheorem} 
 is that $P$ can be seen as a new kind of witness for $g|_{X(\CC)}=0$.
For this, note that $P$ can be represented in finite memory, and by~\eqref{eqRobustMainTheoremIntro_boundsFor_hpi} one can easily assure to pick a suitable $P$ such that the size of this memory is bounded in terms $n$, $m$, as well as the heights and the degrees of $f_1,\ldots,f_m,g$.
Moreover,~\eqref{eqRobustMainTheoremIntro_boundsFor_fiP_and_giP} can be verified efficiently.
To make the last step faster, one may want to switch to interval arithmetic, thus making the proof of~\eqref{eqRobustMainTheoremIntro_boundsFor_fiP_and_giP} numerical, and in this case, one should attach the details (i.e.{} the protocol) of the numerical computation to the witness.
Alternatively, these numerical issues can be avoided altogehter by working with non-archimedian valuations, see Theorem~\ref{thmRobustMainTheorem}.

\paragraph{Criteria for irreducibility.}
An affine variety $X=V(f_1,\ldots,f_m)$, $f_i\in K[x_1,\ldots,x_n]$, is irreducible over $K$ if and only if $\langle f_1,\ldots,f_m\rangle^\rad\subseteq K[x_1,\ldots,x_n]$ is a prime ideal.
Assuming $|K|=\infty$, a useful example of irreducible affine varieties are those with a rational parametrization, i.e. when $X$ is of the form $X=\{(\rho_1(t),\ldots,\rho_n(t))\in K^n\st t\in K^d\}$, where $\rho_1,\ldots,\rho_n \in K(t_1,\ldots,t_d)$ are rational functions.  
More generally, images of irreducible varieties under rational maps are irreducible.
A version of the main theorem for reducible varieties~$X$ is given in Section~\ref{secReducibleMainTheorem}.

\paragraph{Determining the dimension by example.} 
Theorem~\ref{thmRobustMainTheorem_IntroVersion} requires the knowledge of the dimension of~$X$.
In case one has a good guess for $\dim X$, say $d$, then Section~\ref{secDimensionByExample} offers an approach to try prove $\dim X = d$ ``by example''.
Informally, this approach requires a sufficiently generic point $P\in K^n$ that is close enough to $X$ in the sense that $|f_i(P)|$ is small. 
If then a certain determinant can be suitably bounded away from zero it follows that $\dim X=d$.
See Theorem~\ref{thmDimensionByExample} for details.
We will also discuss that the determinant bound is reasonably week, see Corollary~\ref{corDimensionByExampleConverse}.

\paragraph{Deciding whether or not $g|_X=0$.}
The main theorem gives an only sufficient criterion of the form:
If $|g(P)|_v$ is small then $g|_X=0$.
There are situations in which we want to decide from a single sufficiently generic example $P$ close to $X$, whether or not $g|_X=0$.
In Section~\ref{secDichotomy} we discuss a dichotomy theorem that makes this decision based on knowing $|g(P)|$ to a high enough precision.

\paragraph{Short-cut through the paper.}
Readers who are only interested in the proof of the main theorem (Theorem~\ref{thmRobustMainTheorem}) are suggested to read (if necessary) the introductory Section~\ref{secPrerequisitesHeights} on heights, the \bezout{} version of the arithmetic Nullstellensatz (Theorem~\ref{thmArithmNullstellensatzGeneral}), Lemma~\ref{lemXcontainedInCoordinateHyperplane} and its corollary about heights of varieties contained in coordinate hyperplanes, the effective \lojasiewicz{} inequality for empty varieties (Theorem~\ref{thmLojasiewiczForEmptySet}), and finally the proof of the main theorem in Section~\ref{secRobustMainTheorem}.

\paragraph{Proof of Thales' theorem by example.}
In Section~\ref{secThales} we will give two proofs by example of the above mentioned Thales' theorem.
Here is a non-technical summary of the first proof:
\begin{enumerate}
\item
We may assume that the given circle is the unit circle and that the diameter $AB$ lies on the $x$-axis.
\item
As there is one degree of freedom left, we choose the $x$-coordinate of $C$ with $13$ digits behind the comma to make the construction sufficiently generic.
\item
We construct the point $C$ up to 1300 digits of precision to make the construction sufficiently exact.
\item
We verify that the angle at~$C$ 
is close enough to~$90^\circ$.
This finishes the proof by example.
\end{enumerate}
In Section~\ref{secPlaneGeometry} we also discuss more general theorems from plane geometry as examples for the proof by example scheme.

\paragraph{Basic idea.} 
Suppose a polynomial $g\in\ZZ[x]$, $g(x)=\sum_i a_i x^i$, has coefficients of bounded absolute value, say $|a_i|\leq r$.
Suppose our task is to prove that $g=0$.
We claim that for any $p\in\ZZ$ with $p> r+1$, the following implication holds:
If $g(p) = 0$ then $g=0$.
This means that in order to show $g=0$, it is enough to check this equality for the example~$p$.
This  
is a proof by example.
For instance, let $g(x) = 14x^2+4x+4$. Plugging in $100$ yields $g(100) = 140404$.
And we see that one can read of the coefficients of $g$ from its value at $p=100$.
This is made precise by the following root bounds for univariate polynomials.

\subsection{History}
\label{secHistory}

\paragraph{Bounds for roots of univariate polynomials.}
Lagrange (1798) proved the following 
root bound for a monic polynomial $g(x)=\sum_{i=0}^{\deg g-1} a_ix^i\in\CC[x]$: 
If $g(p)=0$ then $|p|\leq \max(1,\sum_{i=0}^{\deg g-1}|a_i|)$.
We can interpret Lagrange's root bound as the world's first proof by example scheme.
But before, let us define the height of a polynomial $g\in\ZZ[x]$ as $h(g):=\log(\max(|a_i|))$ if $g(x)=\sum a_ix^i \neq 0$ and $h(0):=0$.
\emph{Lagrange's proof by example scheme} now reads: 
\begin{quotation}
It $g\in\ZZ[x]$ and $p\in \CC$ satisfy $|p|> e^{h(g)}\deg g$, then $g(p)=0$ implies $g=0$.
\end{quotation}
This means that if the degree and height of $g\in\ZZ[x]$ are bounded and $p$ is sufficiently large in terms of these bounds, then $p$ serves as sufficiently generic example for the statement $g=0$.

\medskip

Similarly, Cauchy (1829) proved an often stronger root bound for monic polynomials $g(x)=\sum_{i=0}^{\deg g-1} a_ix^i\in\CC[x]$: 
If $g(p)=0$ then $|p|\leq 1 + \max \{|a_i|\st 0\leq i <\deg g\}$. We obtain \emph{Cauchy's proof by example scheme}:
\begin{quotation}
If $g\in\ZZ[x]$ and $p\in\CC$ satisfy $|p|>1+e^{h(g)}$, then $g(p)=0$ implies $g=0$.
\end{quotation}
Note that Cauchy's scheme does not need any degree bound for~$g$.
In any case we already see that proofs by example can indeed work, and that the basic ideas in ambient dimension one are classical.

\paragraph{Probabilistic approach.}
The Schwartz--Zippel lemma states: 
If $g\in\FF[x_1,\ldots,x_n]$ is a non-zero polynomial of degree $D$ over any field $F$, and $S$ is a finite subset of $S\subseteq \FF$, and $p_1,\ldots,p_n$ are elements from $S$ that are chosen independently and uniformly at random, then
\[
P[g(p_1,\ldots,p_n)=0] \leq \frac{D}{|S|}.
\] 
That is, the probability that $g$ vanishes at a random point is small.
This was proved independently by Schwartz~\cite{Schwartz80fastProbAlgoForVerificationOfPolyIdentities}, Zippel~\cite{Zippel79probabilisticAlgoForSparsePolynomials}, DeMillo and Lipton~\cite{DeMilloLipton78probabilisticRemarkOnAlgebraicProgramTesting}, and already in 1922 in a special case by Ore~\cite{Ore22ueberHoehereKongruenzen}.
It means that for \emph{disproving} a statement of the form $g=0$, a random example $P$ will yield a  contradiction of the form $g(P)\neq 0$ with high probability. 

At the same time, the Schwartz-Zippel lemma can be used to \emph{prove} a statement of the form $g=0$ if we test all sample points in the finite sample space (for this to work we need $|S|>D$).
This already yields as a corollary a special case of the following combinatorial Nullstellensatz.

\paragraph{Combinatorial Nullstellensatz.}

The (weak) combinatorial Nullstellensatz states the following.
If $g\in \FF[x_1,\ldots,x_n]$ is a polynomial of degree $D$ over any field $\FF$, and $S_1,\ldots,S_n\subseteq \FF$ are finite subsets of $\FF$ of cardinality $|S_i|=D+1$, then
\[
g(S_1\times\cdots\times S_n)=0\quad \Impl\quad g=0.
\] 
Alon~\cite{Alon99combinatorialNullstellensatz} proved an equivalent but more flexible form, which has many nontrivial applications in various branches of combinatorics.

As stated, the weak form can be considered as a proof by \emph{examples} schemes for $X=\AA^n$: 
To prove statements of the form $g=0$ we only need to verify that $g(P)=0$ for the $(D+1)^n$ points in the discrete cube $S_1\times\cdots\times S_n$.

\newcommand{\krsub}{{\textnormal{kr}}}

\paragraph{Kronecker substitution.}
Let $g\in \FF[x_1,\ldots,x_n]$ be a polynomial of degree $\deg g < D$ over any field~$\FF$.
\emph{Kronecker substitution}~\cite{kronecker1882grundzugeEinerArithmetischenTheorie} makes $g$ into a univariate polynomial $g_\krsub$ via the simple substitution
\[
g_\krsub(z):=g(z,z^D,z^{D^2},\ldots,z^{D^{n-1}})\in K[z].
\]
Then one has the equivalence: $g = 0$ if and only if $g_\krsub = 0$.
Using Lagrange's or Cauchy's root bounds (or rather our interpretation as proof by example schemes) together with Kronecker substitution yields a proof by example schemes for multivariate polynomials $g\in\ZZ[x_1,\ldots,x_n]$:
To any such given $g$, take any $p\in\CC$ with $|p|>1+\exp^{h(g_\krsub)}$ (and note that $h(g_\krsub)=h(g)$). Then $g_\krsub(p)=0$ implies $g_\krsub=0$, which implies $g=0$.

The advantage of this over the (weak) combinatorial Nullstellensatz is that we only need to test $g$ at the single example $(p,p^D,\ldots,p^{D^{n-1}})$; the disadvantage is that this example has necessarily large height (see Section~\ref{secPrerequisitesHeights} for definitions).

\paragraph{Gr\"obner bases.}
Let $X=V(I)$, where $I=\langle f_1,\ldots,f_m\rangle$ is an ideal of $K[x_1,\ldots,x_n]$, $K$ a number field.
Showing that a given $g\in K[x_1,\ldots,x_n]$ vanishes on $X(\Qbar)$ is equivalent to $g\in\Radical{I}$, where $\Radical{I}$ is the radical ideal of~$I$.
With Gr\"obner bases (Buchberger~\cite{Buchberger2006phd}) and algorithms for computing radical ideals (e.g.\ Eisenbud--Huneke--Vasconcelos~\cite{eisenbudHunekeVasconcelos92directMethodsPrimDecomp}, Krick--Logar~\cite{krickLogar91algorithmForRadicalIdeal}, and Kemper~\cite{kemper02calculationOfRadicalIdealInPosChar}) one can decide the membership relation $g\in \Radical{I}$ algorithmically.

In general, computing Gr\"obner bases can quickly 
become infeasible:
Mayr and Meyer~\cite{MayrMeyer82complexityOfWordProblemsAndPolynomialIdeals} proved that the uniform word problem for commutative semigroups is exponential space complete and they rephrased that in terms of an ideal membership problem.
For various complexity bounds for computing Gr\"obner bases and ideal membership we refer to 
Mayr and Ritscher~\cite{MayrRitscher13dimDependentBoundsForGroebnerBases} and
Mayr and Toman~\cite{MayrToman17complexityOfMembershipProblemsOfPolynomialIdeals}.

In practice, the Gr\"obner bases approach for showing $g|_{X(\Qbar)}=0$ is presumably still more practical algorithmically in comparison to the proof by example scheme of this paper, as the current bounds for $h(p_d)$ in~\eqref{eqRobustMainTheoremIntro_boundsFor_hpi} and~\eqref{eqRobustMainTheorem_assumptionOnHpi} are not practical, in particular when~$d=\dim X$ becomes large, and as Gr\"obner bases don't require $I$ to be a prime ideal.
If on the other hand the task is to show that $g$ vanishes on a particular component of $X$, the proof by example scheme in Corollary~\ref{corRobustMainTheoremForReducibleX} can become favorable.

\paragraph{Criteria for algebraic independence.}
In this paper we use the same height machinery from diophantine geometry that is used in transcendence theory for example for criteria for algebraic independence.
We have similar aims: Instead of searching for algebraically independent coordinates we search for a relaxation thereof that is more practical for our purposes.
See Section~\ref{secSufficientlyGenericPoints} for more details and pointers to the literature.

\subsubsection{Related topics}

\newcommand{\PIT}{\textnormal{PIT}}
\newcommand{\EVAL}{\textnormal{Eval}}

\paragraph{Polynomial identity testing.}
Computationally, a polynomial $g\in \FF[x_1,\ldots,x_n]$ over a computable field $\FF$ (such as $\QQ$ or finite fields) can be represented as an \emph{arithmetic circuit}, which takes as input variables and constants from $\FF$ and which is allowed to add and multiply any two or more expressions that it already has computed.
The \emph{polynomial identity testing problem} $\PIT_\FF$ is the decision problem of 
whether a given arithmetic circuit over $\FF$ represents the zero polynomial.

The \emph{evaluation problem} $\EVAL_\FF$ is the decision problem of deciding whether a polynomial $g\in \FF[x_1,\ldots,x_n]$ given as an arithmetic circuit evaluates at a given point $P\in \FF^n$ to zero: $g(P)=0$.
Clearly, $\EVAL_\FF$ reduces to $\PIT_\FF$, as one can simply substitute $P$ for the variables in the arithmetic circuit that represents~$g$, to then check whether the resulting circuit represents zero.
Using Kronecker substitution one can reduce $\PIT_\FF$ to $\EVAL_\FF$ if $\FF$ has characteristic~$0$.
In particular $\PIT_K$ and $\EVAL_K$ are polynomial-time equivalent for number fields~$K$; see Allender--B{\"u}rgisser--Kjeldgaard-Pedersen--Miltersen~\cite{allender2009complexityOfNumericalAnalysis}.

On the other hand, assume that the polynomial $g$ is given as an arithmetic circuit.
Checking $g=0$ by expanding the arithmetic circuit has exponential worst-case complexity in the input size, whereas evaluating $g(P)$ at a point of bounded height has polynomial complexity.
The Schwartz--Zippel lemma from above applied with $S = \{1,2,\ldots,3D\}$ yields thus an efficient probabilistic algorithm (in the complexity class co-RP, which is a subset of BPP) to test whether~$g=0$: It simply outputs whether $g(P)=0$ where $p_1,\ldots,p_n$ are drawn independently and uniformly at random from~$S$.
The quest for a suitably generic example in the proof by example scheme is therefore related to the problem whether this probabilistic algorithm can be derandomized:
Does $\PIT_\QQ$ lie in~$P$?
For a comprehensive survey on $\PIT$ and arithmetic circuits in general we refer to to Shpilka and Yehudayoff~\cite{ShpilkaYehudayoff10arithmeticCircuits}.

\paragraph{Derandomization of probabilistic algorithms.}
In computational complexity theory there exists a dichotomy between hardness and randomness.
It roughly states that computationally hard functions exist if and only if one can derandomize any polynomial time random algorithm.
Impagliazzo and Wigderson~\cite{ImpagliazzoWigderson97PeqBPPIfErequiresExpCircuits} proved 
one particular instance: If the satisfiability problem (SAT) cannot be solved by circuits of size $2^{o(n)}$, then BPP$=$P.
An underlying idea is that the existence of hard problems can be turned into pseudo-random generators, which suitably fool efficient algorithms.
Then, a given efficient probabilistic algorithm can be run a polynomial number of times using the pseudo-random generator for different seeds, and a majority vote among the outcomes will yield the desired derandomized deterministic algorithm.

There is also a partial converse by Kabanets and Impagliazzo~\cite{KabanetsImpagliazzo04derandomizingPITs}, who proved that a derandomization of the probabilistic algorithm for $\PIT_\QQ$ embodied by the Schwartz--Zippel lemma would imply certain explicit non-trivial lower complexity bounds. 

\paragraph{Arithmetic complexity.}
Valiant~\cite{Valiant79completenessClassesInAlgebra} defined arithmetic analogues of the classical complexity classes P and NP, called VP$_\FF$ and VNP$_\FF$.
He showed that the determinant is VP$_\FF$-complete respect to quasi-polynomial projections, and that the permanent is VNP$_\FF$-complete. 
Thus $\PIT_\FF$ reduces quasi-polynomially to the decision problem of whether a determinant with entries from the set $\{x_1,\ldots,x_n\}\cup \FF$ is zero.
The latter polynomials form thus natural benchmark for measuring the strength of future improvements of the proof by example scheme.

\paragraph{Fukshansky's theorem.}
Fukshansky~\cite{Fukshansky10algebraicPointsSmallHeightMissingUnionOfVarieties} showed the existence of a point $P$ of bounded height in any given linear space $V\subseteq \AA^n$ that lies in the complement of a union $Z$ of varieties over a number field~$K$.
Thus, if our aim is to show that a polynomial $g\in K[x_1,\ldots,x_n]$ vanishes on $V$, we can choose $Z$ as the union of all hypersurfaces not containing $V$ that are defined by polynomials over $K$ whose heights and degrees are bounded by the same bounds that are known for~$g$.
If Fukshansky's theorem was practically constructive, then the obtained point $P$ could be an example in our proof by example scheme.
For our particular union $Z$, a simple Fukshansky type result could be obtained by applying Siegel's lemma from diophantine approximation to construct a basis for $V(K)$ of bounded height, and then apply the above Kronecker substitution to $g|_{V(K)}$ in order to find a suitably point $P$ via Cauchy's root bound.

\paragraph{Taylor expansion via divided polynomial algebra.}
The divided polynomial algebra $\Gamma_\ZZ[\alpha]$ is a $\ZZ$-algebra, which is generated as a free $\ZZ$-module by $\alpha^{(0)},\alpha^{(1)},\alpha^{(2)},\ldots$, and the multiplication on these generators is defined by
$\alpha^{(k)}\cdot \alpha^{(\ell)} := \binom{k+\ell}{k}\alpha^{(k+\ell)}$.
We also write $\alpha^{(0)}=1$, as it is the multiplicative unit, and $\alpha^{(1)} = \alpha$.
We consider $\Gamma_\ZZ[\alpha]$ as a subalgebra of the polynomial ring $\QQ[t]$ via the identification $\alpha^{(k)} = t^k/k!$.
$\Gamma_\ZZ[\alpha]$ appears naturally in the study of the cohomology ring $H^*(J(S^n);\ZZ)$ of the James reduced product of even dimensional spheres, see e.g.\ Hatcher~\cite[Prop. 3.22]{Hatcher02AT}.
Similarly, let $\Gamma_\ZZ[\alpha_i]$ be a copy of $\Gamma_\ZZ[\alpha]$ with generators $\alpha_i^{(k)}$, $k\geq 0$.
Consider $\Gamma:=\Gamma_{\Qbar}[\alpha_1,\ldots,\alpha_n]:=\Qbar\otimes \Gamma_\ZZ[\alpha_1]\otimes\cdots\otimes\Gamma_\ZZ[\alpha_n]$.

Suppose we are given a polynomial $g\in K[x_1,\ldots,x_n]$ over a number field $K$ (or any field of characteristic zero) as well as a point $P=(p_1,\ldots,p_n)\in\Qbar$.
We substitute in $g$ for $x_i$ the sum $p_i+\alpha_i\in\Gamma$ with $\alpha_i=\alpha_i^{(1)}$, and obtain a simplified Taylor expansion of~$g$ without denominators,
\[
g(p_1+\alpha_1,\ldots,p_n+\alpha_n) = \sum_{{\bf k}} \partial^{({\bf k})}g(P)\; \alpha^{({\bf k})}\quad  \in\Gamma,
\]
where the sum runs over all ${\bf k} = (k_1,\ldots,k_n)\in (\ZZ_{\geq 0})^n$ which has only finitely many non-zero summands, 
$\partial^{({\bf k})}g := \tfrac{d^{{\bf k}}}{dx^{{\bf k}}} g$ and $\alpha^{({\bf k})} := \alpha_1^{(k_1)}\cdots \alpha_n^{(k_n)}$. 
Clearly, $g=0$ if and only if $g(p_1+\alpha_1,\ldots,p_n+\alpha_n)=0$ in $\Gamma$.

In order to work in a finite dimensional $\Qbar$-vector space, let $J=\langle \alpha_1^{(1)},\ldots,\alpha_n^{(1)}\rangle$ be the augmentation ideal of~$\Gamma$.
Suppose we know in advance that $\deg g<D$.
Then if will be enough to consider $g(p_1+\alpha_1,\ldots,p_n+\alpha_n)$ in the quotient $\Gamma/J^D$: Then again, $g=0$ if and only if this congruence class vanishes.

\paragraph{Strong law of small numbers.}
According to Richard Guy, ``there are not enough small numbers to satisfy all the demands placed on them''.
This is usually used in a seemingly contrary sense as in this paper: One cannot draw conclusions from looking at the first few examples. What holds for small numbers may not reflect the general picture for arbitrarily large numbers.
A great example is the (false) conjecture the prime counting function $\pi(x)$ is bounded by the logarithmic integral $\textnormal{Li}(x)$ for all $x\geq 3$? 
Indeed it holds true for all explicit examples where one can test it so far, however the relation has been shown to switch infinitely often, and $\pi(x) > Li(x)$ will hold for some $x\leq 10^{314}$, as proved by Littlewood~\cite{Littlewood14surLaDistributionDesNombresPremiers} and Skewes~\cite{Skewes33piMinusLiI,Skewes33piMinusLiII}. 

On the other hand, on inspecting the sequences $1,2,4,8,16,32$ and $1,1,2,3,5,8,13$, the reader probably thinks about powers of $2$ and the Fibonacci sequence.
Why is this? 
In the uncountable set of all continuations of these two sequences, only countably many can be put in words, and perhaps one tries automatically to find the ``simplest'' continuation.
However this is not a fully satisfactory answer.

Another argument against a universal law of small numbers is that there are many examples of theorems and conjectures that were first thought of because of examples; famous instances are monstrous moonshine, the Birch and Swinnerton-Dyer conjecture, and Montgomery's pair correlation conjecture.

\section{Prerequisites from diophantine geometry}
\label{secPrerequisites}

\subsection{Heights}
\label{secPrerequisitesHeights}

In this section we fix notions of heights for numbers, vectors, matrices, polynomials, and varieties, that we will use in this paper, as well as important relations between these heights.

In the literature, there exist many different notions of the height of a variety, e.g.\ Philippon~\cite{philippon1986criteresIndependanceAlgebrique,philippon1991hauteursAlternativesI,
philippon1994hauteursAlternativesII,philippon1995hauteursAlternativesIII}
and Bost, Gillet, Soul\'e~\cite{BostGilletSoule94heightsOfProjVarietiesAndPosGreenForms}, each with different advantages.
We choose to follow the conventions of Krick, Pardo, Sombra~\cite{KrickPardoSombra01sharpEstimatesForArithmeticNullstellensatz}, because we will use their effective arithmetic Nullstellensatz as well as their version of the arithmetic \bezout{} inequality.

Convention: In order to avoid unnecessary case distinctions, we define $\max\emptyset := -\infty$ and $\log 0 := -\infty$.

\paragraph{Valuations.}
Let $\QQ$ denote the integers, $K$ be a number field, and $\Qbar$ their algebraic closure.

Let $M_\QQ = \{p\st p\textnormal{ prime}\}\cup\{\infty\}$ denote the set of finite and infinite places of~$\QQ$.
They correspond to the standard normalized absolute values of~$\QQ$, given by $|\tfrac{a}{b}p^n|_p = p^{-n}$ if $p\nmid ab$, and $|a|_\infty = \sqrt{a^2}$ the usual absolute value.  

Let $M_K$ be the set of finite and infinite places of~$K$, whose elements $v$ we identify with the standard normalized absolute values $|\ .\ |_v:K\to \RR$ with the convention $|0|_v=0$.
The normalization means that any $v\in M_K$ extends (exactly) one absolute value $v_0\in M_\QQ$, and we write $v\divides v_0$.
Let $M_K^\infty \subset M_K$ denote the subset of infinite places of~$K$, which correspond to the archimedian absolute values, i.e. those $v$ such that $v\mid \infty$.

For $v\in M_K$, let $K_v$ denote the completion of $K$ with respect to $|\ .\ ]_v$.
For $v_0\in M_\QQ$, $\QQ_{v_0}$ is thus a completion of~$\QQ$, and we denote by~$\CC_{v_0}$ the algebraic closure of~$\QQ_{v_0}$.
Then $v$ extends naturally to a valuation on~$K_v$, and $v_0$ to a valuation on $\QQ_{v_0}$ and on~$\CC_{v_0}$.

If $v\divides v_0$, where $v\in M_K$ and $v_0\in M_\QQ$, there is a field homomorphism $\sigma_v:K_v\incl \CC_{v_0}$ such that $|x|_v = |\sigma_v x|_{v_0}$.
Indeed there are $N_v:=[K_v:\QQ_{v_0}]$ such homomorphisms, any we may choose $\sigma_v$ to be any of them.
For a polynomial $f$ with coefficients in $K$ or $K_v$, let $f^{\sigma_{v}}=\sigma_vf$ denote the same polynomial but with coefficients pushed forward to $\CC_{v_0}$, such that $f^{\sigma_{v}}$ has coefficients in~$\CC_{v_0}$.

\paragraph{Logarithmic Weil heights.}

For a non-empty finite set $A\subset K$ and $v\in M_K$ define
\begin{equation}
\label{eqDefOf_Av_and_hvA}
|A|_v := \max_{a\in A} |a|_v \qquad \textnormal{ and } \qquad h_v(A) := \max(0,\log|A|_v).
\end{equation}
If $v\divides v_0$, where $v\in M_K$ and $v_0 \in M_\QQ$, we define as above $N_v := [K_v:\QQ_{v_0}]$. 
There is a product formula $\prod_{v\in M_K} |x|_v^{N_v} = 1$ for each $x\in K\wo\{0\}$, and $[K:\QQ] = \sum_{v\divides v_0} N_v$ for each $v_0\in M_{\QQ}$.

For $A\subset\Qbar$, let $K$ be a number field containing $A$ and define
\[
h(A):= \tfrac{1}{[K:\QQ]}\sum_{v\in M_K} N_v h_v(A),
\]
which does not depend on the choice of~$K$.
If $A=\{x\}$, we also write $h(x):=h(A)$.

If $f_1,\ldots,f_m$ is a set of polynomials (possibly constant) with coefficients in $\Qbar$, let $A$ denote its set of all their coefficients, and define $|f_1,\ldots,f_m|_v := |A|_v$, $h_V(f_1,\ldots,f_m)_v := h_v(A)$, and $h(f_1,\ldots,f_m) := h(A)$.
Similarly, if $u\in K^n$ is a vector, let $A$ be the set of its entries and we define $|u|_v:=|A|_v$, $h_v(u):=h_v(A)$ and $h(u):=h(A)$; and similarly for matrices $M\in K^{n\times m}$.

If $\wt A\subset\CC_{v_0}$ for some $v_0\in M_\QQ$, we define 
\[
|\wt A|_{v_0} := \max_{a\in \wt A} |a|_{v_0} \qquad \textnormal{ and } \qquad h_{v_0}(\wt A) := \max(0,\log|\wt A|_{v_0})
\]
in analogy with~\eqref{eqDefOf_Av_and_hvA}.
Given $A\subset K$ and $v\divides v_0$ with $v\in M_k$ and $v_0\in M_\QQ$, it is convenient to note that $h_v(A) = h_{v_0}(\sigma_v(A))$, and hence 
\[
h(A) = \tfrac{1}{[K:\QQ]}\sum_{v_0\in M_\QQ}\sum_{\sigma:K\to\CC_{v_0}} h_{v_0}(\sigma A),
\]
where the inner sum is over all norm-preserving homomorphisms~$\sigma$ that respect some $v\in M_K$ over~$v_0$.

\paragraph{Chow forms.}
Let $Y\subset \PP^n$ be a $d$-dimensional projective variety defined over~$K$.
Let us consider $d+1$ groups $U_0,\ldots,U_d$ each of $n+1$ variables, $U_i = \{u_{i0},\ldots,u_{in}\}$.
Each $U_i$ parametrizes a general projective hyperplane $\{[x_0:\ldots:x_n\}\in \PP^n \st u_0x_0+\ldots+u_nx_n = 0\}$.
Thus $U_0,\ldots,U_d$ parametrize $d+1$ general projective hyperplanes in~$\PP^n$.
Let $Z\subseteq\AA^{(d+1)(n+1)}(\Qbar)$ the the locus of all $(d+1)$-tuples of hyperplanes with $\Qbar$-coefficients whose intersection contains a point of $X(\Qbar)$.
Chow and van der Waerden~\cite{ChowVanDerWaerden37zurAlgebraischenGeometrie} proved that~$Z$ is a hypersurface, which is defined by a multihomogeneous polynomial $f_Y\in K[U_0,\ldots,U_d]$, which is homogeneous of degree $\deg Y$ in each group of variables~$U_i$.
Furthermore, $f_Y$ is uniquely determined by $Y$ up to a non-zero scalar factor in~$\Qbar$, and each such polynomial is called a \emph{Chow form} for~$Y$.
Moreover, $Y$ is irreducible over~$K$ if and only if $f_Y$ is irreducible over~$K$.
If $Y=\overline{X}$ is the projective closure of an affine variety $X\subseteq \AA^n$, then we also write $f_X := f_Y$.

For more background on Chow forms, a rich interpretation thereof as general resultants, and relations to elimination theory, we refer to Gelfand, Kapranov, Zelevinsky~\cite{GelfandKapranovZelvinski94greenBible}.

\paragraph{Mahler measure.}
Let $f\in\CC[X_0,\ldots,X_d]$ be a polynomial in $d+1$ groups of $n+1$ variables each, having degree at most $D$ in each group of variables.
Let $S_{n}:=S^{2n-1}$ be the unit sphere in $\CC^{n}$, and let $\mu_n$ be the $U(n)$-invariant probability measure on~$S_n$.
Let $S_{n+1}^{d+1} = (S_{n+1})^{d+1}$ be the cartesian product of $d+1$ $(2n+1)$-dimensional unit spheres.
Philippon~\cite{philippon1991hauteursAlternativesI} defined the $S_{n+1}^{d+1}$-Mahler measure of $f$ as
\[
m(f;S_{n+1}^{d+1}) := \int_{S_{n+1}^{d+1}} \log |f|\ d\mu_{n+1}^{d+1},
\]
where $|~.~|$ here denotes the usual norm of~$\CC$, and $\mu_{n+1}^{d+1}:=(\mu_{n+1})^{\otimes (d+1)}$ is the Haar measure on~$S^{d+1}_{n+1}$.  
The (usual, logarithmic) Mahler measure of $f$ as defined by Mahler~\cite{Mahler62someInequalitiesForPolynomialsInSeveralVariables} is then the $S_1^{(n+1)(d+1)}$-Mahler measure of $f$, when considering each variable as a single group:
\[
m(f) := m\big(f;S_1^{(n+1)(d+1)}\big) 
= \int_{[0,1]^{(n+1)(d+1)}} \log\big|f\big(e^{2\pi t_1},\ldots,e^{2\pi t_{(n+1)(d+1)}}\big)\big|\ dt_1\cdots dt_{(n+1)(d+1)}.
\]
There are inequalities by Lelong~\cite{Lelong94mesureDeMahler},
\begin{equation}
\label{eqDifferenceBetweenDifferentMahlerMeasures}
0 \leq m(f) - m(f;S_{n+1}^{d+1}) \leq D(d+1)\sum_{i=1}^n\tfrac{1}{2i}.
\end{equation}
Philippon~\cite[Lem. 1.13]{philippon1986criteresIndependanceAlgebrique} and Krick, Pardo, Sombra~\cite[(1.1)]{KrickPardoSombra01sharpEstimatesForArithmeticNullstellensatz} showed that for 
$v\in M_K^\infty$,
\begin{equation}
\label{eqDifferenceBetweenMahlerMeasureAndLocalHeight}
|m(f^{\sigma_v}) - \log |f|_v| \leq D(d+1)\log(n+2).
\end{equation}

\paragraph{Height of a variety.}
A possible definition of a height of a projective variety would be the height of one of its Chow forms that is normalized to have at least one coefficient equal to~$1$.
It turns out that there are slightly more technical definitions that are easier to work with in practice.
We follow the definition from~\cite{KrickPardoSombra01sharpEstimatesForArithmeticNullstellensatz}.

Let $X\subseteq \PP^n$ be an irreducible $d$-dimensional projective variety defined over~$K$.
The height of $X\subseteq \PP^n$ is defined as
\[
h(X) := \tfrac{1}{[K:\QQ]}\bigg(\sum_{v\in M_K^\infty} N_v m(f_X^{\sigma_v};S_{n+1}^{d+1})
+ \sum_{v\not\in M_K^\infty} N_v \log |f_X|_v \bigg) + (d+1)\deg X\sum_{i=1}^n\tfrac{1}{2i}.
\]
For a reducible projective subvariety of $\PP^n$ the height is defined as the sum of the heights of its components.

\begin{remark}
The definition of $h(X)$ does not depend on the choice of~$K$:
We only need to check this for an irreducible affine variety $X\subset \AA^n$ over~$K$ and a finite field extension $L/K$.
If $X$ is still irreducible over $L$, then the height $h(X)$ does not depend on whether one takes $K$ or $L$ as the base field, the argument being the same as the one for the well-definedness of the height of a polynomial.
Let $X_1,\ldots,X_k$ be the irreducible components of $X$ over~$L$.
Then $f_X = \prod_i f_{X_i}$.
The Mahler measure behaves multiplicatively, $m(f_X^{\sigma_v};S_{n+1}^{d+1}) = \sum_i m(f_{X_i}^{\sigma_v};S_{n+1}^{d+1})$, and analogously $\log |f_X|_v = \sum_i \log |f_{X_i}|_v$ by the Gauss lemma.
Finally $\deg X = \sum_i \deg X_i$.
Thus $h(X)$ is indeed the same over~$K$ and over~$L$.
\end{remark}

\subsection{Fundamental height inequalities}

\begin{lemma}
\label{lemHeightInequalityNormalization}
Let $a_1,\ldots,a_n\in K^*$.
Then $h(1,a_2/a_1,\ldots,a_n/a_1) \leq h(a_1,\ldots,a_n)$.
\end{lemma}

\begin{proof}
This follows from the product formula.
\end{proof}

The following elementary height inequalities are mostly taken from Krick, Pardo, and Sombra~\cite[Lem. 1.2]{KrickPardoSombra01sharpEstimatesForArithmeticNullstellensatz}.
In the lemma, $\delta_{v\mid\infty}$ denotes a Kronecker delta, i.e.\ $\delta_{v\mid\infty}=1$ if $v$ is an archimedian valuation, and $\delta_{v\mid\infty}=0$ if $v$ is non-archimedian. 

\begin{lemma}
\label{lemElementaryLocalHeightInequalities}
Let $K$ be a number field, $v\in M_K$, and $f_1,\ldots,f_m\in K[x_1,\ldots,x_n]$.
Then
\begin{enumerate}
\item 
$h_v(\sum_{i=1}^m f_i) \leq \max_{i} h_v(f_i) + \delta_{v\mid\infty}\log m$.
\item 
$h_v(\prod_{i=1}^m f_i) \leq \sum_{i=1}^m h_v(f_i) + \delta_{v\mid\infty}\log (n+1)\sum_{i=1}^{m-1}\deg f_i$.
\item 
Let $F\in K[y_1,\ldots,y_m]$, and let $D:=\max_i\deg f_i$. 
Then
\[
h_v(F(f_1,\ldots,f_m)) \leq h_v(F) + \deg F\, \max_i h_v(f_i) + \delta_{v\mid\infty} \deg F\, \big(\log(m+1)+D\log(n+1)\big).
\]
\item 
Let $f\in K[y_1,\ldots,y_n]$ and $P\in K^n$. 
Then
\[
h_v(f(P)) \leq h_v(f) + \deg f\, h_v(P) + \delta_{v\mid\infty} \deg f\, \log(n+1).
\]
\item
Let $P,Q\in K^n$ and $g\in K[x_1,\ldots,x_n]$.
Then
\[
\log|g(P)-g(Q)|_v \leq \log |P-Q|_v + h_v(g) + (\deg g-1)\, h_v(P,Q) + \delta_{v\mid\infty}\deg g\,\log (n+2).
\]
\end{enumerate}
\end{lemma}

\begin{proof}
The first three inequalities are taken from Krick, Pardo, and Sombra~\cite[Lem. 1.2]{KrickPardoSombra01sharpEstimatesForArithmeticNullstellensatz}.
The fourth inequality is a special case of the third one.

It remains to show the last inequality.
We assume that $P=(p_1,\ldots,p_n)$ and $Q=(q_1,\ldots,q_n)$ differ only in one coordinate, say the $k$'th one; the general case follows then via the triangle inequality with an additional error term $\delta_{v\mid\infty}\log n$.
Consider the specialization $g_1(x):=g(p_1,\ldots,p_{k-1},x,q_{k+1},\ldots,q_n)$ and write $g_1(x) = \sum_i c_i x^i$.
Then $g(P)-g(Q) = g_1(p_k)-g_1(q_k) = (p_k-q_k)g_2(p_k,q_k)$, where $g_2(p_k,q_k):=\sum_{i\geq 1} c_i\sum_{j=0}^{i-1} p_k^jq_k^{i-1-j}$.
Hence 
\[
\log|g(P)-g(Q)|_v\leq \log|P-Q|_v + h_v(g_2(p_k,q_k)).
\]
If we interpret $g_2(p_k,q_k)$ as an evaluation of a polynomial of degree $\deg g-1$ in $n+1$ variables at the point $(p_1,\ldots,p_k,q_k,\ldots,q_n)$, we obtain 
\[
h_v(g_2(p_k,q_k))\leq h_v(g) + (\deg g-1)h_v(P,Q) + \delta_{v\mid\infty}(\deg g-1)\log(n+2).
\]
Putting these bounds together yields the claimed inequality.
\end{proof}

They imply the following global height inequalities.

\begin{lemma}
\label{lemElementaryHeightInequalities}
Let $f_1,\ldots,f_m\in \Qbar[x_1,\ldots,x_n]$.
Then
\begin{enumerate}
\item 
$\max_i h(f_i) \leq h(f_1,\ldots,f_m) \leq \sum_{i=1}^m h(f_i)$.
\item
$h(\sum_{i=1}^m f_i) \leq h(f_1,\ldots f_m) + \log m \leq \sum_{i=1}^m h(f_i) + \log m$. 
\item
$h(\prod_{i=1}^m f_i) \leq \sum_{i=1}^m h(f_i) + \log(n+1)\sum_{i=1}^{m-1}\deg f_i$.
\item
Let $F\in\Qbar[y_1,\ldots,y_m]$, and let $D:=\max_i \deg f_i$.
Then
\[
h(F(f_1,\ldots,f_m)) \leq h(F) + \deg F\,\big(h(f_1,\ldots,f_m) + \log (m+1) + D\log(n+1)\big)
\]
\item 
Let $f\in \Qbar[x_1,\ldots,x_n]$ and $P\in\Qbar^n$.
Then
\[
h(f(P)) \leq h(f) + \deg f\,(h(P) +\log(n+1)).
\]
\end{enumerate}
\end{lemma}

The following is a fundamental Liouville inequality, bounding valuations from below via the height. 

\begin{lemma}[Liouville inequality]
\label{lemLiouville}
Let $K$ be a number field, $a\in K\wo\{0\}$, and $v$ a normalized valuation of~$K$.
Then
\[
h_v(a)\geq \log |a|_v \geq -\frac{[K:\QQ]}{N_v}h(a).
\]
\end{lemma}

As $N_v\geq 1$, some authors prefer to omit this factor for simplicity.

\begin{lemma}[{\cite[Lemma 1.3]{KrickPardoSombra01sharpEstimatesForArithmeticNullstellensatz}}]
\label{lemHeightInequalityKvsOK}
Let $A\subset K$ be a finite subset.
Then there exists $b\in\ZZ\wo\{0\}$ such that $bA\subset\OK$ and
\[
h(A) \leq h(\{b\}\cup bA) \leq [K:\QQ]h(A).
\]
\end{lemma}

\begin{lemma}[Height of a hypersurface]
\label{lemHeightOfHypersurfaceInequality}
Suppose $X = V(f)\subset\AA^n$ is a hypersurface defined by a non-zero square-free polynomial $f\in\Qbar[x_1,\ldots,x_n]$ of degree~$\delta$. Assume that $f$ is normalized in the sense that at least one of its coefficients equals~$1$. 
Then
\[
-2\delta\log(n+2) \leq h(X) - h(f) \leq 2\delta\log(n+2) + \delta\sum_{i=1}^{n-1}\sum_{j=1}^i\frac{1}{2j}. 
\]
\end{lemma}

\begin{proof}
Let $F$ denote the homogenization of~$f$.
For the sake of Mahler measures as defined above, we consider $F$ as a `multi-homogeneous' polynomial in one group of $n+1$ variables.
Using
\cite[Sect. 1.2.3, 1.2.4]{KrickPardoSombra01sharpEstimatesForArithmeticNullstellensatz},
\[
h(X) = \tfrac{1}{[K:\QQ]}\sum_{v\in M_K} N_v h_v(X),
\]
where $h_v(X):= h_v(F)$ if $v\ndivides\infty$, and 
\[
h_v(X) := m(\sigma_v F;S_{n+1}) + \delta\sum_{i=1}^{n-1}\sum_{j=1}^i\tfrac{1}{2j},\quad \textnormal{if $v\divides\infty$.}
\]
For $v\divides\infty$, \eqref{eqDifferenceBetweenMahlerMeasureAndLocalHeight} reads
\[
|m(\sigma_v F;S_{n+1})-\log|F|_v|\leq 2\delta\log(n+2).
\]
As $f$ is normalized, so is $F$, and hence $\log|F|_v = h_v(F)$ for all~$v$.
Putting everything together yields the assertion.
\end{proof}

\begin{lemma}[Height of a point]
\label{lemHeightOfAPointInequality}
Suppose $X\subset\AA^n$ is a zero-dimensional variety of degree one, with $X(\Qbar) = \{P\}$, $P\in\Qbar^n$.
Then
\[
h(P) \leq h(X) \leq h(P) + \tfrac{1}{2}\log (n+1).
\]
\end{lemma}

\begin{proof}
For general affine varieties $V$ which satisfy \cite[Assumption 1.5]{KrickPardoSombra01sharpEstimatesForArithmeticNullstellensatz},
we have
\[
h(V) = \tfrac{1}{[K:\QQ]}\sum_{v\in M_K} N_v h_v(V)
\]
with $h_v(V):= h_{v_0}(\sigma_v(V))$ if $v\divides v_0 \in M_\QQ$, where
\[
h_\infty(\sigma_v V) := m(\sigma_v f_V^h;S_{n+1}^{d+1}) + \deg V \sum_{i=1}^{n-1}\sum_{j=1}^i \tfrac{1}{2j}
\]
(here, $f_V^h$ denotes the homogenization of a normalized~$f_V$) and
\[
h_p(\sigma_v V):=h_p(\sigma_v f_V)\quad \textnormal{for $p$ a rational prime.}
\]
For the given $X$, their Assumption 1.5 is satisfied, and by~\cite[Prop. 4]{philippon1991hauteursAlternativesI} or~\cite[Sect. 1.2.3]{KrickPardoSombra01sharpEstimatesForArithmeticNullstellensatz} we obtain
\[
h_\infty(\sigma_v X) = \tfrac{1}{2}\log(1+|\sigma_v p_1|^2+\ldots +|\sigma_v p_n|^2)
\]
and
\[
h_p(\sigma_v X) = h_p(\sigma_v P)\quad \textnormal{for $p$ a rational prime.}
\]
Thus the lemma follows from
\[
h_\infty(\sigma_v P)\leq h_\infty(\sigma_v X)\leq h_\infty(\sigma_v P) + \tfrac{1}{2}\log(n+1),
\]
which translates to
\[
\max \{1,|x_i|\st i=1,\ldots,n\} \leq (1+\sum_i |x_i|^2)^{1/2} \leq (n+1)^{1/2}\max\{1,|x_i|\st i=1,\ldots,n\},
\]
which is the elementary comparison between the $\ell^2$ and $\ell^\infty$ norms in $\RR^{n+1}$.
\end{proof}

More generally, the height of a general affine variety and the height of its Chow form bound each other as follows.

We call a Chow form $f_X$ normalized if one of its coefficients is equal to~$1$.
By Lemma~\ref{lemHeightInequalityNormalization}, among all possible Chow forms $f$ of a given projective variety, the normalized ones minimize the height $h(f)$. (However they are in general not the only minimizers, e.g. one can scale by roots of unity.)

\begin{lemma}[Height of a variety]
\label{lemHeightOfXvsChowFormOfX}
Let $X\subseteq\AA^n$, be a $d$-dimensional affine variety of degree~$D$, and let $f_X$ be a normalized Chow form of the projective closure of~$X$.
Then
\[
-D(d+1)\log(n+2) \leq h(f_X) - h(X) \leq D(d+1)\Big(\log(n+2)+\sum_{i=1}^n\tfrac{1}{2i} \Big). 
\]
\end{lemma}

\begin{proof}
This follows from~\eqref{eqDifferenceBetweenDifferentMahlerMeasures} and~\eqref{eqDifferenceBetweenMahlerMeasureAndLocalHeight}.
\end{proof}

\begin{lemma}[Height of an affine image {\cite[Prop. 2.4]{KrickPardoSombra01sharpEstimatesForArithmeticNullstellensatz}}]
\label{lemHeightOfAffineImage}
Let $X\subseteq\AA^n$ be a $d$-dimensional affine variety of degree~$\delta$.
Let $\phi:\AA^n\to\AA^{n'}$ be an affine map.
Then
\[
h(\phi(X))\leq h(X) + (d+1)\delta (h(\phi)+8\log(n+n'+1)).
\]
\end{lemma}

\subsection{Geometric and arithmetic \bezout{} inequalities}

The degree $\deg X$ of an irreducible affine variety $X\subseteq \AA^n$ is the number of points in the intersection of $X(\Qbar)$ with a generic affine plane of complementary dimension.
For a general affine variety $X\subseteq\AA^n$ the degree $\deg X$ is defined as the sum of the degrees of its irreducible components.

\begin{theorem}[Geometric \bezout{} inequality]
\label{thmGeometricBezout}
Let $V_1,\ldots V_m\subset\AA^n$ be affine varieties.
Then 
\[
\deg(V_1\cap\ldots \cap V_m) \leq \deg V_1 \cdot \ldots \cdot \deg V_m.
\]
\end{theorem}

For a proof as well as more background, see e.g.\ Fulton~\cite[Ex. 8.4.6]{Fulton98IntersectionTheory}.

\begin{corollary}
\label{corDegreeBoundForVarietyXinTermsOfFis}
Let $X=V(f_1,\ldots,f_m)\subseteq\AA^n$, and $D=\max_i \deg f_i$.
Then
\[
\deg X \leq D^{\min(n,m)}.
\]
\end{corollary}

\begin{proof}
If $m=\min(n,m)$ then the assertion follows immediately from Theorem~\ref{thmGeometricBezout}.
If $n=\min(n,m)$, let $g_1,\ldots,g_n$ be $n$ generic $\Qbar$-linear combinations of $f_1,\ldots,f_m$.
Then $\deg g_i = D$ for each $1\leq i\leq n$, and $X=V(g_1,\ldots,g_n)$.
The assertion now follows from Theorem~\ref{thmGeometricBezout} applied to $V(g_1),\ldots,V(g_n)$.
\end{proof}

Similarly, we will need to bound the height of an intersection of two varieties by the heights and degrees of the intersecting varieties.
Such arithmetic \bezout{} inequalities were established by Bost, Gillet, Soul\'e~\cite{BostGilletSoule91analogueOfArithmeticBezout,BostGilletSoule94heightsOfProjVarietiesAndPosGreenForms}, see also Philippon~\cite[Ch. 6]{NesterenkoPhilippon01introductionToAlgebraicIndependenceTheory}.
We will use the following version by Krick, Pardo, Sombra~\cite{KrickPardoSombra01sharpEstimatesForArithmeticNullstellensatz}.

\begin{theorem}[Arithmetic \bezout{} inequality~\cite{KrickPardoSombra01sharpEstimatesForArithmeticNullstellensatz}]
\label{thmArithmeticBezout}
Let $X\subseteq\AA^n$ be an affine variety, and let $f_1,\ldots,f_m\in\Qbar[x_1,\ldots,x_n]$ with $D_i:=\deg f_i$ and $D_1\geq \ldots\geq D_m$.
Let $n_0:=\min(\dim X,m)$.
Then
\[
h(X\cap V(f_1,\ldots,f_m)) \leq 
\Big(h(X)+h(f_1,\ldots,f_m)\deg X\sum_{i=1}^{n_0}\tfrac{1}{D_i} + n_0\log(n+1)\deg X \Big)\prod_{i=1}^{n_0}D_i.
\]
\end{theorem}

\begin{corollary}[{\cite{KrickPardoSombra01sharpEstimatesForArithmeticNullstellensatz}}]
\label{corHeightOfXboundedInTermsOfHeightOfFis}
Let $f_1,\ldots,f_m\in\Qbar[x_1,\ldots,x_n]$ with $D_i:=\deg f_i$ and $D_1\geq \ldots \geq D_m$.
Let $n_0=\min(n,m)$.
Then
\[
h(V(f_1,\ldots,f_m)) \leq 
\Big(h(f_1,\ldots,f_m)\sum_{i=1}^{n_0}\tfrac{1}{D_i}+(n+n_0)\log(n+1) \Big)\prod_{i=1}^{n_0}D_i
\]
\end{corollary}

\subsection{Arithmetic Nullstellens\"atze}

Krick, Pardo and Sombra~\cite{KrickPardoSombra01sharpEstimatesForArithmeticNullstellensatz}
proved the following essentially sharp arithmetic version of Hilbert's Nullstellensatz for polynomials without common zeros.

\begin{theorem}[Arithmetic Nullstellensatz -- \bezout{} version~\cite{KrickPardoSombra01sharpEstimatesForArithmeticNullstellensatz}]
\label{thmArithmNullstellensatzBezout}
Let $K$ be a number field, and suppose that $f_1,\ldots,f_m\in \OK[x_1,\ldots,x_n]$ are polynomials such that $V(f_1,\ldots,f_m)=\emptyset$ (in~$\Qbar^n$).
Let $D := \max_i \deg f_i$.
Then there exist $a\in \OK\wo 0$ and $\lambda_1,\ldots,\lambda_m\in\OK[x_1,\ldots,x_n]$ such that
\begin{enumerate}
\item
$a = \lambda_1f_1+\ldots+\lambda_mf_m$,
\item
$\deg \lambda_i \leq 4nD^n$, and
\item
$h(a,\lambda_1,\ldots,\lambda_m)\leq 4n(n+1)D^n\Big(h(f_1,\ldots,f_m)+\log m + (n+7)\log(n+1)D\Big)$.
\end{enumerate}
\end{theorem}

The \bezout{} version of Hilbert's Nullstellensatz for polynomials without common zeros implies the general Nullstellensatz via the classical Rabinowitsch trick.
The same trick works for arithmetic Nullstellens\"atze, proving the following theorem.

\begin{theorem}[Arithmetic Nullstellensatz -- general version]
\label{thmArithmNullstellensatzGeneral}
Let $K$ be a number field, and suppose $f_1,\ldots,f_m,g\in \OK[x_1,\ldots,x_n]$ are polynomials such that $g$ vanishes on the common zeros of $f_1,\ldots,f_m$.
Let $D := \max\{\deg g+1,\deg f_i \st i=1,\ldots,m\}$.
Then there exist an $a\in\OK\wo\{0\}$ and $\lambda_1,\ldots,\lambda_m\in\OK[x_1,\ldots,x_n]$ such that with $N:=4(n+1)D^{n+1}$,
\begin{enumerate}
\item
$ag^N = \lambda_1f_1 + \ldots + \lambda_mf_m$,
\item
$\deg \lambda_i \leq N(\deg g+1)$, and
\item
$h(a,ag^N,\lambda_1,\ldots,\lambda_m)\leq N(n+3)\Big(h(f_1,\ldots,f_m,g)+\log(m+1)+(n+8)\log(n+2)D\Big)$.
\end{enumerate}
\end{theorem}

\begin{proof}
Following Rabinowitsch's trick, note that $f_1,\ldots,f_m,1-x_{n+1}g\in \OK[x_1,\ldots,x_{n+1}]$ have no zeros in~$\AA^{n+1}$ in common.
Note that $1-x_{n+1}g$ has the same height and the same local heights as $g$, and its degree is larger by~$1$.

Using Theorem~\ref{thmArithmNullstellensatzBezout} we obtain polynomials $\mu_1,\ldots,\mu_{m+1}\in \OK[x_1,\ldots,x_{n+1}]$ of degree at most $N=4(n+1)D^{n+1}$ and an $a\in\OK\wo\{0\}$ such that
\[
a = \mu_1f_1 + \ldots + \mu_mf_m + \mu_{m+1}(1-x_{n+1}g).
\] 
If we specialize to $x_{n+1}:=1/g$ and clear denominators by multiplying the equation with~$g^N$, we obtain
\[
a g^N = \lambda_1f_1 + \ldots + \lambda_mf_m,
\]
where $\lambda_i:=g^N \mu_i|_{x_{n+1}=1/g} \in \OK[x_1,\ldots,x_n]$.

Let $h:=h(f_1,\ldots,f_m,g)$.
By Theorem~\ref{thmArithmNullstellensatzBezout}, we may assume
\[
h(a,g_1,\ldots,g_m) \leq H:= N(n+2)(h+\log(m+1) + (n+8)\log(n+2)D).
\]

To obtain the claimed height bounds, write $\mu_i = \sum_{\ell=0}^N \mu_{i\ell}x_{n+1}^\ell$ with $\mu_{i\ell}\in \OK[x_1,\ldots,x_n]$.
Then $\lambda_i = \sum_{\ell=0}^N \mu_{i\ell}g^{N-\ell}$.
For the summands,
\begin{align*}
h(\mu_{i\ell}g^{N-\ell}) 
& \leq H + (N-\ell)h(g) + \log(n+2)(N + (N-\ell)\deg g) \\
& \leq H + N(h(g) + \log(n+2)(\deg g + 1)).
\end{align*}
Thus, 
\[
h(\lambda_i) \leq H + N(h(g) + \log(n+2)(\deg g + 1)) + \log N.
\]
The same bound holds for $h(a g^N)$ and $h(a)$.
These height arguments work in the same lines locally via Lemma~\ref{lemElementaryLocalHeightInequalities}, and they hold at any $v\in M_K$ with equal bounds for $a,ag^N,\lambda_1,\ldots,\lambda_m$. 
Hence
\begin{align*}
h(a,ag^N,\lambda_1,\ldots,\lambda_m)
~\leq\ & H + N(h(g) + \log(n+2)(\deg g + 1)) + \log N \\
 \leq\ & N\Big((n+2)\big(h+\log(m+1)+(n+8)\log(n+2)D\big) + \\ 
       & ~~~~~~~~~~~~~~~~~~~~~~~~~~h(g) + \log(n+2)(\deg g + 1) + 1\Big)\\
 \leq\ & N(n+3)\Big(h+\log(m+1)+(n+8)\log(n+2)D\Big).
\end{align*}
This finishes the proof.
\end{proof}

\subsection{Degree bounds}

The following lemma gives a degree bound for the field extension obtained by augmenting the coordinates of the $\Qbar$-rational points of a zero-dimensional variety to the base field.

\begin{lemma}
\label{lemDegreeOfPointsOf0dimVariety}
Let $Z/K$ be a zero-dimensional affine variety of degree~$\delta$ defined over a number field~$K$.
Then there exists a field extension $L/K$ of degree $[L:K]\leq \delta!$, such that $Z(\Qbar) = Z(L)$, i.e.\ the closed geometric points of $Z$ have coordinates in~$L$.
\end{lemma}

\begin{proof}
The Galois group $G=\Gal(\Qbar/K)$ acts on $Z$. Let $H\leq G$ be the closed and normal subgroup that fixes the finite set $Z(\Qbar)$. Let $L = \Qbar^H$ be the fixed field associated to~$H$.
Then $Z(\Qbar) = Z(L)$.
Now $G/H$ can be identified with a subgroup of the group of permutations of the set $Z(\Qbar)$, and $|Z(\Qbar)|\leq \delta$.
Therefore, $[L:K] = |G/H|$ divides $\delta!$.
\end{proof}

The following lemma gives a degree bound for the field extension obtained by augmenting the coordinates of a single $\Qbar$-rational point of a zero-dimensional variety to the base field.

\begin{lemma}
\label{lemDegreeOfPointOf0dimVariety}
Let $Z/K$ be a zero-dimensional affine variety of degree~$\delta$ defined over a number field~$K$.
Let $P\in Z(\Qbar)$.
Then $K(P) = K(p_1,\ldots,p_n)$ has degree at most $\delta$ over~$K$.
\end{lemma}

\begin{proof}
A point $P\in Z(\Qbar)$ can be regarded as a point $x\in Z$ together with a $K$-algebra homomorphism $\phi: k(x)\to\Qbar$, see Liu~\cite[Prop. 3.2.18]{Liu02AlgebraicGeometryAndArithmeticCurves}, where $k(x)$ is the residue field of~$x$.
In this correspondence, $K(P)=\im\phi$.
Now, the degree $[K(P):K]=[k(x):K]$ equals the separable degree $[k(x):K]_{\textnormal{sep}}$ (as $\Char K=0$), which equals the number of $K$-algebra homomorphisms $k(x)\to\Qbar$.
The latter ones determine mutually distinct points of $Z(\Qbar)$.
Therefore, $[K(P):K]\leq |Z(\Qbar)|\leq \delta$.
\end{proof}

\subsection{Containment in coordinate hyperplanes}

The following lemma states that affine varieties cannot be contained in coordinate hyperplanes of sufficiently larger height.
This will be used iteratively in the proof of the main theorem.

\begin{lemma}[Containment in coordinate hyperplane]
\label{lemXcontainedInCoordinateHyperplane}
Let $X\subseteq \AA^n$ be a non-empty affine variety over~$\Qbar$, whose projective closure has Chow form~$f_X$.
Let $U=\{x_i=a\}$ be a coordinate hyperplane for some $i\in\{1,\ldots,n\}$ and some $a\in\Qbar$.
If 
\[
h(a) \geq h(f_X) + \log 2
\]
then
\[
X\not \subseteq U.
\]
\end{lemma}

\begin{proof}
We may assume $i=1$.
Let $K$ be a number field that contains $a$ and over which $X$ is defined.
Let $d=\dim X$.
Let $\overline X\subseteq \PP^n$ be the projective closure of~$X$, and similarly $\overline U$ the projective closure of~$U$.
Assume for a contradiction that $X\subseteq U$.
Then also $\overline X\subseteq \overline U$.
The Chow form $f_X$ is a multihomogeneous polynomial in $K[U_0,\ldots,U_d]$ in $d+1$ groups of $n+1$ variables $U_j = (u_{j0},\ldots,u_{jn})$.
Let $\alpha$ be a new variable and let $g(\alpha)\in K[U_0,\ldots,U_d][\alpha]$ be the polynomial obtained from specializing $U_0$ to the projective coordinates of the hyperplane $\{x_1 = \alpha\}$, that is,
\[
g(\alpha) := f_X|_{U_0=(-\alpha,1,\ldots,0)}.
\]
We claim that $a$ is a root of $g(\alpha)$, in the sense that $g(a) = f_X|_{U_0=(-a,1,0,\ldots,0)}$ is the zero polynomial:
Recall that $f_X$ takes as parameters the projective coordinates of $d+1$ hyperplanes in $\PP^n$, at which $f_X$ vanishes if and only if the intersection of $\overline X$ and all these hyperplanes is non-empty.
This intersection is at least zero-dimensional if one of the hyperplanes contains~$\overline X$, which follows from a simple dimension count (intersecting a projective variety with a projective hyperplane decreases its dimension by at most one).
Since $\overline X\subseteq \overline U$, this proves the claim $g(a)=0$.

Write
\[
g(\alpha) = \sum_{j=0}^D g_j\alpha^j,\qquad g_j\in K[U_0,\ldots,U_d],\qquad D=\deg_\alpha g(\alpha).
\]

We claim that $D\geq 1$:
If $D\leq 0$, then $g(\alpha)$ is constant, hence zero (as it has a root~$a$).
This means that any projective coordinate hyperplane $\overline{\{x_1=\alpha\}}$ for $\alpha\in \Qbar$ needs to contain $\overline X$.
As $\overline{\{x_1 = a\}}\subseteq \PP^n$ already contains $\overline X$, this means that $\overline X$ lies in the hyperplane at infinity, i.e. $X = \overline X\cap \AA^n = \emptyset$, which was excluded by assumption.
This proves the claim $D\geq 1$. 

Let $m$ be a non-zero monomial of $g_D$, and let $c_j':=g_j[m]\in K$ be its coefficient in~$g_j$.
Thus $c_D'\neq 0$, and $\sum_{j=0}^D c_j' a^j = 0$.
Let $c_j:=c_j' / c_D' \in K$.
Thus $c_D = 1$ and $a^D = -\sum_{j=0}^{D-1} c_j a^j$.

Next we claim that
\[
h(f_X) \geq h(c_0',\ldots,c_D') \geq |c_0',\ldots,c_D'|_v = |c_0,\ldots,c_D|_v = h(c_0,\ldots,c_D).
\]
The first inequality holds because $\{c_0',\dots,c_D'\}$ is a subset of the set of coefficients of $f_X$:
We replaced $u_{01}$ by the new variable $\alpha$, and some other variables by $0$ and $1$, thus only deleting a few monomials.
The second inequality is trivial.
The next equality holds by the product formula $\prod_{v\in M_K} |x|_v^{N_v} = 1$.
The last equality holds since $c_D=1$.
This proves the claim.
 
Thus in order to prove the lemma, it suffices to bound $h(c_0,\ldots,c_D) > h(a)-\log 2$.
To do this, we bound each local height $h_v$ separately.

For finite~$v$, let $A = |a|_v$.
We obtain $A^D \leq \max_{j=0}^{D-1} |c_j|_v A^j$. In case $A > 1$, this implies $A\leq |c_j|_v$ for some $j=0,\ldots,D-1$.
Hence $h_v(a)\leq h_v(\{c_0,\ldots,c_D\})$.

For infinite~$v$, let $A = |a|_v$ and $C = \max_{j=0}^{D-1} |c_j|_v$.
We obtain $A^D \leq \sum_{i=0}^{D-1} |c_j|_v A^j$.
If $A>1$, it follows $A^D\leq C\sum_{i=0}^{D-1} A^j = C\frac{A^D}{A-1}$.
Thus, $C\geq \frac{A^{D+1}-A^D}{A^D-1} = \frac{A^{D+1}-1}{A^D-1} - 1 \geq A-1$.
Hence $\log_+ C\geq \log_+ A-\log 2$, i.e. $h_v(a)\leq h_v(\{c_0,\ldots,c_D\})+\log 2$.

Summing these inequalities up for all $v\in M_K$, we obtain $h(a) < h(c_0,\ldots,c_D) + \log 2$, which remained to show.
\end{proof}

\begin{remark}[Alternative proof for Lemma~\ref{lemXcontainedInCoordinateHyperplane} (sketch)]
A less computational proof for a weaker statement can be obtained using the arithmetic \bezout{} inequality as follows.
Let $H$ be an affine plane of codimension $\dim X$ and small height such that $X\cap H$ is non-empty and zero-dimensional.
By the arithmetic \bezout{} theorem, $X\cap H$ has bounded height.
On the other hand, if $X\subseteq U$ then $X\cap H$ contains a point of height at least $h(a)$, which yields a contradiction if the height of $a$ is large.
Therefore $X\not\subseteq U$.
\end{remark}

A variety is called equidimensional if each its irreducible components is of the same dimension.

\begin{corollary}[Containment in coordinate hyperplane]
\label{corContainmentInCoordinateHyperplane}
Let $X\subseteq \AA^n$, $n\geq 1$, be a non-empty affine variety over $\Qbar$.
Let $U=\{x_i=a\}$ be a coordinate hyperplane for some $a\in\Qbar$.
If 
\[
h(a) \geq h(X) + 2\deg X(\dim X+1)\log(n+2) + \log 2
\]
then
\[
X\not\subseteq U.
\]
If furthermore $X$ is equidimensional, then so is $X\cap U$, and if moreover $X \cap U\neq\emptyset$ then $\dim (X\cap U) = \dim X - 1$.
\end{corollary}

\begin{proof}
The first part follows immediately from Lemmas~\ref{lemXcontainedInCoordinateHyperplane} and~\ref{lemHeightOfXvsChowFormOfX}, as well as the trivial estimate $\sum_{i=1}^n\tfrac{1}{2i}\leq \tfrac{1}{2}+\tfrac{1}{2}\log(n)\leq \log(n+2)$.

To prove the second statement about dimensions, suppose that $X$ is equidimensional and let $X'$ be a geometrically irreducible component of~$X$.
As the height of $X$ equals the heights of its components, $h(X')\leq h(X)$.
Thus also $X'\not\subseteq U$, and hence any component of $X\cap U$ has dimension at most $\dim X-1$.
That they have dimension at least $\dim X-1$ follows from the affine dimension theorem, see e.g.~\cite[Prop. I.7.1]{Hartshorne77algebraicGeometry}, which states that if $X',U\subseteq\AA^n$ are two irreducible affine varieties then $\codim W\leq \codim X' + \codim U$ for any component~$W$ of~$X'\cap U$.
\end{proof}

\section{Effective algebraic \lojasiewicz{} inequalities}
\label{secLojasieviczInequalities}

Let $X=V(f_1,\ldots,f_m) \subseteq \Qbar^n$ be a variety and let $P\in\Qbar^n$ be a point.
We are interested in statements of the following form:
If for each $i=1,\ldots,m$, $f_i(P)$ is close to zero in a certain sense, then $P$ must be close to $X$ in a similar sense.
In case $X(\Qbar)$ is empty, this should mean that 
the $f_i(P)$ cannot all be close to zero.
Statements of this type are called (algebraic) \emph{\lojasiewicz{} inequalities}.

Originally, \lojasiewicz{}~\cite[\S{}17]{Lojasiewicz59problemeDeLaDivision},~\cite{Lojasiewicz71ensemblesSemiAnalytiques}, \cite{Lojasiewicz58divisionDuneDistribution}
proved that for a real analytic function $f:U\to \RR$ for some open set $U\subseteq\RR^n$ and for any compact set $A\subset U$ there exist positive constants $\alpha$ and $C$ such that $|f(p)|^\alpha\geq C\,\dist(p,Z_f)$ for $p\in A$, where $Z_f = f^{-1}(0)$.
It was the main ingredient of \lojasiewicz's proof of Schwarz's division conjecture. 

Brownawell's Theorem A in~\cite{Brownawell88localDiophantineNullstellenInequalities} can be interpreted as an effective algebraic \lojasiewicz{} inequality for the complex points of varieties that are defined over~$\QQ$, where the distance is measured with respect to the standard norm.
With his theorem he extends an algebraic \lojasiewicz{} inequality for the empty variety that appeared implicitly in Masser--W\"ustholz~\cite{MasserWustholz83fieldsOfLargeTranscendenceDegree} (see Brownawell's introduction). Moreover, Brownawell~\cite[Thm. A']{Brownawell88localDiophantineNullstellenInequalities} gives an effective generalization for varieties over number fields $K=\QQ[\zeta]$, where his explicit bounds involve a notion of (local) height that depend on the chosen primitive element~$\zeta$.

Before that, Brownawell~\cite[Prop 8, 8']{Brownawell87boundsForDegreeInNullstellensatz} already proved a \lojasiewicz{} inequality for the empty variety over $\CC$, which is dependent on a non-explicit constant.

Ji, \kollar{}, and Shiffman~\cite{JiKollarShiffman92LojasiewiczInequality} proved a version of Brownawell's algebraic \lojasiewicz{} inequality for general algebraically closed fields and absolute values $v$ while at the same time improving Brownawell's exponents, however at the expense of loosing the effectiveness of the involved constant.
Let $P$ and $X$ be as above, and let $v\in M_K$ be a normalized valuation of~$K$.
Define the $v$-distance of $P$ to $X$ as $\dist_v(P,X):=\inf_{Q\in X(\Qbar)}|P-Q|_v$ (in~\cite{JiKollarShiffman92LojasiewiczInequality} an $\ell^2$-norm is used instead).
Ji, \kollar{}, and Shiffman derived an inequality of the form
\[
\dist_v(P,X)^M  \leq C\cdot \max_i{|f_i(x)|} \cdot (1+|P|_v)^{\overline B}
\]
for some $1\leq M\leq \overline B$, where $\overline B$ is a positive integer that depends explicitly on $n$ and $\deg f_1$, $\ldots$, $\deg f_m$, and where $C$ is a non-explicit constant that depends on $f_1,\ldots,f_m$.
Their focus was to optimize the exponent $M$ in terms of the degrees of the $f_i$.
We could immediately use this inequality for our purpose, if the constant $C$ was effectively computable. Hickel~\cite[Thm. 2.1.(iii)]{Hickel01solutionDuneConjectureDeBerensteinYger}
provided a very similar \lojasiewicz{} inequality, which again depends on a non-explicit multiplicative constant.
\kollar~\cite[Thm. 7.6]{Kollar99effectiveNullstellensatzForArbitraryIdeals}
proved a generalized \lojasiewicz{} inequality for intersecting varieties of the form
\[
\dist_v(P,X_1\cap\ldots\cap X_m) \ll \max_i \dist_v(P,X_i).
\]
His estimate involves a non-explicit multiplicative constant.
Cygan~\cite{Cygan98intersectionTheoryAndSperationExponent} and
Cygan, Krasi\'{n}ski, and Tworzewski~\cite{CyganKrasinskiTworzewski99separationOfAlgebraicSetsAndLojasiewiczExponent}
proved separation theorems in for complex varieties in the spirit of \kollar's result, involving non-explicit multiplicative constants.

In Section~\ref{secLojasieviczForNonEmptyVarieties} we give an effective \lojasiewicz{} inequality for non-empty varieties, which together with its proof arose from the attempt to make the proof of Ji, \kollar{}, and Shiffman effective.
Furthermore, the special cases when $V$ is empty or zero-dimensional are important.
In these cases, there are similar yet simpler effective proofs, which are given in Sections~\ref{secLojasiewiczForEmptyVariety} and~\ref{secLojasiewiczForZeroDimensionalVarieties}.

\subsection{\ldots\ for the empty variety} 
\label{secLojasiewiczForEmptyVariety}

\begin{theorem}[Effective algebraic \lojasiewicz{} inequality for $\emptyset$]
\label{thmLojasiewiczForEmptySet}
Let $K$ be a number field.
Let $f_1,\ldots,f_m\in K[x_1,\ldots,x_n]$ be polynomials without common zeros in $\Qbar^n$, i.e. $X=V(f_1,\ldots,f_m)\subseteq \AA^n$ is the empty variety.
Let $v$ be a normalized absolute value of~$K$.
Suppose $P\in \AA^n(\Qbar)$ lies in the ball $|P|_v\leq R$ 
for some $R\geq 1$.
Let $D:=\max_i \deg f_i$.
Then for some $i\in\{1,\ldots,m\}$,
\mathsmall{
\begin{equation}
\label{eqLojasiewiczForEmptySet_boundFor_fiPv}
\log |f_i(P)|_v \geq -4\tfrac{[K:\QQ]}{N_v}(n+1)^2D^n
\bigg([K:\QQ]h(f_1,\ldots,f_m) + \log m + (n+7)\log(n+1)D + \tfrac{\log R}{n+1}   \bigg)
\end{equation}
}
\end{theorem}

\begin{remark}
In case the coefficients of $f_1,\ldots,f_m$ lie in $\OK$, then the inner factor $[K:\QQ]$ can be omitted in Theorem~\ref{thmLojasiewiczForEmptySet}.
\end{remark}

\begin{proof}[Proof of Theorem~\ref{thmLojasiewiczForEmptySet}]
Let $H:=[K:\QQ]h(f_1,\ldots,f_m)$.
By Lemma~\ref{lemHeightInequalityKvsOK}, we choose $b\in\ZZ\wo\{0\}$ such that all $bf_i$ have coefficients in $\OK$ and $h(bf_1,\ldots,bf_m)\leq H$.
By the arithmetic Nullstellensatz~\ref{thmArithmNullstellensatzBezout}, there exists $a\in \OK\wo\{0\}$ and $\lambda_1,\ldots,\lambda_m\in\OK[x_1,\ldots,x_n]$ such that
\[
a = \sum_i \lambda_i bf_i, 
\]
$\deg \lambda_i \leq 4nD^n$, and $h(a,\lambda_1,\ldots,\lambda_m) \leq 4n(n+1)D^n(H+\log m + (n+7)\log(n+1)D) =: \wt H$.
By Lemma~\ref{lemHeightInequalityNormalization} it follows $h(1,\lambda_1/a,\ldots,\lambda_m/a)\leq \wt H$.

\medskip

If $v$ is infinite, $|a|_v = |\sum_i b\lambda_i(P) f_i(P)|_v \leq m\max_i |b\lambda_i(P)f_i(P)|_v = m|b\lambda_{i_0}(P)f_{i_0}(P)|_v$ for some~$i_0$.
Hence
\[
|f_{i_0}(P)|_v \geq 1/(m|b|_v \cdot|\lambda_{i_0}(P)/a|_v).
\]
By the definition of the height, $\log|b|_v\leq \tfrac{[K:\QQ]}{N_v}h(b)\leq \tfrac{[K:\QQ]}{N_v}H$.
Similarly, using $\max_i |p_i|_v\leq R$ and that $\lambda_{i_0}$ has at most $(n+1)^{\deg \lambda_{i_0}}$ monomials, we obtain $\log|\lambda_{i_0}(P)/a|_v\leq h_v(\lambda_{i_0}/a) + \deg\lambda_{i_0}\log R + \deg\lambda_{i_0}\log(n+1) \leq \tfrac{[K:\QQ]}{N_v}\wt H + 4nD^n(\log R+\log(n+1))$.
Putting bounds together we obtain
\[
\log |f_{i_0}(P)|_v \geq - \tfrac{[K:\QQ]}{N_v}(H+\wt H) - 4nD^n(\log R+\log(n+1)) - \log m.
\]

If $v$ is finite, we obtain the same bound without the summands $\log m$ and $\log(n+1)$.
In both cases, the claimed bound follows from trivial estimates.
\end{proof}

\begin{remark}[Comparison with a Liouville estimate]
In case one has a height bound for~$P$, one can obtain possibly stronger bounds for $\max_i |f(P)|_v$ simply by using the Liouville inequality:
As $f_1,\ldots,f_m$ have no common zero, $f_i(P)\neq 0$ for some~$i$, and thus for this $i$ we get
\begin{equation}
\label{eqLojasiewiczForEmptySet_LiouvilleEstimateForComparison}
\log |f_i(P)|_v\geq -\tfrac{[K:\QQ]}{N_v} h(f(P)) \geq -\tfrac{[K:\QQ]}{N_v}\big(h(f)+\deg f\,(h(P)+\log(n+1))\big).
\end{equation}
The caveat is that the right hand side depends on the \emph{global} height of~$P$, and not only on $h_v(P)$ as in Theorem~\ref{thmLojasiewiczForEmptySet}. 

In applications this is a crucial difference:
Suppose the aim is to show that $X$ is non-empty.
Using Theorem~\ref{thmLojasiewiczForEmptySet}, one would find a point $P$ in some ball $v$-close to $X$ in the sense that all $|f_i(P)|_v$ are small, such that it contradicts~\eqref{eqLojasiewiczForEmptySet_boundFor_fiPv} and thus proves $\dim X\geq 0$.
This has a chance to work in practice.
On the other hand, trying to construct a contradiction to~\eqref{eqLojasiewiczForEmptySet_LiouvilleEstimateForComparison} by finding a point $P$ that is $v$-close to~$X$ is in general less promising, as one would need an additional control over $h(P)$.
If for some $P$ $v$-close to $X$, $h(P)$ is too large to contradict~\eqref{eqLojasiewiczForEmptySet_LiouvilleEstimateForComparison}, and if one tries to solve this by moving $P$ even $v$-closer to $X$, one may run in circles.
\end{remark}

\subsection{\ldots\ for zero-dimensional varieties}
\label{secLojasiewiczForZeroDimensionalVarieties}

This section can be skipped by the reader.
We give a conceptually simple proof of an effective \lojasiewicz{} inequality for zero-dimensional varieties, which is an important special case -- however the bound we obtain here is quite weak in comparison to the more general bound from the next section.

Let $K$ be a number field, and let $v$ be an absolute value of~$\Qbar$.
For $x\in\Qbar^n$, let $||x||_v := \max_i |x_i|_v$.
For two affine varieties $X,Y\subseteq \AA^n$, let $\dist_v(X,Y):=\inf \{||x-y||_v\st x\in X(\Qbar), y\in Y(\Qbar)\}$.

\begin{theorem}[Effective algebraic \lojasiewicz{} inequality for zero-dimensional varieties]
\label{thmLojasiewiczFor0dim}
Let $K$ be a number field.
Let $f_1,\ldots,f_m\in K[x_1,\ldots,x_n]$ define a zero-dimensional affine variety $X=V(f_1,\ldots,f_m)\subseteq \AA^n$.
Let $v$ be a normalized absolute value of~$\Qbar$.
Suppose $P\in \AA^n(\Qbar)$ lies in the ball $|P|_v\leq R$ 
for some $R\geq 1$.
Let $D:=\max_i \deg f_i$.
Let $H:=[K:\QQ]h(f_1,\ldots,f_m)$ and $N:=4(n+1)D^{n+1}$.
Then for some $i\in\{1,\ldots,m\}$,
\[
\log \dist_v(P,X) \leq \max(1,H)[K:\QQ]^2(n+3)^3(D^n+2)! + 2\log R + \log m + \tfrac{1}{N\deg X}\log|f_{i}(P)|_v.
\]
\end{theorem}

\begin{remark}
The factor $(D^n)!$ in the bound is unreasonably large.
\end{remark}

\begin{proof}
According to Lemma~\ref{lemHeightInequalityKvsOK} there exists $B\in\ZZ\wo\{0\}$ such that $Bf_1,\ldots,Bf_m$ lie in $\OK[x_1,\ldots,x_n]$ and $h(B,Bf_1,\ldots,Bf_m)\leq H$. 
In particular,
\[
\log |B|_v \leq H. 
\]

Write $X(\Qbar) =\{Q_1,\ldots,Q_\delta\}$, with $\delta=\deg X\leq D^n$.
Following Lemma~\ref{lemDegreeOfPointsOf0dimVariety}, the points $Q_1,\ldots,Q_\delta$ have coordinates in an extension field $L/K$ of degree at most $\delta!$ over~$K$.
Using Lemma~\ref{lemHeightOfAPointInequality},
\begin{equation}
\label{eqIneqHQivsHX}
h(Q_1,\ldots,Q_\delta) \leq \sum_{i=1}^\delta h(Q_i) \leq h(X).
\end{equation}
Write $P=(p_1,\ldots,p_n)$ and $Q_i = (q_{i1},\ldots,q_{in})$.
By Lemma~\ref{lemDegreeOfPointOf0dimVariety}, $[K(Q_i):K]\leq \delta$, and thus $[K(Q_i):\QQ]\leq \delta[K:\QQ]$. 
Hence, using Lemma~\ref{lemHeightInequalityKvsOK} we can write $q_{ij} = a_{ij}/b_{i}$ with $a_{ij}\in\OL$, $b_{i}\in\ZZ$, and
\begin{equation}
\label{eqIneq_Hbi_vs_HQi}
h(b_i)\leq h(a_{i1},\ldots,a_{in},b_i)\leq \delta[K:\QQ] h(Q_i).
\end{equation}
For $i=1,\ldots,\delta$, $j=1,\ldots,n$, let $g_{ij}(x):= b_{i}(x_j-q_{ij}) = b_{i}x_j - a_{ij}$.
For each $i=1,\ldots,\delta$, choose $j_i\in\{1,\ldots,n\}$ such that $||P-Q_i||_v = |p_j-q_{ij}|_v$.
Then
\[
\dist_v(P,X)^\delta \leq \prod_{i=1}^\delta |\tfrac{1}{b_i}g_{ij_i}(P)|_v = |g(P)|_v/|b|_v,
\]
where $g:=\prod_{i=1}^\delta g_{ij_i}\in\OL[x_1,\ldots,x_n]$ and $b:=\prod_{i=1}^\delta b_{i}\in\ZZ$.
By construction, $g$ vanishes on $X(\Qbar)$.
Thus the arithmetic Nullstellensatz~\ref{thmArithmNullstellensatzGeneral} yields a relation
\[
ag^N = \lambda_1 Bf_1 + \ldots + \lambda_m Bf_m
\]
with $N=4(n+1)D^{n+1}$, $\lambda_i\in \OL[x_1,\ldots,x_n]$, $\deg \lambda_i\leq N(\deg g+1) = N(\delta+1)$, and
\[
h(a,ag^N,\lambda_1,\ldots,\lambda_m) \leq N(n+3)\big(h(Bf_1,\ldots,Bf_m,g)+\log(m+1)+(n+8)\log(n+2)D\big) =: \wt H.
\]

Next we bound several valuations.
By Lemma~\ref{lemHeightInequalityNormalization}, $h(\lambda_i/a)\leq h(a,\lambda_i) \leq \wt H$ for any $i=1,\ldots,m$.
With $[L:\QQ]\leq \delta![K:\QQ]$, the Liouville inequality (Lemma~\ref{lemLiouville}) implies
\[
\log |\lambda_i/a]_v \leq \delta![K:\QQ]\wt H.
\]
By Corollary~\ref{corHeightOfXboundedInTermsOfHeightOfFis},
\begin{equation}
\label{eqIneqForhX}
h(X)\leq D^n(nh(f_1,\ldots,f_m) + 2n\log(n+1))\leq nD^n(H/[K:Q]+2\log(n+1)).
\end{equation}
With~\eqref{eqIneq_Hbi_vs_HQi}, \eqref{eqIneqHQivsHX} and~\eqref{eqIneqForhX},
\begin{align*}
\log|b^{-1}|_v & \leq h(b)\leq \sum_i h(b_i) + \log \delta \leq \delta [K:\QQ] \sum_i h(Q_i)+\log \delta \\
& \leq \delta[K:\QQ] h(X)+\log \delta. \\
& \leq n\delta D^{n}(H+2[K:\QQ]\log(n+1)) + \log \delta.
\end{align*}
Furthermore, by~\eqref{eqIneq_Hbi_vs_HQi}, 
\[
h(g_i)\leq \delta[K:\QQ]h(Q_i),
\]
and thus with~\eqref{eqIneqHQivsHX} and~\eqref{eqIneqForhX},
\[
h(g)\leq \delta^2[K:\QQ]h(X)\leq n\delta^2 D^{n}(H+2[K:\QQ]\log(n+1)).
\]
With this we can bound
\[
h(Bf_1,\ldots,Bf_n,g)\leq H + h(g) \leq H+n\delta^2 D^{n}(H+2[K:\QQ]\log(n+1)).
\]
If $v$ is infinite, 
\begin{align*}
\log|\lambda_i(P)/a|_v 
& \leq h_v(\lambda_i/a) + \deg\lambda_i\log R + \deg\lambda_i\log(n+1) \\
& \leq [L:\QQ] h(\lambda_i/a) + N(\delta+1)(\log R + \log(n+1)).
\end{align*}
Next, 
\[
|g^N(P)|_v = \big|\sum_i \tfrac{B}{ab}\lambda_i(P)f_i(P)\big|_v \leq m\max_i \big|\tfrac{B}{ab}\lambda_i(P)f_i(P)\big|_v = m \big|\tfrac{B}{ab}\lambda_{i_0}(P)f_{i_0}(P)\big|_v,
\]
for some $i_0$.
Putting bounds together, we obtain
\begin{align*}
\log \dist_v(P,X)^{\delta N} 
& \leq \log |g^N(P)/b^N|_v \\
& \leq \log (m|\lambda_{i_0}(P)/a\cdot f_{i_0}(P)B/b^N|_v) \\
& \leq [L:\QQ]\wt H + N(\delta+1)(\log R+\log(n+1)) + H + \log |f_{i_0}(P)|_v + \\
& \ \ \ \ \ Nn\delta D^{n}(H+2[K:Q]\log(n+1))+ N\log \delta  + \log m
\end{align*}
Rather crude estimates yield
\[
\log \dist_v(P,X) \leq \max(1,H)[K:\QQ]^2(n+3)^3(D^n+2)! + 2\log R + \log m + \frac{1}{\delta N}\log|f_{i_0}(P)|_v.
\]
For finite $v$, all bounds work equally and in parts even stronger, without certain summands $\log(n+1)$ and $\log(m)$.
\end{proof}

\subsection{\ldots\ for non-empty varieties}
\label{secLojasieviczForNonEmptyVarieties}

Let $K$ be a number field, and let $v$ be an absolute value of~$\Qbar$.
For $x\in\Qbar^n$, let $||x||_v := \max_i |x_i|_v$.
For two affine varieties $X,Y\subseteq \AA^n$, let
\[
\dist_v(X,Y):=\inf \{||x-y||_v\st x\in X(\Qbar), y\in Y(\Qbar)\}.
\]

\begin{theorem}[Effective algebraic \lojasiewicz{} inequality for non-empty varieties]
\label{thmLojasiewiczForNonEmptyVarieties}
Let $K$ be a number field.
Let $f_1,\ldots,f_m\in K[x_1,\ldots,x_n]$ define a non-empty affine variety $X=V(f_1,\ldots,f_m)\subseteq \AA^n$. 
Let $v$ be a normalized absolute value of~$\Qbar$.
Suppose $P\in \AA^n(\Qbar)$ lies in the ball $|P|_v\leq R$ 
for some $R\geq 1$ such that $|f_i(P)|_v\leq 1$ for all $1\leq i\leq m$.
Let $D:=\max_i \deg f_i$.
Let $H:=h(f_1,\ldots,f_m)$.
Then
\mathsmall{
\begin{equation}
\label{eqLojasievicForNonEmptyVarieties_distanceInequality}
\log \dist_v(P,X) \leq \frac{\max_i \log|f_{i}(P)|_v}{4(n+1)(D^n+1)^{n+2}} + \tfrac{[K:\QQ]^2}{N_{v|K}}(n+7)^2(D^n+1)\big(H + \log(mnD^{2n}) + 21 \big)
 + 2\log R.
\end{equation}
}
\end{theorem}

\begin{remark}[General bound]
The assumption $|f_i(P)|_v\leq 1$ may seem artificial, and indeed it is.
The only purpose is to avoid a case distinction in the statement of the theorem and in the proof.
In case some $|f_i(P)|_v$ is larger than $1$, then the estimate~\eqref{eqLojasievicForNonEmptyVarieties_distanceInequality} will hold if the denominator under $\max_i \log|f_i(P)|_v$ is replaced by $4(n+1)D^{n+1}$, as this is a lower bound for the term $N\delta'$ from the proof.
\end{remark}

\begin{remark}[Algorithmic bound]
If one uses this theorem in a computer, one can obtain a significantly stronger upper bound for $\dist_v(P,X)$. To do this, one simply lets the algorithm follow the chain of inequalities from the proof of this theorem, and whenever possible one replaces estimates by explicit computations.
\end{remark}

We need the following technical lemma, which loosely speaking claims the existence of a  linear projection from $n$ to $1$ dimensions that does not contract a given finite set of vectors by more than an explicit factor.

\begin{lemma}
\label{lemTechnicalProjectionDirection}
Let $n\geq 1$.
Suppose we are given a finite set $U\subset\Qbar^n\wo\{0\}$ of non-zero vectors, as well as a normalized valuation of~$\Qbar$ denoted~$v$.
Put $R:=(n-1)|U|$.
If $v\divides\infty$, let $c:=1$, otherwise let $c:=1/R$.
Then there exists a vector $s\in\ZZ^n$ with $s_1 = 1$ and $\max_i |s_i|_\infty\leq R$, such that
\begin{equation}
\label{eqTechnicalProjectionDirectionInequality}
|s^tu|_v \geq c\max_{2\leq i\leq n} |u_i|_v \quad \textnormal{ for all } \quad u\in U. 
\end{equation}
\end{lemma}

\begin{proof}[Proof of Lemma~\ref{lemTechnicalProjectionDirection}]
For $n=1$ we can trivially choose $s=1$, hence from now on assume that $n\geq 2$.
If $v\divides\infty$, let $\eps = 2$, otherwise $\eps = 1$.

We will use a simple counting argument.
The discrete cube
\[
C=\{1\}\times \{R,R-\eps,R-2\eps,\ldots,R-R\eps\}^{n-1}
\]
parametrizes certain relevant vectors~$s\in C$ with $s_1=1$ and $\max_i |s_i|_\infty\leq R$.
Clearly, $|C| = (R+1)^{n-1}$.

Fix $2\leq i\leq n$ and $u\in U$.
Choose $2\leq k\leq n$ such that $|u_k|_v = \max_{2\leq \ell\leq n} |u_\ell|_v$.
Suppose that there are two distinct lattice points $s,s'\in C$ that differ only in the $k$'th coordinate and such that both satisfy $|s^tu|_v < c |u_i|_v$ and $|(s')^tu|_v < c |u_i|_v$.
Then by the triangle inequality,
\[
|(s_k-s'_k)|_v|u_k|_v  = |(s-s')^tu|_v \leq \eps\max(|s^tu|_v,|(s')^tu|_v)< \eps c |u_i|_v \leq \eps c |u_k|_v.
\]
Thus $|s_k-(s')_k| < \eps c$, which is impossible and hence contradicts the existence of $s$ and $s'$ as above.
This means that if we fix all but the $k$'th coordinate of $s$, there is at most one such $s\in C$ with $|s^tu|_v < c |u_i|_v$.
Therefore the number of $s\in C$ with $|s^tu|_v < c |u_i|_v$ is at most $(R+1)^{n-2}$.

Repeating this count for all $(n-1)|U|$ choices for $i$ and $u$, we see that at most $(n-1)|U|(R+1)^{n-2}$ points $s\in C$ do not satisfy~\eqref{eqTechnicalProjectionDirectionInequality}.
Since $|C| = (R+1)^{n-1}$ and $R+1>(n-1)|U|$, there is at least one $s\in C$ that satisfies~\eqref{eqTechnicalProjectionDirectionInequality}.
\end{proof}

\begin{proof}[Proof of Theorem~\ref{thmLojasiewiczForNonEmptyVarieties}]

Let $d:=\dim X$ and $\delta := \deg X$.
Let $\overline X$ denote the projective closure of $X$ in $\PP^n$.
Let $Z := \overline X \cap H_\infty$ denote the intersection of $\overline X$ with the hyperplane at infinity $H_\infty=(\PP^n\wo \AA^n)=\{x_0=0\}$.
By \bezout's theorem, Theorem~\ref{thmGeometricBezout}, $\deg Z\leq \delta\leq D^n$.
We regard $Z$ as a subvariety of $\PP^{n-1}$ and consider its Chow form $f_Z$, which is a multihomogeneous polynomial of degree $\deg Z$ in in each of the $d' = \dim_{\PP^{n-1}} Z + 1$ groups of variables.
Hence $\deg f_Z\leq n\delta$.
By the weak combinatorial Nullstellensatz, see Section~\ref{secHistory}, we can choose values among $\{0,\ldots,n\delta\}$ for each variable of $f_Z$, such that $f_Z$ is non-zero at this point.
Since $\dim Z < \dim X$, $d'\leq d$.
Thus there are $d'$ hyperplanes $h_i:\sum_{j=1}^n\mu_{ij} x_j = 0$ in $\PP^{n-1}$, $i=1,\ldots,d'$, with $\mu_{ij}\in \{0,\ldots,n\delta\}$, such that their intersection $U:=V(h_1,\ldots,h_{d'})$ is disjoint from $Z(\Qbar)$.
We can consider $h_1,\ldots,h_{d'}\in\ZZ[x_0,\ldots,x_n]$, with coefficient~$0$ for~$x_0$, and in this way, these hyperplanes and $U$ are projective subspaces of $\PP^n$ that go through the origin~$e_0=(1,0,\ldots,0)^t$.
Using $h_1,\ldots,h_{d'}$, we will construct a linear surjection $\phi_d:\AA^n\to \AA^d$ of bounded height, such that $\ker \phi_d$ is contained in~$U$.

We may assume that the $\QQ$-span of $h_1,\ldots,h_{d'}$ is spanned by $h_1,\ldots,h_{d''}$, with $d''\leq d'$.
Let $M_{d''}$ be the $d''\times n$-matrix $(\mu_{ij})_{1\leq i\leq d'',1\leq j\leq n}$, which is of full rank.
By relabeling the coordinates, we may assume that the left $d''\times d''$-submatrix is invertible.
Let $M_{n,s}$ be the $n\times n$-matrix obtained from $M_{d''}$ by appending rows $(e_{d''+1})^t,\ldots,(e_n)^t$, except that in row $d+1$ we put $(0,\ldots,0,s_{1},\ldots,s_{n-d})$, for a yet to be specified $s=(s_1,\ldots,s_{n-d})^t\in\ZZ^{n-d-1}$ with $s_1 = 1$.
\[
M_{n,s} =  
\begin{pmatrix}
\mu_{11}   & & & \dots & ~~~~~~~~\mu_{1n} \\
\vdots         & & &       & ~~~~~~~~\vdots \\
\mu_{d''1} & & & \dots & ~~~~~~~~\mu_{d''n} \\
\hline
& & I_{d-d''} &  \\
\hline
& & & 1 & s_{2}\cdots s_{n-d}\\
& & & & I_{n-d-1} \\
\end{pmatrix}
\]
Clearly, $M_{n,s}$ is invertible.
For $k=1,\ldots,n$, let $M_{k,s}$ denote the $k\times n$ submatrix of $M_{n,s}$ obtained from taking the first $k$ rows.
We also write $M_d:=M_{d,s}$, since it does not depend on~$s$.
Let $\phi_{n,s},\phi_{d+1,s},\phi_d$ denote the linear maps associated to the matrices $M_{n,s},M_{d+1,s},M_d$, respectively.
Furthermore, let $\pi_k$ denote the projection to the first $k$ coordinates, for any suitable domain and any suitable~$k$.
The following diagram commutes.
\[
\xymatrix{
\AA^n \ar[rr]^{\phi_{n,s}} \ar[drr]^{\phi_{d+1,s}} \ar[ddrr]_{\phi_d} & & \AA^n \ar[d]^{\pi_{d+1}} \\
& & \AA^{d+1} \ar[d]^{\pi_d} \\
& & \AA^d 
}
\]

We claim that $\phi_{d}(X) = \AA^n$: 
Assume that some point $Q\in\AA^n$ is not in the image $\phi_d(X)$.
Let $W:=\phi_d^{-1}(Q)$ be its preimage in $\AA^n$, which is an affine plane of dimension $n-d$ with $W\cap X=\emptyset$.
Since $X$ is non-empty, its degree $\delta$ is positive, and thus $\overline X\cap \overline W$ is non-empty.
This intersection must occur within $H_\infty$, hence $Z\cap \overline W\neq\emptyset$.
The equations for $W\cap H_\infty$ in $\PP^{n-1}=H_\infty$ are given by the polynomials $h_1,\ldots,h_{d''},x_{d''+1},\ldots,x_d$, and thus $W\cap H_\infty \subseteq U$.
Thus $Z\cap U\neq\emptyset$, which is a contradiction by how we constructed~$U$.
This proves the claim.

Moreover we see that the fibers of the map $\phi_d|X:X\to\AA^d$ are $0$-dimensional varieties:
We just showed that they are non-empty.
If they were of dimension larger than~$0$, then using their non-zero degree, their projective closure would intersect $H_\infty$.
These intersection points again witness a non-empty intersection of $Z$ and $U$, which is impossible.

Let $Q_{k,s}:=\phi_{k,s}(P)\in\AA^k(\Qbar)$ for $k\in\{d+1,n\}$, and $Q_{d}:=\phi_d(P)\in\AA^d(\Qbar)$.

Let $Y_{1},\ldots,Y_{\delta}\in X(\Qbar)$ denote the preimages of $Q_d$ under $\phi_{d}$ in~$X$ (counted with multiplicity).
Let $Z_{k,s,i}:= \phi_{k,s}(Y_{i})\in\AA^{k}(\Qbar)$.

Next we apply Lemma~\ref{lemTechnicalProjectionDirection} to the set of vectors $u_i\in\Qbar^{n-d}$ ($1\leq i\leq\delta$), where $u_i = (P-Y_i)_{d+1,\ldots,n}$ is $(n-d)$-vector obtained from the last $n-d$ coordinates of $P-Y_i$.
By Lemma~\ref{lemTechnicalProjectionDirection}, there exists an $s\in\ZZ^{n-d}$ with $s_1=1$ and $|s|_\infty\leq n\delta$ such that 
\begin{equation}
\label{eqDistanceQZn_vs_QZ_dplus1}
\\dist_v(Q_{n,s},Z_{n,s,i}) \leq n\delta\cdot \dist_v(Q_{d+1,s},Z_{d+1,s,i}) \quad \textnormal{for all} \quad 1\leq i\leq \delta.
\end{equation}
The initial factor $n\delta$ on the right-hand can be omitted if $v\divides\infty$.

The entries of $M_{n,s}$ lie in $[-n\delta,n\delta]\cap\ZZ$, hence $\det M_{n,s}\in\ZZ$ satisfies $|\det M_{n,s}|_v\leq n!(n\delta)^n$ (estimate is not optimal).
The adjugate matrix $\adj(M_{n,s})$ of $M_{n,s}$ also has entries in $\ZZ$ of $v$-norm bounded by $n!(n\delta)^n$.
Thus (no matter whether $v\divides\infty$ or not) the entries of $(M_{n,s})^{-1} = \tfrac{1}{\det M_{n,s}}\adj(M_{n,s})$ satisfy the same $v$-norm bound.
Thus we obtain
\begin{equation}
\label{eqDistancePY_vs_QZn}
\dist_v(P,Y_{i,s}) \leq n(n^2\delta)^n\cdot \dist_v(Q_{n,s},Z_{n,s,i}) \quad \textnormal{for all} \quad 1\leq i\leq \delta.
\end{equation}
The initial factor $n$ on the right-hand side comes from the triangle inequality and can be omitted if $v\ndivides\infty$.

\medskip

Let $X_s:=\phi_{d+1,s}(X)\subset\AA^{d+1}$. 
Then $X_s$ is a hypersurface, $X_s = V(g_s)$ for some $g_s\in K[x_1,\ldots,x_{d+1}]$.
Let $\delta':=\deg g_s = \deg X_s$, which satisfies $\delta'\leq \delta$. 

Let $\ell$ be the coefficient of $x_{d+1}^{\delta'}$ in~$g_s$.
We claim that $\ell\neq 0$:
To prove this claim, note that $\ell$ is the coefficient of $x_{d+1}^{\delta'}$ in $\wt g_0(x_d):=g_s(0,\ldots,0,x_{d+1})\in\Qbar[x_{d+1}]$.
That is, we need to show that $\deg \wt g_0 = \delta'$.
Geometrically, $\deg \wt g_0$ equals the number of intersections (counted with multiplicity) of $X_s$ with the $d+1$'st coordinate axis $W_{d+1}:=V(x_1,\ldots,x_d)\subset \AA^{d+1}$.
That is, we need to show that $\overline{X_s}$ does not intersect $\overline{W_{d+1}}$ at infinity.
Pulling this back via $\phi_{d+1,s}$, this is equivalent to the fact that $\overline{X}$ and $\overline{\ker \phi_d}$ do not meed at the hyperplane at infinity.
The latter is assured by the way we constructed~$U$.
Thus the claim is proved.

As $\ell\neq 0$, Liouville's inequality yields
\begin{equation}
\label{eqEllv_intermsof_hga}
\log |\ell|_v \geq -\tfrac{[K:\QQ]}{N_{v|K}}h(g_s).
\end{equation}
Let $\wt g(x_{d+1}) := g_s(Q_{i,d},x_{d+1})\in\Qbar[x_{d+1}|$. 
Let $\xi_i$ ($1\leq i\leq \delta'$) be the zeros of $\wt g$ (counted with multiplicity).
Then $\wt g(x_{d+1}) = \ell\prod_{i=1}^{\delta'} (x_{d+1}-\xi_i)$.
Each $\xi_i$ is the $(d+1)$'th coordinate of $Z_{d+1,s,j}$ for some $1\leq j\leq \delta$, and for simplicity we may assume $j=i$ (simply reorder the indices for that). 
Further note that $|(Q_{d+1,s})_i - \xi_i|_v = \dist_v(Q_{d+1,s},Z_{d+1,i,s})$.

Let $G_s:= g_s\circ\phi_{d+1,s} \in K[x_1,\ldots,x_n]$.
Clearly, $G_s$ vanishes on~$X$.
Furthermore,
\[
|G_s(P)|_v = |g_s(Q_{d+1,s})|_v = |\wt g((Q_{d+1,s})_i)|_v = |\ell|_v\prod_{i=1}^{\delta'}\dist_v(Q_{d+1,s},Z_{d+1,i,s}). 
\]
Using estimates~\eqref{eqDistanceQZn_vs_QZ_dplus1} and~\eqref{eqDistancePY_vs_QZn}, this implies
\begin{align}
\dist_v(P,X)^{\delta'}
&\ \leq\ \prod_{i=1}^{\delta'} \dist_v(P,Y_{i,s}) 
\ \leq\ \prod_{i=1}^{\delta'} (n^2\delta)^{n+1} \dist_v(Q_{d+1,s},Z_{d+1,i,s}) \nonumber\\
&\ \leq\  (n^2\delta)^{(n+1)\delta'}\frac{|G_s(P)|_v}{|\ell|_v}. \label{eqBoundDistvPX_intermsof_GaP}
\end{align}
It remains to bound $|G_s(P)|_v$ in terms of $\max_i |f_i(P)|_v$.
For this we make use of the arithmetic Nullstellensatz (Theorem~\ref{thmArithmNullstellensatzGeneral}).

First, according to Lemma~\ref{lemHeightInequalityKvsOK} there exists $B\in\ZZ\wo\{0\}$ such that $Bf_1,\ldots,Bf_m,BG_s\in\OK[x_1,\ldots,x_n]$ and 
\[
h(B,Bf_1,\ldots,Bf_m,BG_s)\leq [K:\QQ](H + h(G_s)) =: H_g.
\]
In particular,
\begin{equation}
\label{eqBoundForBv}
\log |B|_v \leq H_g. 
\end{equation}
Since $G_s$ vanishes on $X$, Theorem~\ref{thmArithmNullstellensatzGeneral} yields a relation
\[
a G_s^N = \lambda_1 B f_1 + \ldots \lambda_m B f_m,
\]
with $N=4(n+1)D_g^{n+1}$, $D_g = \max(D,\delta'+1)$, $\deg\lambda_i \leq N(\delta'+1)$, and
\[
h(a,aG_s^N,\lambda_1,\ldots,\lambda_m)\leq N(n+3)
\Big(H_g+\log(m+1)+(n+8)\log(n+2)D_g\Big) =: \wt H.
\]
Next we estimates heights.
Using Corollary~\ref{corHeightOfXboundedInTermsOfHeightOfFis},
\[
h(X)\leq nD^n(H+2\log(n+1)).
\]
Clearly, 
\[
h(\phi_{d+1,s})=h(M_{d+1,s})\leq h(1,2,\ldots,n\delta)\leq \log(n\delta).
\]
Hence, using Lemma~\ref{lemHeightOfAffineImage},
\[
h(X_s) \leq h(X) + (d+1)\delta(\log(n\delta) + 8\log(n+d+2)) \nonumber 
\]
Using Lemma~\ref{lemHeightOfHypersurfaceInequality},
\begin{align}
h(g_s)
&\ \leq\ h(X_s) + 2\delta'\log(d+3) \nonumber \\
&\ \leq\ nD^n(H + \log\delta + 11\log(n+1) + 8). \label{eqhga_intermsof_H}
\end{align}
Using Lemma~\ref{lemElementaryHeightInequalities}, 
\begin{align}
h(G_s) \ & \leq \ h(g_s) + (d+1)\delta\big(\log(n\delta) + \log((n+1)(d+1)+1) + \log(n+1)\big) \nonumber \\
& \leq \ h(g_s) + n\delta(\log\delta + 4\log(n+1)) \nonumber\\
& \leq \ nD^n(H + 2n\log D + 15\log(n+1) + 8). \label{eqhGa_intermsof_H}
\end{align}
This bounds also $H_g$,
\begin{equation}
\label{eqHg_intermsof_Hgstar}
H_g \leq [K:\QQ](nD^n+1)(H + 2n\log D + 15\log(n+1) + 8) =: H_g^*.
\end{equation}
Note that via~\eqref{eqhga_intermsof_H} we can also simply bound
\begin{equation}
\label{eqhga_intermsof_Hgstar}
h(g_s) \leq H_g^*.
\end{equation}

Next we bound certain valuations.
By Lemma~\ref{lemHeightInequalityNormalization}, $h(\lambda_i/a)\leq h(a,\lambda_i) \leq \wt H$ for any $i=1,\ldots,m$.
Liouville's inequality (Lemma~\ref{lemLiouville}) implies
\[
h_v(\lambda_i/a) \leq [K:\QQ]\wt H.
\]
Local height inequalities (Lemma~\ref{lemElementaryLocalHeightInequalities}) yield
\begin{align}
\log|\lambda_i(P)/a|_v 
& \leq h_v(\lambda_i/a) + \deg\lambda_i\log R + \deg\lambda_i\log(n+1) \nonumber \\
& \leq [K:\QQ]\wt H + N(\delta'+1)(\log R + \log(n+1)). \label{eqLambdaiPav_intermsof_wtH}
\end{align}
Next, 
\begin{equation}
\label{eqGaNPv_via_Nullstellensatz}
|G_s^N(P)|_v
= \big|\sum_i \tfrac{\lambda_i(P)}{a}Bf_i(P)\big|_v 
\leq m\max_i \big|\tfrac{\lambda_i(P)}{a}Bf_i(P)\big|_v 
= m \big|\tfrac{\lambda_{i_0}(P)}{a}Bf_{i_0}(P)\big|_v,
\end{equation}
for some $i_0$.
Putting bounds together, from~\eqref{eqBoundDistvPX_intermsof_GaP} using~\eqref{eqGaNPv_via_Nullstellensatz} we obtain
\mathsmall{
\begin{align*}
\log \dist_v(P,X)^{\delta' N} 
& \leq N(n+1)\delta'\log (n^2\delta)+\log |G_s^N(P)|_v-N\log|\ell|_v \\
& \leq N(n+1)\delta'\log (n^2\delta) + \log m + \log |\tfrac{\lambda_{i_0}(P)}{a}|_v + \log |B|_v + \log |f_{i_0}(P)|_v -N\log|\ell|_v.
\end{align*}
}%
Using~\eqref{eqBoundForBv} to bound $\log |B|_v$, \eqref{eqLambdaiPav_intermsof_wtH} to bound $\log |\lambda_{i_0}(P)/a|_v$, \eqref{eqEllv_intermsof_hga} and~\eqref{eqhga_intermsof_Hgstar} to bound $-\log |\ell|_v$, \eqref{eqHg_intermsof_Hgstar} to bound each $H_g$ by $H_g^\star$, 
as well as trivial estimates such $\delta'\leq\delta\leq D^n$ and $d+1\leq n$, we obtain
\mathsmall{
\begin{align*}
\log\dist_v(P,X)
& \, \leq \, \frac{\log|f_{i_0}(P)|_v}{N\delta'} + (n+1)\log(n^2D^n)+2\log R+2\log(n+1)\ + \\
&  \ \ \ \ \ \frac{1}{\delta'}\Big( H_g^* + \tfrac{1}{N}\log m +\tfrac{H_g^*}{N} + \tfrac{[K:\QQ]}{N_{v|K}}(n+3)\big(H_g^*+\log(m+1)+(n+8)\log(n+2)D_g\big)\Big).
\end{align*}
}%
The only non-positive summand on the right hand side is $\log|f_{i_0}(P)|_v$, thus we bound its denominator trivially from above, $N\delta'\leq 4(n+1)(D^n+1)^{n+2}$, and we bound the second denominator $\delta'$ trivially from below, $\tfrac{1}{\delta'}\leq 1$.
Plugging in the definition for $H_g^*$ and a few further trivial estimates for the sake of simplicity yields the asserted inequality.
\end{proof}

\begin{remark}[Remarks to the proof of Theorem~\ref{thmLojasiewiczForNonEmptyVarieties}]
The reader may wonder why we argued using a somewhat technically constructed $\wt g$, instead of simply projecting $X$ via $\phi_{d+1,s}$ to $\AA^{d+1}$, which becomes a hypersurface, for which we can relatively simply bound the distance from $Q_{d+1,s}$, from which in turn via~\eqref{eqDistanceQZn_vs_QZ_dplus1} and~\eqref{eqDistancePY_vs_QZn} we obtain an upper bound for $\dist_v(P,X)$.
This argument would be closer to the original one from~\cite{JiKollarShiffman92LojasiewiczInequality}, and perhaps simpler.
The disadvantage is that $\phi_{d+1,s}(X)$ may have smaller degree than $X$, and this degree enters the distance estimates as an exponent. 
\end{remark}

\section{Proof by example scheme}
\label{secProofByExampleScheme}

\subsection{Robust main theorem}
\label{secRobustMainTheorem}

\begin{theorem}[Robust main theorem]
\label{thmRobustMainTheorem}
Let $f_1,\ldots,f_m\in K[x_1,\ldots,x_n]$ be polynomials over a number field $K$, such that the associated 
variety $X = V(f_1,\ldots,f_m)\subseteq \AA^n$ is $d$-dimensional and irreducible over~$K$. 
Let $g\in K[x_1,\ldots,x_n]$.
Let the maximal degree be denoted by $D:=\max(\deg f_1,\ldots,\deg f_m,\deg g)$.
Choose any normalized valuation $v \in M_K$ of~$K$.
Let $P\in K^n$ be a point in the ball $|P|_v\leq R$ for some $R\geq 1$.
Suppose the coordinates of $P=(p_1,\ldots,p_n)$ satisfy
\begin{equation}
\label{eqRobustMainTheorem_assumptionOnHpi}
h(p_i) \geq nD^{m+1}\Big(h(f_1,\ldots,f_m,g,p_1,\ldots,p_{i-1}) + 4\log(n+2)\Big), \quad\textnormal{for $i=1,\ldots,d$}.
\end{equation}
Let
\begin{equation}
\label{eqRobustMainTheorem_eps}
\begin{split}
\log\eps := -4\tfrac{[K:\QQ]}{N_v}(n+1)^2D^n
\bigg([K:\QQ]& h(f_1,\ldots,f_m,g,p_1,\ldots,p_d)\ + \\ & \log m + (n+7)D\log(n+1) + \tfrac{\log R}{n+1}   \bigg).
\end{split}
\end{equation}
If 
\begin{equation}
\label{eqRobustMainTheorem_boundsForFiPv_and_gPv}
|f_i(P)|_v \leq \eps\quad  (\textnormal{for }1\leq i\leq m)\quad \textnormal{ and }\quad |g(P)|_v \leq \eps,
\end{equation}
then $g|_{X(\Qbar)} = 0$.
\end{theorem}

Note that $g|_{X(\Qbar)}=0$ if and only if $g\in \langle f_1,\ldots,f_m\rangle^\rad$ if and only if $g|_{X(L)} = 0$ for any extension field $L\supseteq \Qbar$.
Thus in this case we may simply write $g|_X=0$.

\begin{remark}[Dependence of $\eps$ on $P$]
The strength of the theorem comes from the fact that $\eps$ depends only on the \emph{local} height $h_v(P)$ of~$P$, as we can always choose the radius $R=\exp h_v(P)$ in the theorem.
In practice one chooses $R$ a bit larger, such that one can use the same $\eps$ for several ``examples''~$P$ that are chosen closer and closer to $X$, such that eventually one $P$ satisfies~\eqref{eqRobustMainTheorem_boundsForFiPv_and_gPv}.
\end{remark}

\begin{remark}[Weakening the height assumptions]
\label{remWeakeningHeightAssumptionInRobustMainTheorem}
A closer inspection of the proof reveals that the theorem will still hold if~\eqref{eqRobustMainTheorem_assumptionOnHpi} is replaced by the following weaker hypothesis.
Let $n_i:=\min(n,m+i)$, and let $(D_1,D_2,D_3,\ldots)$ be the sequence $(\deg f_1,\ldots,\deg f_m,\deg g,1,1,1,\ldots)$ ordered in non-increasing order. 
Then instead of \eqref{eqRobustMainTheorem_assumptionOnHpi}, we only need that for each $1\leq i\leq d$,
\begin{equation}
\label{eqRobustMainTheorem_assumptionOnHpi_moreCarefulEstimate}
\begin{split}
\qquad h(p_i)\geq \Big(h(f_1,\ldots,f_m, g,p_1,\ldots, & p_{i-1}) \sum_{j=1}^{n_i}\tfrac{1}{D_j} + (n+n_i) \log(n+1)\Big) \prod_{j=1}^{n_i}D_j \ + 
\\ & D^{n_i}(d-i+1)\big(\log(n+2)+{\textstyle\sum_{j=1}^m}\tfrac{1}{2j}\big) + \log 2.
\end{split}
\end{equation}
\end{remark}

\begin{remark}[Relation to Hilbert's Nullstellensatz]
We remark that Theorem~\ref{thmRobustMainTheorem} can be interpreted as a ``robust one-point Nullstellensatz'', in the sense that the inclusion $g\in I(X)$ can be verified by checking the error tolerant condition~\eqref{eqRobustMainTheorem_boundsForFiPv_and_gPv} for a single suitable point~$P$, instead of checking $g(P)=0$ for every $P\in X(\Qbar)$ according to Hilbert's Nullstellensatz.
Also, $P$ can be seen as a witness or a certificate that $g|_X=0$.
We refer to the discussion below Theorem~\ref{thmRobustMainTheorem_IntroVersion} from the introduction.
\end{remark}

\begin{proof}[Proof of Theorem~\ref{thmRobustMainTheorem}]
Let $Y_0:=X\cap V(g)$.
By Corollary~\ref{corDegreeBoundForVarietyXinTermsOfFis}, $\deg X\leq D^m$ and $\deg Y_0\leq D^{m+1}$.
Our aim is to show $\dim Y_0=d$. (This will also show the extension in Section~\ref{secReducibleMainTheorem}.)

So assume for a contradiction that $\dim Y_0 \leq d-1$.
For $i=1,\ldots,d$, define recursively $Y_i := Y_{i-1}\cap V(x_i-p_i)$.
By Theorem~\ref{thmGeometricBezout}, $\deg Y_i\leq D^{m+1}$.
Let $H_i:=h(f_1,\ldots,f_m,g,p_1,\ldots,p_i)$.
Using Corollary~\ref{corHeightOfXboundedInTermsOfHeightOfFis} we bound the height of $Y_i$ as follows,
\[
h(Y_{i-1}) \leq nD^{m+1}\big(H_{i-1} + 2\log(n+1)\big).
\]
By assumption~\eqref{eqRobustMainTheorem_assumptionOnHpi},
\[
h(p_i)\geq nD^{m+1}(H_{i-1}+4\log(n+2)),
\]
hence
\[
h(p_i)-h(Y_{i-1}) \geq nD^{m+1}\big(2\log(n+2)+2\log\tfrac{n+2}{n+1}\big) \geq 2D^{m+1}n\log(n+2) + \log 2.
\]
Also note that $\dim Y_{i-1} \leq \dim Y_0 \leq n-1$.
Then from Corollary~\ref{corContainmentInCoordinateHyperplane} we see that no component (over $K$) of $Y_{i-1}$ is contained in the coordinate hyperplane $V(x_i-p_i)$, and that thus $\dim Y_i\leq \dim Y_{i-1}-1$ if $Y_{i-1}$ is not already empty.

As $\dim Y_0\leq d-1$ it follows that $Y_d(\Qbar)=\emptyset$.
Therefore we can apply our effective \lojasiewicz{} inequality for empty varieties, Theorem~\ref{thmLojasiewiczForEmptySet}, to~$Y_d = V(f_1,\ldots,f_m,g,x_1-p_1,\ldots,x_d-p_d)$ and to the point~$P$ with $R=h_v(P)$.
It follows that either $|g|_v>\eps$ or $|f_i(P)|_v>\eps$ for some $i\in\{1,\ldots,m\}$, which in either case yields the desired contradiction.
\end{proof}

\subsection{Reducible case}
\label{secReducibleMainTheorem}

If $X$ is not assumed to be irreducible over~$K$, then the same theorem holds if the assertion 
``$g|_{X(\Qbar)} = 0$''
 is replaced by 
\begin{quotation}
``$g$ vanishes on some $d$-dimensional irreducible component of $X$ over $K$''.
\end{quotation}
The proof of Theorem~\ref{thmRobustMainTheorem} 
works equally well for this seemingly more general statement.
The next question is, on \emph{which} $d$-dimensional component does $g$ vanish?

To answer this, consider the variety $Y_d = X\cap V(g,x_1-p_1,\ldots,x_d-p_d)$, which is non-empty as the proof of the main theorem showed.
Thus we may apply the effective \lojasiewicz{} inequality, Theorem~\ref{thmLojasiewiczForNonEmptyVarieties}, to~$Y_d$.
We obtain a point $Q\in X(\Qbar)$ such that $\log\dist_v(P,Q)\leq\eps'$, where 
\mathsmall{
\begin{equation}
\label{eqEpsPrimeForRobustMainTheoremForReducibleX}
\begin{split}
\log\eps'\ :=\  &
\frac{\log\max(|f_{i}(P)|_v,|g(P)|_v))}{4(n+1)(D^n+1)^{n+2}}\ + \\
& \tfrac{[K:\QQ]^2}{N_v}(n+7)^2(D^n+1)\big(h(f_1,\ldots,f_m,g,p_1,\ldots,p_d) + \log((m+1)nD^{2n}) + 21 \big)
 + 2\log R.
\end{split}
\end{equation}
}%
Let $Y_0'$ be a component of $Y_0 = X\cap V(g)$ that contains~$Q$.
We apply Corollary~\ref{corContainmentInCoordinateHyperplane} as in the proof of the main theorem and deduce that $\dim Y_0'\geq d$.
Since $\dim Y_0'\leq\dim X = d$, $Y_0'$ is a component of $X$ on which $g$ vanishes.
Summarizing, we obtain the following corollary.

\begin{corollary}[Robust main theorem for reducible~$X$]
\label{corRobustMainTheoremForReducibleX}
Let $f_1,\ldots,f_m$, $g$, $v$, $P$, and $R$ be given under the same hypothesis as in Theorem~\ref{thmRobustMainTheorem}, except that now $X$ is allowed to be reducible.
Extend $v$ to $\Qbar$.
Let $\eps'$ be given as in~\eqref{eqEpsPrimeForRobustMainTheoremForReducibleX}.
Then $g|_{X'}=0$ for some $d$-dimensional component $X'$ (irreducible over~$K$) of~$X$ with 
\[
\dist_v(P,X')\leq \eps'.
\]
\end{corollary}

\subsection{Dichotomy}
\label{secDichotomy}

Given the variety $X$ and the polynomial $g$ as in the proof by example scheme, there are two mutually exclusive cases:
Either $g|_X=0$ or $g|_X\neq 0$.
In this section we study the problem of deciding within the proof by example framework, which of these cases hold.

The following dichotomy theorem separates both cases $g|_X=0$ and $g|_X\neq 0$ far enough from each other such that one can use numerical methods in order to decide in which of the two cases a given $g$ falls.

\begin{theorem}[Dichotomy in the proof by example scheme]
\label{thmDichotomy}
Let $K$, $n$, $f_1,\ldots,f_m$, $g$, $X$, $d$, $P$, $v$, $R$ be given as in Theorem~\ref{thmRobustMainTheorem}.
Choose $D_f\geq \max(\deg(f_1,\ldots,\deg f_m)$ and $D_g\geq\deg g$.
Suppose the coordinates of $P=(p_1,\ldots,p_n)$ satisfy
\begin{equation}
\label{eqDichotomy_assumptionOnHpi}
h(p_i) \geq nD_f^{m}\Big(h(f_1,\ldots,f_m,p_1,\ldots,p_{i-1}) + 4\log(n+2)\Big), \quad\textnormal{for $i=1,\ldots,d$}.
\end{equation}
Choose $H\geq h(f_1,\ldots,f_m,p_1,\ldots,p_d)+h(g)$ and
\begin{align*}
\log \eps_f := &-4\tfrac{[K:\QQ]^2}{N_{v}}(n+7)^3(D_f^n+1)^{n+4}D_g\Big(H+2\log n+\log m +\log R+12\Big),\\
\log \eps_g := & -\tfrac{[K:\QQ]}{N_{v}}nD_f^{2n}D_g\Big(H+4\log(n+1)\Big).
\end{align*}
Suppose that 
\begin{equation}
\label{eqDichotomy_boundsForFiPv}
|f_i(P)|_v\leq\eps_f\quad \textnormal{for all } 1\leq i\leq m.
\end{equation}
Then exactly one of the following two cases holds.
\begin{enumerate}[label={Case \arabic*.},labelindent=\parindent,leftmargin=*]
\item
$|g(P)|_v \leq \eps_g$ and $g|_X = 0$.
\item
$|g(P)|_v \geq 2\eps_g$ and $g|_X \neq 0$. 
\end{enumerate}
\end{theorem}

Thus in order to decide whether or not $g|_X=0$, using this dichotomy we only need to compute $|g(P)|_v$ up to an additive precision of~$\eps_g$.
The point $P$ at which we evaluate $g$ needs to be sufficiently generic in a slightly weaker sense~\eqref{eqDichotomy_assumptionOnHpi} compared to~\eqref{eqRobustMainTheorem_assumptionOnHpi} in the main theorem, and additionally $P$ needs to be sufficiently close to $X$ in the sense of~\eqref{eqDichotomy_boundsForFiPv}.
Note that here, the allowed margin of error $\eps_f$ is considerably smaller than the $\eps$~\eqref{eqRobustMainTheorem_eps} from the main theorem; compare in particular the power of~$D$ in $\eps$ and the power of $D_f$ in~$\eps_f$.

In the theorem we consider separate degree bounds $D_f$ and $D_g$, unlike the bound $D=\max(D_f,D_g)$ we used in Theorem~\ref{thmRobustMainTheorem}, to obtain finer estimates, as in Section~\ref{secDimensionByExample} we apply Theorem~\ref{thmDichotomy} in a situation when $D_g$ is considerably larger than~$D_f$.

\begin{proof}[Proof of Theorem~\ref{thmDichotomy}]
We consider $K$ as a subfield of $\Qbar$ and extend $v$ to a valuation on~$\Qbar$.
Define $X_d:=X\cap V(x_1-p_1,\ldots,x_d-p_d)$.
Theorem~\ref{thmLojasiewiczForEmptySet} applied to $X_d$ at the point $P$ and using assumption~\eqref{eqDichotomy_boundsForFiPv} implies that $X_d(\Qbar)$ cannot be empty, $\dim X_d\geq 0$. 
Then as in the proof of the main theorem, the genericity assumption~\eqref{eqDichotomy_assumptionOnHpi} together with Corollary~\ref{corContainmentInCoordinateHyperplane} and $\dim X=d$ imply $\dim X_d=0$.
Define
\begin{equation}
\label{eqDichotomy_epsPQ}
\begin{split}
\log \eps_{PQ} :=\ & \tfrac{\log\eps_f}{4(n+1)(D_f^n+1)^{n+2}}\ + \\
& \tfrac{[K:\QQ]^2}{N_{v|K}}(n+7)^2(D_f^n+1)\Big(H+\log((m+d)nD_f^{2n})+21\Big) + 2\log R.
\end{split}
\end{equation}
We apply Theorem~\ref{thmLojasiewiczForNonEmptyVarieties} to $X_d$ and the point $P$, and obtain the existence of a point $Q\in X(\Qbar)$ such that 
\[
|P-Q|_v\leq \eps_{PQ}.
\]
One easily bounds $v\leq 1$, and hence 
\[
|Q|_v\leq |P|_v + |P-Q|_v \leq R+1\leq 2R \quad \textnormal{ and thus } \quad h_v(P,Q)\leq 2R.
\]
Moreover $h_v(g)\leq \tfrac{[K:\QQ]}{N_{v|K}}h(g)$.
Using Lemma~\ref{lemElementaryLocalHeightInequalities},
\[
\log |g(P)-g(Q)|_v \leq \log\eps_{PQ} + \tfrac{[K:\QQ]}{N_{v|K}}h(g) + (D_g-1)\log 2R + D_g\log(n+2) =: C_1.
\]
A simply but lengthier chain of estimates yields
\begin{equation}
\label{eqdichotomy_epsg_vs_gPminusgQ}
\log \eps_g\geq C_1 \geq \log |g(P)-g(Q)|_v.
\end{equation}

Let $L:=K(Q)$ be the field extension of $K$ generated by the coordinates of~$Q$.
Since $\dim X_d=0$, Lemma~\ref{lemDegreeOfPointOf0dimVariety} and Corollary~\ref{corDegreeBoundForVarietyXinTermsOfFis} imply
\[
[L:K]\leq \deg X \leq D_f^n.
\]
The exponent of $D_f$ on the right hand side can be taken to be $m$ (using Theorem~\ref{thmGeometricBezout} together with Corollary~\ref{corDegreeBoundForVarietyXinTermsOfFis}), but again $n$ will turn out to be easier to simplify later.
We get $[L:\QQ]\leq D_f^n[K:\QQ]$, and clearly $N_{v|L}\geq N_{v|K}$.  

Also since $\dim X_d=0$, Lemma~\ref{lemHeightOfAPointInequality} and Corollary~\ref{corHeightOfXboundedInTermsOfHeightOfFis} imply
\[
h(Q)\leq h(\{Q\}) \leq h(X_d) \leq nD_f^n(H + 2\log(n+1)),
\]
where the exponent of $D_f$ could again be replaced by~$m$, but $n$ leads to easier simplifications later.
With Lemma~\ref{lemElementaryHeightInequalities} we continue to bound
\[
h(g(Q))\leq h(g) + D_g(h(Q)+\log(n+1)).
\]
We consider the two mutually exclusive cases:

\medskip 

\noindent
Case 1: $g(Q)=0$. In this case we apply Theorem~\ref{thmRobustMainTheorem} to $X$, $g$ and the point $Q$ (in place of~$P$), and we immediately obtain $g|_X=0$.
Moreover, from~\eqref{eqdichotomy_epsg_vs_gPminusgQ} we obtain $|g(P)|_v\leq \eps_g$.

\medskip

\noindent
Case 2: $g(Q)\neq 0$. In this case, obviously $g|_X\neq 0$.
We use Lemma~\ref{lemLiouville} to bound
\[
\log |g(Q)|_v \geq -\tfrac{[L:\QQ]}{N_{v|L}}h(g(Q)) \geq -D_f^n\tfrac{[K:\QQ]}{N_{v|K}}h(g(Q)).
\]
Putting bounds together and using simple estimates, we obtain
\[
\log |g(Q)|_v \geq \log (2\eps_g).
\]
The latter inequality together with~\eqref{eqdichotomy_epsg_vs_gPminusgQ} yields
\[
|g(P)|_v \geq |g(Q)|_v - |g(P)-g(Q)|_v  \geq 2\eps_g - \eps_g = \eps_g. 
\]
\end{proof}

\paragraph{How to use the dichotomy.}
In applications, this dichotomy can be used as follows:
Suppose our task is to determine whether or not $g|_X=0$, where we are in the setting of Theorem~\ref{thmRobustMainTheorem}.
To do this, the first step is to construct a point $P\in\QQ^n$ that is sufficiently generic in the sense of~\eqref{eqDichotomy_assumptionOnHpi} and that is sufficiently close to $X$ in the sense of $\max_i|f_i(P)|_v\leq \eps_f$.
This is the difficult part.
After having constructed such a point~$P$, one only needs to compute $|g(P)|_v$ up to sufficient precision (namely up to an additive error of $\eps_g$) such that one knows which of the two exclusive cases hold, $|g(P)|_v\leq \eps_g$ or $|g(P)|_v\geq 2\eps_g$.
By the dichotomy, these two cases correspond exactly to the cases $g|_X=0$ and $g|_X\neq 0$, which solves the task.

\paragraph{Reducible case.}
Let us also consider the dichotomy theorem in the seemingly more general case when $X$ is not assumed to be irreducible over~$K$.
Using essentially the same proof we obtain the following corollary.

\begin{corollary}[Dichotomy for reducible~$X$]
\label{corDichotomyForReducibleX}
Let $K$, $n$, $f_1,\ldots,f_n$, $g$, $X$, $d$, $D_f$, $D_g$, $P$, $v$, $R$ be given as in Theorem~\ref{thmDichotomy}, except that $X$ is not required to be irreducible.
Suppose that $P$ satisfies~\eqref{eqDichotomy_assumptionOnHpi} and~\eqref{eqDichotomy_boundsForFiPv}, and let $\eps_{PQ}$ be given as in~\eqref{eqDichotomy_epsPQ}.
Consider $K\subseteq \Qbar$ and extend $v$ to~$\Qbar$.
Then exactly one of the following two cases holds.
\begin{enumerate}[label={Case \arabic*.},labelindent=\parindent,leftmargin=*]
\item
$|g(P)|_v \leq \eps_g$ and $g|_{X'}=0$ for some $d$-dimensional $K$-component $X'$ of $X$ with $\dist_v(P,X')\leq \eps_{PQ}$.
\item
$|g(P)|_v \geq 2\eps_g$ and $g|_X \neq 0$. More precisely, $g(Q)\neq 0$ for some $Q\in X(\Qbar)$ with $\dist_v(P,Q)\leq\eps_{PQ}$.
\end{enumerate}
\end{corollary}

\subsection{Exact main theorem}
\label{secExactMainTheorem}

As a corollary of Theorem~\ref{thmRobustMainTheorem} we obtain the following exact main theorem, ``exact'' because we require that all $f_i(P)$ and $g(P)$ vanish exactly. 

\begin{corollary}[Exact main theorem]
\label{corExactMainTheorem}
Let $f_1,\ldots,f_m\in K[x_1,\ldots,x_n]$ be polynomials over a number field~$K$, and let the associated affine variety $X=V(f_1,\ldots,f_m)\subseteq\AA^n$ be $d$-dimensional and irreducible over~$K$. 
Let $g\in K[x_1,\ldots,x_n]$.
Let $P\in X(K)$ be a $K$-rational point with coordinates $P=(p_1,\ldots,p_n)$ that satisfy the height bounds~\eqref{eqRobustMainTheorem_assumptionOnHpi} (or alternatively, the weaker bounds~\eqref{eqRobustMainTheorem_assumptionOnHpi_moreCarefulEstimate})
and $g(P) = 0 $.
Then $g|_X = 0$.
\end{corollary}

\begin{proof}[Proof of Corollary~\ref{corExactMainTheorem}]
Choose an arbitrary valuation $v\in M_K$ and apply Theorem~\ref{thmRobustMainTheorem}.
\end{proof}

\paragraph{Reducible case.}
If $X$ is not assumed to be irreducible over $K$, then the same corollary holds if the assertion 
``$g|_X = 0$''
 is replaced by 
\begin{quotation}
``$g$ vanishes on any $d$-dimensional component of $X$ that passes through~$P$''.
\end{quotation}

The proof is analogous to the one of Corollary~\ref{corRobustMainTheoremForReducibleX}.
Note that in this exact setting, it works indeed for \emph{any} choice of component of $Y_0$ that contains~$P$.

\subsection{Measuring dimension by example}
\label{secDimensionByExample}

In the robust and exact main theorems, Theorem~\ref{thmRobustMainTheorem} and Corollary~\ref{corExactMainTheorem}, we assumed that the underlying variety $X$ is irreducible and $d$-dimensional.
We furthermore saw natural generalizations of the main theorems in case $X$ is only reducible and $d$-dimensional.
In this section we discuss how to proceed without the knowledge of the dimension of $d$, as it may indeed not be known in practice.
In particular the following theorem can be used similarly to the proof by example scheme to determine $\dim X$ by checking a certain inequality~\ref{eqDimByExample_CriterionOnDeterminantBeingVlarge} for a point~$P$ that is sufficiently generic~\eqref{eqDimByExample_HeightAssumptionForHpi} and sufficiently close to~$X$~\eqref{eqDimByExample_fiP_isVcloseTo0}, as in the robust main theorem.

\begin{theorem}[Dimension by example]
\label{thmDimensionByExample}
Let $K$ be a number field and $v\in M_K$. 
Let $f_1,\ldots,f_m\in K[x_1,\ldots,x_n]$ define an irreducible affine variety $X=V(f_1,\ldots,f_m)\subseteq\AA^n$.
Let $D_f:=\max(\deg f_1,\ldots,\deg f_m)$.
Let $0\leq d\leq n$.
Let $P=(p_1,\ldots,p_n)\in K^n$ be a point in the ball $|P|_v \leq R$ for some $R\geq 1$, whose first $d$ coordinates satisfy the height bounds 
\begin{equation}
\label{eqDimByExample_HeightAssumptionForHpi}
h(p_i) \geq nD_f^{m}\Big(h(f_1,\ldots,f_m,p_1,\ldots,p_{i-1}) + 4\log(n+2)\Big), \quad\textnormal{for $i=1,\ldots,d$}.
\end{equation}
Let $H:=h(f_1,\ldots,f_m,p_1,\ldots,p_d)$ and 
\begin{align*}
\log \eps_f'\ :=\ &-4\tfrac{[K:\QQ]^2}{N_{v}}(n+7)^3(D_f^n+1)^{n+5}\Big(H+2\log n+\log m +\log R+12\Big),\\
\log \eps_{\det}\ :=\ & -\tfrac{[K:\QQ]}{N_{v}}nD_f^{3n}\Big(H+4\log(n+1)\Big).
\end{align*}
Suppose
\begin{equation}
\label{eqDimByExample_fiP_isVcloseTo0}
|f_i(P)|_v\leq \eps_f'\quad\textnormal{for all $i=1,\ldots,m$}
\end{equation}
and
\begin{equation}
\label{eqDimByExample_CriterionOnDeterminantBeingVlarge}
|\det(e_1,\ldots,e_d,\nabla f_1(P),\ldots,\nabla f_{n-d}(P))|_v > \eps_{\det}.
\end{equation}
Then $\dim X = d$.
\end{theorem}

\begin{remark}[Reducible case]
If $X$ is not assumed to be irreducible over~$K$, then the same theorem holds if the assertion ``$\dim X = d$'' is replaced by
\begin{quotation}
``$X$ has a $d$-dimensional component $X'$ with $\dist_v(P,X')\leq \eps_{PQ},$''
\end{quotation}
where $\eps_{PQ}$ is given in~\eqref{eqDichotomy_epsPQ}.
\end{remark}

\begin{proof}[Proof of Theorem~\ref{thmDimensionByExample}]
Let $M$ be the $n\times n$ matrix $(e_1,\ldots,e_d,\nabla f_1,\ldots,\nabla f_{n-d})$ with entries in $K[x_1,\ldots,x_n]$.
Consider the polynomial $g:=\det M$.
Note that $|\tfrac{d}{dx_j} f_i|_v\leq \deg f_i\,|f_i|_v$.  
Hence $h_v(M_{(x_1,\ldots,x_n)}) \leq h_v(f_1,\ldots,f_m) + \log D_f$.
Using Lemma~\ref{lemElementaryLocalHeightInequalities}, we obtain
\[
h_v(g)\leq (n-d)(h_v(f_1,\ldots,f_m)+\log D_f) + \delta_{v\mid\infty}(n-d)(\log(n-d) + D_f\log(n+1)).
\]
and hence
\[
h(g)\leq (n-d)\big(h(f_1,\ldots,f_m)+\log D_f + \log(n-d) + D_f\log(n+1)\big).
\]
Clearly, $\deg g \leq (D_f-1)^{n-d}\leq D_f^n=:D_g$.
We apply the dichotomy results from Section~\ref{secDichotomy} to~$X$, $g$, and~$P$, in particular Corollary~\ref{corDichotomyForReducibleX}.
By assumption~\ref{eqDimByExample_CriterionOnDeterminantBeingVlarge}, $g(P)> \eps_{\det}$, which implies that we are in Case 2 of the dichotomy of Corollary~\ref{corDichotomyForReducibleX}.
Hence there exists $Q\in X(\Qbar)$ such that $g(Q)\neq 0$ (and $\dist_v(P,Q)\leq \eps_{PQ}$, with $\eps_{PQ}$ as in~\eqref{eqDichotomy_epsPQ}).
Following the proof of the Corollary~\ref{corDichotomyForReducibleX}, we may assume that $Q\in X_d(\Qbar)$.
This implies that $\nabla f_1(Q),\ldots,\nabla f_{n-d}(Q)$ are linearly independent over~$\Qbar$.
Therefore any component $X'$ of $X$ through $Q$ has dimension at most~$d$: $\dim X'\leq d$.

Let $X_d':=X'\cap X_d$.
Using Corollary~\ref{corContainmentInCoordinateHyperplane} together with assumption~\eqref{eqDimByExample_HeightAssumptionForHpi}, we see that either $X_d'=\emptyset$ or $\dim X' = \dim X_d'+d$.
As $Q \in X_d'(\Qbar)$, the latter holds and $\dim X'\geq d$.
Together both inequalities prove that $\dim X'=d$.
Since $X$ is irreducible we have $X=X'$, which finishes the proof.
\end{proof}

\paragraph{On the strength of Theorem~\ref{thmDimensionByExample} and a converse.}

The strength of Theorem~\ref{thmDimensionByExample} comes from the fact that the required bound on the valuation of the determinant~\eqref{eqDimByExample_CriterionOnDeterminantBeingVlarge} is reasonably weak. If this bound was too large, the theorem would rarely apply and thus become useless.
In fact the bound is so weak that it is ``almost equivalent'' to $\dim X=d$, in the following sense.

Assume that for any permutation of the $f_1,\ldots,f_m$, \eqref{eqDimByExample_HeightAssumptionForHpi} and~\eqref{eqDimByExample_fiP_isVcloseTo0} hold but not~\eqref{eqDimByExample_CriterionOnDeterminantBeingVlarge}.
Then by the dichotomy, there exists a point $Q\in X_d(\Qbar)$ with $\dist_v(P,Q)\leq \eps_{\det}$ such that the $\Qbar$-linear span $N_Q (X):=\Qbar\{\nabla f_1(Q),\ldots,\nabla f_m(Q)\}$ is at most $(n-d-1)$-dimensional.

Assuming that $X$ has indeed dimension $d$ and that $f_1,\ldots,f_m$ generate $I(X)$ (i.e. that $\langle f_1,\ldots,f_m\rangle$ is a radical ideal), then $\dim N_Q(X) \leq n-d-1$ means that $Q$ is a singular point of~$X$. 
In other words: 

\begin{corollary}[A converse of Theorem~\ref{thmDimensionByExample}]
\label{corDimensionByExampleConverse}
If $\langle f_1,\ldots,f_m\rangle = I(X)$, if $X$ is smooth, and if~\eqref{eqDimByExample_HeightAssumptionForHpi} and~\eqref{eqDimByExample_fiP_isVcloseTo0} hold, then $\dim X=d$ is equivalent to  the statement that the equations~\eqref{eqDimByExample_CriterionOnDeterminantBeingVlarge} hold for all permutations of $f_1,\ldots,f_m$.
\end{corollary}

A version of this last statement can also be proved in case $X$ is not smooth:
If $X$ is not smooth, then its singular points form a subvariety $S$ of $X$ which is cut-out by the vanishing of all $(n-d)$-minors of the matrix $(\nabla f_1,\ldots,\nabla f_m)$.
In this way, $h(S)$ can be explicitly bounded in terms of $h(X)$, $m$, and~$n$.
Thus under a suitably stronger version of assumption~\eqref{eqDimByExample_HeightAssumptionForHpi}, one can argue via Corollary~\ref{corContainmentInCoordinateHyperplane} that $S\cap X_d$ is empty, and hence $Q$ a regular point of~$X$.

\begin{remark}[Robustness of~\eqref{eqDimByExample_CriterionOnDeterminantBeingVlarge}]
A last point to note is that~\eqref{eqDimByExample_CriterionOnDeterminantBeingVlarge} can be checked numerically: Under the assumptions~\eqref{eqDimByExample_HeightAssumptionForHpi} and~\eqref{eqDimByExample_fiP_isVcloseTo0}, the dichotomy implies that $|\det(\ldots)|_v>\eps_{\det}$ holds if and only if $|\det(\ldots)|_v\geq 2\eps_{\det}$.
Thus $|\det(\ldots)|_v$ needs to be computed only up to an additive precision of~$\eps_{\det}$. 
\end{remark}

\subsection{Sufficiently generic points}
\label{secSufficientlyGenericPoints}

In the robust main theorem, the chain of height inequalities~\eqref{eqRobustMainTheorem_assumptionOnHpi} ensured that the example~$P$ is sufficiently generic.
This criterion is easy to check and valid for a set of points~$P$ which is obviously dense with respect to the $v$-norm for any $v\in M_K$, however the upper Banach density of this set is close to zero.
Now we ask, whether this criterion can be generalized, such that in applications one becomes more flexible in constructing a suitably generic example.

\begin{definition}
Let $H\in\RR_{>0}$, $D\in\ZZ_{\geq 1}$, $K$ be a number field, $v\in M_K$, and $d\in\ZZ_{\geq 1}$.
We call a point $Q\in (\CC_v)^d$ \textdef{$(H,D,K)$-generic} if $f(Q)\neq 0$ for any polynomial $f\in K[x_1,\ldots,x_d]$ with $h(f)\leq H$ and $\deg f\leq D$.
\end{definition}

Note that $(\infty,\infty,K)$-generic points are exactly the points in $(\CC_v)^d$ with algebraically independent components, independently of the number field~$K$.
Thus $(H,D,K)$-genericity is a relaxation of algebraic independence.
Also, for given $H,D,K$, the non-$(H,D,K)$-generic points lie in the union of finitely many hypersurfaces $V(f)$ with $h(f)\leq H$, $\deg f\leq D$, and this union is an algebraic variety.
In this sense, most of the points of $(\CC_v)^d$ are $(H,D,K)$-generic. 

\begin{corollary}
\label{corRobustMainTheorem_forHDKgenericPoints}
Theorem~\ref{thmRobustMainTheorem} also holds when the genericity condition~\eqref{eqRobustMainTheorem_assumptionOnHpi} is replaced by the requirement that $Q := (p_1,\ldots,p_d)$ is $(H_0,D_0,K)$-generic with $D_0 := dD^{m+1}$ and
\[
H_0 := nD^{m+1}\Big(h(f_1,\ldots,f_m,g) + 3\log(n+2)\Big).
\]
\end{corollary}

\begin{proof}
As in the proof of the main theorem, we assume that $Y_0=X\cap V(g)$ has $\dim Y_0=d-1$.
Let $f_{Y_0}$ be a normalized Chow form for~$Y_0$.
Let $e_1,\ldots,e_d$ be the standard basis vectors of $\QQ^d$.
Let $f(\alpha_1,\ldots,\alpha_d):=f_{Y_0}((-\alpha_1,e_1),\ldots,(-\alpha_d,e_d)) \in K[\alpha_1,\ldots,\alpha_d]$. 
By definition of $f_{Y_0}$, any zero $(\alpha_1,\ldots,\alpha_d)$ of $f$ corresponds to a $d$-dimensional hyperplane $\{x \st x_1=\alpha_1,\ldots,x_d=\alpha_d\}$ whose projective closure intersects the projective closure of~$Y_0$.
Thus if we can show that $f(Q)\neq 0$, then $Y_d:=Y_0\cap V(x_1-p_1,\ldots,x_d-p_d)$ is the empty variety and we can proceed with $Y_d$ as in the proof of the main theorem.

Let $\delta_0 := \deg Y_0 \leq D^{m+1}$, where as before, $D=\max(\deg f_1,\ldots,\deg f_m,\deg g)$.
Hence 
\[
\deg f\leq \deg f_{Y_0}\leq d\delta_0 \leq dD^{m+1} = D_0.
\]
As for the height, 
\[
h(f)\leq h(f_{Y_0}) \leq h(Y_0) + \delta_0d \big(\log(n+2) + \sum_{i=1}^n \tfrac{1}{2i}\big)
\]
by Lemma~\ref{lemHeightOfXvsChowFormOfX}.
Furthermore, 
\[
h(Y_0) \leq nD^{m+1}\big(h(f_1,\ldots,f_m,g)+2\log(n+1)\big).
\]
Putting both bounds together, we obtain $h(f)\leq H_0$.
As $Q$ was assumed to be $(H_0,D_0,K)$-generic, it follows $f(Q)\neq 0$, which remained to be proved.
\end{proof}

Corollary~\ref{corRobustMainTheorem_forHDKgenericPoints} offers a considerably more general genericity assumption compared to~\eqref{eqRobustMainTheorem_assumptionOnHpi} from Theorem~\ref{thmRobustMainTheorem}, however in practice it will be more difficult to check whether a given point~$P$ satisfies this condition.
The trivial algorithm to check $(H,D,K)$-genericity would verify $f(Q)\neq 0$ for all $f$ from the definition, however this is far from practical.
Whether a practical algorithm exists is unclear.

In cases where~\eqref{eqRobustMainTheorem_assumptionOnHpi} is too restrictive, one may try to borrow methods from transcendence theory that are usually build to prove algebraic independence.
We in particular refer to Philippon~\cite{philippon1986criteresIndependanceAlgebrique},
Jabbouri~\cite{Jabbouri92approximationDiophantiennes},
Laurent and Roy~\cite{LaurentRoy01criteriaOfAlgebraicIndependence},
and Philippon~\cite[Ch. 8]{NesterenkoPhilippon01introductionToAlgebraicIndependenceTheory}.
For example, \cite[Ch. 8]{NesterenkoPhilippon01introductionToAlgebraicIndependenceTheory} contains a quite general ``criterion for algebraic independence'' with several adjustable parameters that in fact offers sufficient criteria for $(H,D,K)$-genericity.

\begin{remark}[Sufficiently generic points for the \emph{exact} proof by example scheme.]
Suppose $X\subseteq \AA^n$ is a given irreducible $d$-dimensional variety, and $g\in K[x_1,\ldots,x_n]$ is a given polynomial for which we want to show $g|_X=0$.
In the \emph{exact} proof by example scheme, we need a point $P\in X(\Qbar)$ such that $g(P)=0$ implies $g|_X=0$.
If $g$ satisfies bounds $h(g)\leq H$ and $\deg g\leq D$ that are known to us, then a point $P\in X(\Qbar)$ is sufficiently generic for $g$ if $P$ does not lie in the union of all $X\cap V(q)$ with $q\in K[x_1,\ldots,x_n]\wo I(X)$, $h(q)\leq H$, and $\deg q\leq D$.
By Theorems~\ref{thmGeometricBezout} and~\ref{thmArithmeticBezout}, all $\Qbar$-rational points in $X\wo \bigcup W$ are therefore sufficiently generic for $g$, when the union runs over all $(d-1)$-dimensional subvarieties $W$ of $X$ with
\[
\deg W\leq D\deg X \quad \textnormal{ and } \quad h(W)\leq D\big(h(X)+\deg X(H+\log(n+1))\big).
\]
\end{remark}

\section{Proving theorems in plane geometry}
\label{secPlaneGeometry}

In this section we consider theorems from plane geometry in order to discuss the strengths and weaknesses of the proof by example scheme.
Disclaimer: Proof by example will \emph{not} be a serious contender for the currently best geo\-metric automated theorem prover, as it is limited to statements of the type $g\in I(X)$, has obstructions such as the irreducibility of $X$, and may not be computationally efficient enough; compare with Section~\ref{secLimitsOfProofsByExample}.

Automated theorem proving in geometry has a long history. 
The first such algorithm to prove (or disprove!) statements in elementary euclidean geometry was obtained by Tarski by showing that his first-order axiomatization of elementary geometry is decidable.
There are various more practical approaches today:
Wu's method~\cite{Wu78decisionProblemAndTheoremProving, Wu84mechanicalTheoremProvingInElementaryGeometries} is based on Ritt's characteristic set method~\cite{Ritt50differentialAlgebra}, see also Chou~\cite{Chou85thesis_wusMethod}.
Gr\"obner bases by Buchberger~\cite{Buchberger2006phd} were used for example by Kapur~\cite{Kapur86groebnerBasesForGeometryProblems}.
 Chou, Gao, and Zhang developed an area method~\cite{ChouGaoZhang96automatedProofsI} and a full angle method~\cite{ChouGaoZhang96automatedProofsII}.
This is only a small extract of a vast literature.
Two of the current implementations of automated theorem provers in geometry are GCLC~\cite{Janivic06GCLC} and the OpenGeoProver within GeoGebra~\cite{geogebra5}.

The strength to prove many theorems is only one criterion for automated theorem provers.
Others are 
for example the human readability of the constructed proofs. 
The author will let the reader decide whether proofs by example are particularly readable for humans.
Here, important aspects are not only the easiness of verification of a proof, but also the conceptual understanding of why a statement is true.

\subsection{Example: Thales' theorem}
\label{secThales}

We choose Thales' theorem not only because it is a minimal non-trivial example for our proof by example scheme, but also for its historic value.
Thales of Miletus ($\sim$\,600 BC) is the first known individual person to whom a mathematical discovery is attributed, 
and his most famous one is the theorem that is named after him.
In Euclid's ``Elements''~\cite{Euclid02elements}, it  
is stated as follows: ``The angle in a semi-circle is a right-angle.'' See the first sentence of Proposition 31 in the third book of~\cite{Euclid02elements}.
Purportedly, Thales sacrificed an ox for his discovery.

\begin{theorem}[Thales]
If three distinct points $A$, $B$, $C$ lie on a circle $k$ such that the line segment $\overline{AB}$ is a diameter of~$k$, then $\angle BCA$ measures $90^\circ$.
\end{theorem}

We will give two similar proofs by example, one using $v=\infty$, the other with~$v=7$.

\begin{proof}[First proof of Thales' theorem by example, via $v=\infty$]
We may assume that $A$ and $B$ have Cartesian coordinates $A=(-1,0)$ and $B=(1,0)$, as this can be achieved by a similarity transformation.
The circle $k$ with diameter $\overline{AB}$ is thus the unit circle.

\paragraph{Setting up the proof by example scheme.}
Let $C$ have coordinates $C=(p_1,p_2)$.
As $C\in k$, $(p_1,p_2)$ is a root of the polynomial $f(p_1,p_2) := ||C||^2 - 1$, where $||(x,y)||^2 = x^2+y^2$ denotes the standard Euclidean norm of~$\RR^2$.
As there are no other steps in the construction and the points $A$, $B$ have been fixed in advance, the variety $X:=V(f)\subseteq \AA^2$ parametrizes all possible constructions, in the sense that each $(p_1,p_2)\in X(\RR)$ corresponds to a point $C\in k$.

\begin{figure}[htb]
\centering
\begin{minipage}[h]{0.75\textwidth}
\centering
\includegraphics[scale=0.8]{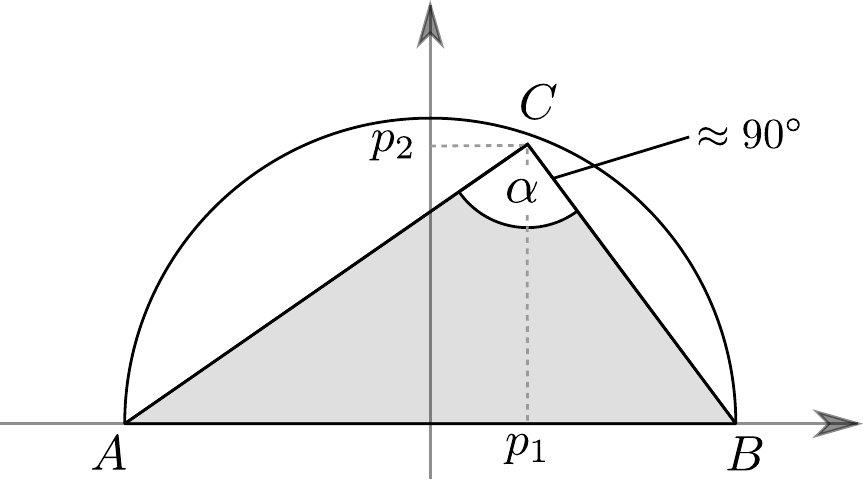}
\caption{A non-exact sketch for Thales' theorem. If it is sufficiently generic and precise, it leads to valid proof by example.}
\label{figThales}
\end{minipage}
\end{figure}

That $\angle BCA$ is a right angle can be encoded in saying that $(p_1,p_2)$ is a root of the polynomial $g(p_1,p_2) := \langle C-A,C-B\rangle$, where $\langle (x,y),(x',y')\rangle = xx'+yy'$ denotes the standard euclidean inner product of~$\RR^2$.

Thales' theorem now states that $g$ vanishes on $X(\RR)$.
We will show this by proving $g\in I(X)$ using the proof by example scheme as follows.

\begin{remark}
Of course in this simple setting one can directly compute $g(p_1,p_2) = p_1^2 + p_2^2 -1 = f(p_1,p_2)$, which thus yields an analytic proof of Thales' theorem.
\end{remark}

\paragraph{Irreducibility and dimension of $X$.}
Note that $f$ is non-zero and irreducible, hence $X$ is irreducible over~$\QQ$ and of dimension~$1$.

\paragraph{Bounding heights.}
As all non-zero coefficients of $f$ are $\pm 1$, we see that
\[
h(f)=0.
\]
For the sake of a good example, let us bound the height of $g$ naively, ignoring that in fact $h(g)=0$.
Suppose that we compute $(p_1+1)(p_1-1) = p_1^2-1$, which has height $0$, and that $(p_2-0)(p_2-0)=p_2^2$ has height $0$ as well, but when we take their sum we naively refer to Lemma~\ref{lemElementaryHeightInequalities}, to obtain the bound
\[
h(g)\leq \log 2.
\]
Together we get 
\[
h(f,g)\leq h(f)+h(g) \leq \log 2.
\]

\paragraph{Choosing $p_1$.}
In order to obtain a sufficient lower bound for~$h(p_1)$ we use~\eqref{eqRobustMainTheorem_assumptionOnHpi_moreCarefulEstimate} from Remark~\ref{remWeakeningHeightAssumptionInRobustMainTheorem} instead of the more restrictive bound~\eqref{eqRobustMainTheorem_assumptionOnHpi} from the main theorem, simply in order to get a proof with smaller numbers.
According to~\eqref{eqRobustMainTheorem_assumptionOnHpi_moreCarefulEstimate}, we may choose any $p_1$ with height at least
\[
H := \big(h(f,g)2\tfrac{1}{2} + (2+2)\log(3)\big)2^2 + 2^2\cdot 1(\log 4+3/4) + \log 2.
\]
Using the bound $h(f,g)\leq \log 2$ from above, we compute that $H/\log 10 \leq 13$, 
and thus taking
\[
p_1 = 0.1234567890123 = 1234567890123/10^{13}
\]
will have a sufficiently large height.
One may choose any other $p_1\in\Qbar$ of sufficiently large height, notably also those with $|p_1|>1$, but let us keep things real.

We will use the valuation $v=\infty$, that is, $|\ .\ |_v$ is the usual norm of~$\CC$.
We need to declare a bound $R\geq 1$ for $|P|$.
As we are only interested in points $P$ on the unit circle, we choose $R:=2$, allowing for some comfortable numerical margin of error.
For this $p_1$ and $R$, we compute
\[
\eps \geq 10^{-1300},
\]
which corresponds to a precision of $1300$ decimal digits.

\paragraph{Solving for $p_2$.}
Given this $p_1$ and $f$, we have to find an at least approximate root of $f(p_1,\blank)$.
Since $1300\log_2 10\approx 4318.5$, we compute numerically $\wt{p_2}$ as $\sqrt{1-p_1^2}$ in interval arithmetic with $4330$ bits of precision (our computations were done in the computer algebra system Sage~\cite{sage}, see Figure~\ref{figAlgorithmInfinity} for the relevant part of our implementation).
That is, $\wt p_2$ is an exactly computed interval which contains the correct value of $\sqrt{1-p_1^2}$, and whose end-points are floating point numbers with $4330$ bits of relative precision.
In what follows we will compute that any choice for $p_2\in\wt{p_2}$ yields a suitable example $P=(p_1,p_2)$ for the proof by example scheme.
This freedom is depicted in Figure~\ref{figThales}.

\paragraph{Bounding $|f(P)|$, $|g(P)|$, and~$R$.}
We plug $(p_1,\wt{p_2})$ into $g$, again using interval arithmetic, and obtain an exactly computed interval
\[
g(p_1,\wt{p_2}) \subseteq [-1.1\cdot 10^{-1303},1.1\cdot 10^{-1303}].
\]
Thus we are certain that plugging in any value $p_2\in\wt{p_2}$ yields
\[
|g(p_1,p_2)| \leq 1.1\cdot 10^{-1303} \leq \eps.
\]
We similarly check that indeed $|(p_1,\wt{p_2})| \leq 1.0 \leq R$, and hence $|(p_1,p_2)|\leq R$ for any $p_2\in \wt{p_2}$.
We could do the same for $f(p_1,\wt{p_2})$, but for simplicity we can just take the correct root $p_2=\sqrt{1-p_1^2}$ (without computing it numerically), such that automatically $f(p_1,p_2)=0$.

\paragraph{Finishing the proof.}
We apply Theorem~\ref{thmRobustMainTheorem} and obtain that $g$ vanishes on $X=V(f)$.
This finishes the proof of Thales' theorem.
\end{proof}

\begin{figure}[thb]
 \centering
  \begin{minipage}{0.45\textwidth}
   ~~
   \begin{minipage}{0.95\textwidth}
    \footnotesize
   \begin{verbatim}
K = RealIntervalField(4330)      
p1 = K(1234567890123/10^13) 
p2 = sqrt(1-p1^2) 
gP = (p1-1)*(p1+1)+(p2-0)*(p2-0) 
print "g(p1, p2) >=", gP.lower() 
print "g(p1, p2) <=", gP.upper()\end{verbatim}
   \end{minipage}
   \caption{Sage code that is used to bound $|g(p_1,p_2)|_\infty$ for our first proof of Thales' theorem.}
   \label{figAlgorithmInfinity}
  \end{minipage}
 \qquad
  \begin{minipage}{0.45\textwidth}
   ~~
   \begin{minipage}{0.95\textwidth}
    \footnotesize
    \begin{verbatim}
K = Qp(7, prec = 1525) 
p1 = K(7*1234567890123) 
p2 = sqrt(1-p1^2) 
gP = (p1-1)*(p1+1)+(p2-0)*(p2-0) 
print "g(p1, p2) =", gP
    \end{verbatim}
   \end{minipage}
   \caption{Sage code that is used to bound $|g(p_1,p_2)|_7$ for our second proof of Thales' theorem.}
   \label{figAlgorithmAt7}
  \end{minipage}
\end{figure}

In the proof we chose $v=\infty$, since it corresponds to the norm in $\Qbar$ coming from the usual norm in~$\CC$.
Because of this choice, we ran into a mild numerical challenge of bounding errors, which we solved using interval arithmetic.
The advantage of using $p$-adic norms instead is that the numerics become easier as one can work in a fixed precision (essentially because when adding two $p$-adic digits, carries have smaller norm, as opposed to larger norm when using $v=\infty$).

\begin{proof}[Second proof of Thales' theorem by example, via $v=7$]
The proof works in exactly the same way as the first one, constructing $f$, $g$, $R=2$, $H$, up to the point where we choose $p_1$.
We work over the field of $7$-adic numbers~$\QQ_7$. 
It will turn out that for the $p_1$ chosen below we have 
$\eps \geq 7^{-1525}$, and so we do computations within a precision of~$1525$ $7$-adic digits, i.e. the computations are done~in the ring $\QQ_7/7^{1525}\ZZ_7$: Since in this example all numbers turn out to be $7$-adic integers, i.e.~numbers from~$\ZZ_7$, taking this precision means that we simply do computations in $\ZZ/7^{1525}\ZZ$. See Figure~\ref{figAlgorithmAt7} for the most relevant part of our implementation in Sage.
We choose a value for $p_1$ with $|p_1|_7<1$ in order to avoid going into an extension field of~$\QQ_7$ when taking the square root below.
So let us take $p_1 = 7\cdot 1234567890123$, which has sufficient height, $h(p_1)\geq H$, and which in $\ZZ_7$ up to precision $1525$ reads
\mathsmall{
\[
7 + 5 \cdot 7^{2} + 5 \cdot 7^{3} + 3 \cdot 7^{5} + 6 \cdot 7^{6} + 2 \cdot 7^{7} + 7^{8} + 5 \cdot 7^{9} + 3 \cdot 7^{10} + 2 \cdot 7^{11} + 7^{12} + 5 \cdot 7^{13} + 5 \cdot 7^{14} + 7^{15} + O(7^{1525}).
\]
}%
Next we compute $p_2=\sqrt{1-p_1^2}$ within $\ZZ_7$ up to precision $1525$, the first and last few digits of which are
\[
p_2 = 1 + 3\cdot 7^2 + 5\cdot 7^3 
+ \ldots +  2\cdot 7^{1523} + 6\cdot 7^{1524} + O(7^{1525}).
\]
Next we plug $p_1$ and $p_2$ into $f$ and $g$ and obtain
\[
f(p_1,p_2) = g(p_1,p_2) = 0 + O(7^{1525}).
\]
In other words, we could simply pick the rational integers (i.e.\ elements in~$\ZZ$)
\[
p_1 = 7\cdot 1234567890123 \quad \textnormal{ and } \quad p_2 = 1 + 3\cdot 7^2 + \ldots + 6\cdot 7^{1524},
\]
for which
\[
|f(p_1,p_2)|_7 \leq \eps, \quad |g(p_1,p_2)|_7 \leq \eps, \quad \textnormal{ and clearly }\quad |(p_1,p_2)|_7 = 1 \leq R.
\]
Applying Theorem~\ref{thmRobustMainTheorem} shows that $g\in I(X)$, which proves Thales' theorem.
\end{proof}

\subsection{More general ruler and compass constructions.}
\label{secGeneralElementaryGeometricConstructions}

We treat geometric constructions here in a way, which is not most efficient but easy to study.

For us a construction consists of the following three building blocks.
\begin{enumerate}
\item
Free objects, such as: points, lines, and circles.
\item
Equalities, obtained from requiring incidences such as ``point $A$ lies on line~$\ell$'' or relations such as ``lines $\ell_1$ and $\ell_2$ are parallel''.
\item
Inequalities, obtained from requirements of distinction such as ``points $P_1$ and $P_2$ are distinct'' or ``line $\ell$ is not a tangent of circle~$k$''.
\end{enumerate}

For the sake of this section, a geometric construction is an iterative recipe using these building blocks:
One generates free objects and then puts (algebraic) equalities and inequalities on them.
In this way, we build up a configuration space $X=V(f_1,\ldots,f_m)\subset \AA^n$ of all possible outcomes of a geometric construction.
At the beginning of the construction, $(n,m)=0$ and $X=\AA^0$.

\paragraph{Free objects.}
Let $X=V(f_1,\ldots,f_m)\subseteq \AA^n$ be the configuration space for a construction.
We can add a free point $P=(a,b)$ to the construction by introducing two new coordinates $p_{n+1},p_{n+2}$ to the ambient space, which correspond to the coordinates of~$P$.
The new configuration space is $X'=I(f_1,\ldots,f_m)\subseteq \AA^{n+2}$, determined by the same equations as~$X$.

Similarly, adding a line with equation $L: y=ax+b$ to the construction adds as well two new coordinates.
Note that lines of this form are exactly those that are not parallel to the $y$-axis.
We restrict to such lines in order to keep $X$ an affine variety.
Statements about constructions that involve lines parallel to the $y$-axis can be reduced to the one which avoids such lines altogether by one of the following two arguments:
\begin{enumerate}
\item Suppose the statement in question is rotation invariant, in the sense that if it holds for one construction then it also holds for any rotated copy of this construction.
Then, if a construction involves lines parallel to the $y$-axis, one can reduce it to the case without lines parallel to the $y$-axis by applying a suitable rotation.

\item Suppose the space of all constructions is the set $Y(\RR)$ of real points on some suitable projective variety $Y$ that contains $X$ as an open subvariety, such that $Y=\overline X$ is the projective closure of~$X$, and that the statement in question is that $\wt g = 0$ for some regular function $\wt g$ on~$Y$. 
Then let $g=\wt g|_X$, and suppose one can prove the statement $g=0$ on $X$ by example.
This means that $g\in I(X)$, which by continuity implies that $\wt g=0$ on~$Y$. 
\end{enumerate}

Similarly, adding a circle of the form $k: (x-a)^2 + (y-b)^2 = r^2$ or a conic of the form $C:ax^2+bxy+cy^2 + dx + ey + 1 = 0$ adds $3$ and $5$ new coordinates, respectively.
Note that circles and conics of this form may be degenerate; if this is not desired one can add an inequality (e.g. $r\neq 0$) as explained below.
Also note that conics of this form cannot go through the origin; if this restriction is not desired then one can argue as with lines above.

\paragraph{Equalities of the form $f=0$.}
Let $X=V(f_1,\ldots,f_m)\subseteq \AA^n$ be the configuration space for a construction.
We can impose any algebraic relations on the parameters of the free objects in the construction, simply by appending this relation as a further equation $f_{m+1}$ on $X$, thus obtaining a new configuration space $X'=I(f_1,\ldots,f_{m+1})\subseteq \AA^n$.

For example we can impose that two lines are parallel, or that two circles are tangent to each other.
Slightly less trivial are angles:
Let $\angle(\ell_1,\ell_2)$ denote the oriented mod-$\pi$ angle between two lines $\ell_i:y=a_ix+b_i$ ($i=1,2$).
This angle is $\alpha_2-\alpha_1$ mod $\pi$ if $\tan\alpha_i = b_i$ ($i=1,2$), and it is thus determined by $\tan(\alpha_2-\alpha_1) = (b_2-b_1)/(1+b_1b_2)$ (which we allow to be $\infty$ in the obvious way).
For four such lines $\ell_i$ ($i=1,2,3,4$), we obtain that $\angle(\ell_1,\ell_2) = \angle(\ell_3,\ell_4)$ mod $\pi$ if and only if $(b_2-b_1)(1+b_3b_4) - (b_4-b_3)(1+b_1b_2) = 0$, which is an equation of degree~$3$.

\paragraph{Inequalities of the form $f\neq 0$.}
Let $X=V(f_1,\ldots,f_m)\subseteq \AA^n$ be the configuration space for a construction.
For any polynomial $f(p_1,\ldots,p_n)$ in the $n$ coordinates of $\AA^n$ we can add the constraint $f(p_1,\ldots,p_n)\neq 0$ to the construction as follows: We append one ``dummy variable'' $p_{n+1}$ and the equation $f_{m+1}(p_1,\ldots,p_{m+1})=0$, where $f_{m+1}(p_1,\ldots,p_{m+1}) := 1 + p_{n+1}f(p_1,\ldots,p_n)$.
Let $X':= V(f_1,\ldots,f_{m+1})\subseteq \AA^{n+1}$, and let $\pi:\AA^{n+1}\to\AA^n$ be the projection to the first $n$ coordinates.
Then $\pi|_{X'}: X' \to X\wo\{f=0\}$ is an isomorphism of quasi-affine varieties, the inverse of which is given by $P\mapsto (P,-1/f(P))$.

For example we can require that two points are distinct.
Moreover for every equality from above that we could impose on the construction, we could equally well impose its negation on the construction instead.

\subsection{The associated proof by example scheme}

In the previous section we described how we obtain the configuration space 
\[
X=V(f_1,\ldots,f_m)\subseteq \AA^n.
\]
The statement $g|_X = 0$ will usually come from a polynomial $g$ that corresponds to one of the equalities as in the previous section.

The ambient dimension~$n$ and the number $m$ of generators of $I(X)$ depends on the number of building blocks of the construction.
Let $n_f$ be the number of parameters of all free objects.
Let $n_e$ be the number of equalities, and $n_i$ the number of inequalities.
Let $D_e$ be the maximal degree of an equality, and $D_i$ the maximal degree of an inequality.
Then $n = n_f+n_i$, $m = n_e+n_i$, and $D:=\max(\deg f_1,\ldots,\deg f_m,\deg g)  = \max(D_e,D_i+1,\deg g)$.
Plugging this into the main theorem, we obtain an explicit proof by example scheme, except that the dimension of $X$ and its irreducibility may be non-obvious.

\paragraph{Irreducibility.}
To achieve irreducibility, we usually have to remove the trivial components by adjoining further inequalities to the construction.
In case there are more than one irreducible components $X_i$ of interest for which we want to show $g|_{X_i}=0$, one has to apply the main theorem for each of these components separately, each time blending out the other components by adjoining appropriate inequalities. 
As all depends on the context, a general method for achieving irreducibility in a practical way is unlikely to exist.
Irreducibility can indeed be a major obstacle for the applicability of the proof by example scheme.

\paragraph{Dimension.}
Determining the dimension $d=\dim X$ will in general be easier than achieving irreducibility.
If the chosen equalities are seemingly independent from each other, the natural guess for the dimension of $X$ is $d=n-m$.
To check that this guess is indeed correct, we can make use of Theorem~\ref{thmDimensionByExample}.

A lower bound for the determinant~\eqref{eqDimByExample_CriterionOnDeterminantBeingVlarge} in Theorem~\ref{thmDimensionByExample} can be quickly obtained without large matrix computations if we restrict to certain natural subclasses of constructions: For the general constructions of Section~\ref{secGeneralElementaryGeometricConstructions} we allowed to simply add polynomial equalities and inequalities to the construction. 
If instead we start the construction from some free objects and create new objects only from previous ones via ruler and compass constructions, the matrix in Theorem~\ref{thmDimensionByExample} will be block-diagonal, with one block for every construction step, that is, for every newly constructed object.
The determinant of the whole matrix is the product of the determinants of the blocks.
So if we choose the initial free objects in such a way that every diagonal block has determinant with $v$-norm at least $(\eps_{\det})^{1/(n-d)}$, then condition~\eqref{eqDimByExample_CriterionOnDeterminantBeingVlarge} will hold.
A $v$-large determinant of the diagonal block of a newly constructed object just means that this new object was constructed in a sufficiently transversal way.
For example if this new object is the intersection of two already constructed lines, this means that the lines should not be `too parallel'.

\section{Discussion}
\label{secDiscussion}

We have seen that Theorem~\ref{thmRobustMainTheorem} provides a way to prove statements of the form $g\in I(X)$ by example.
However so far this scheme has 
practical limits.
 
\subsection{Limits of the proof by example scheme}
\label{secLimitsOfProofsByExample}

\paragraph{Non-algebraic statements.}
So far, the proof by example scheme works only in the algebraic setting of the main theorem. Already if we introduce inequalities of the form $f\geq 0$ to a construction in plane geometry, the scheme stops working.

\paragraph{Irreducibility.}
Checking that $X$ is irreducible or computing its connected component over $K$ might be very difficult. 
There does not seem to exist a natural way to show irreducibility ``by example'' due to the local nature of the proof scheme.
One might use Gr\"obner bases, however in cases where these are easily computable, one can equally well check $g\in I(X)$ using Gr\"obner bases directly.

\paragraph{Precision.}
In practice, the required precision might just be too high in order to do computations within reasonable time and memory constraints. 

\paragraph{Finding a suitable example.}
When constructing the example $P=(p_1,\ldots,p_n)$, the first $d$ coordinates can be chosen freely up to the height constraints.
But then $p_{d+1},\ldots,p_n$ have to be determined at least numerically up to a given precision.
The problem is to compute at least one point up to a given precision of the zero-dimensional variety $X\cap V(x_1-p_1,\ldots,x_d-p_d)$.
\emph{Numerical algebraic geometry} provides robust methods to solve this problem, see e.g. 
Sommese, Verschelde, Wampler~\cite{SommeseVerscheldeWampler05intoToNumericalAlgebraicGeometry}.
However it still can be a bottleneck in the proof by example scheme.

\subsection{Future problems}
\label{secFutureProblems}

\paragraph{Strengthening bounds.}
Can the bounds in Theorem~\ref{thmRobustMainTheorem} be tightened?
The natural aim are practical bounds, which also means that we may restrict this question to ideals $I(X)$ that appear in practice.
For example for sparse systems of polynomials $(f_1,\ldots,f_m,g)$, Krick--Pardo--Sombra~\cite[Cor. 4.12]{KrickPardoSombra01sharpEstimatesForArithmeticNullstellensatz} provide a suitable 
arithmetic Nullstellensatz with bounds that can be considerably tighter than their general bounds, and we can directly use it in place of their Theorem~\ref{thmArithmNullstellensatzBezout} above.

\paragraph{Computer-generated bounds.}
It would be useful to have an algorithm that computes very good bounds that are needed in the proof by example scheme.
The algorithm should not only implement the bounds from the theorems, but also from the proofs, which can be considerably improved by the more knowledge of the $f_i$'s and $g$ one uses.
One should not only do this for the proofs of this paper, but also perhaps more importantly for the arithmetic Nullstellensatz.

\paragraph{Sufficiently generic points.}
Are there more general sufficient criteria for a point $P$ to be sufficiently generic for the proof by example scheme to work?
Such criteria should be easy to verify for a human, or at least for a computer.

\paragraph{Non-algebraic settings.}
Are there settings other than the algebraic one where proofs by examples may reasonably exist?

\small

\bibliographystyle{plain}

\bibliography{../../mybib07}

\begin{thebibliography}{10}

\bibitem{allender2009complexityOfNumericalAnalysis}
Eric Allender, Peter B{\"u}rgisser, Johan Kjeldgaard-Pedersen, and Peter~Bro
  Miltersen.
\newblock On the complexity of numerical analysis.
\newblock {\em SIAM Journal on Computing}, 38(5):1987--2006, 2009.

\bibitem{Alon99combinatorialNullstellensatz}
Noga Alon.
\newblock Combinatorial {N}ullstellensatz.
\newblock {\em Combin. Probab. Comput.}, 8(1-2):7--29, 1999.
\newblock Recent trends in combinatorics (M\'{a}trah\'{a}za, 1995).

\bibitem{BostGilletSoule91analogueOfArithmeticBezout}
Jean-Beno\^{i}t Bost, Henri~A. Gillet, and Christophe Soul\'{e}.
\newblock Un analogue arithm\'{e}tique du th\'{e}or\`eme de {B}\'{e}zout.
\newblock {\em C. R. Acad. Sci. Paris S\'{e}r. I Math.}, 312(11):845--848,
  1991.

\bibitem{BostGilletSoule94heightsOfProjVarietiesAndPosGreenForms}
Jean-Beno\^{i}t Bost, Henri~A. Gillet, and Christophe Soul\'{e}.
\newblock Heights of projective varieties and positive {G}reen forms.
\newblock {\em J. Amer. Math. Soc.}, 7(4):903--1027, 1994.

\bibitem{Brownawell87boundsForDegreeInNullstellensatz}
W.~Dale Brownawell.
\newblock Bounds for the degrees in the {N}ullstellensatz.
\newblock {\em Ann. of Math. (2)}, 126(3):577--591, 1987.

\bibitem{Brownawell88localDiophantineNullstellenInequalities}
W.~Dale Brownawell.
\newblock Local {D}iophantine {N}ullstellen inequalities.
\newblock {\em J. Amer. Math. Soc.}, 1(2):311--322, 1988.

\bibitem{Buchberger2006phd}
Bruno Buchberger.
\newblock {B}runo {B}uchberger’s {PhD} thesis 1965: {A}n algorithm for
  finding the basis elements of the residue class ring of a zero dimensional
  polynomial ideal.
\newblock {\em J. Symbol. Comput.}, 41(3):475--511, 2006.
\newblock Engl. translation by Michael P. Abramson.

\bibitem{Chou85thesis_wusMethod}
Shang-Ching Chou.
\newblock {\em Proving and discovering geometry theorems using {W}u's method}.
\newblock ProQuest LLC, Ann Arbor, MI, 1985.
\newblock PhD thesis, University of Texas at Austin.

\bibitem{ChouGaoZhang96automatedProofsI}
Shang-Ching Chou, Xiao-Shan Gao, and Jing-Zhong Zhang.
\newblock Automated generation of readable proofs with geometric invariants.
  {I}. {M}ultiple and shortest proof generation.
\newblock {\em J. Automat. Reason.}, 17(3):325--347, 1996.

\bibitem{ChouGaoZhang96automatedProofsII}
Shang-Ching Chou, Xiao-Shan Gao, and Jing-Zhong Zhang.
\newblock Automated generation of readable proofs with geometric invariants.
  {II}. {T}heorem proving with full-angles.
\newblock {\em J. Automat. Reason.}, 17(3):349--370, 1996.

\bibitem{ChowVanDerWaerden37zurAlgebraischenGeometrie}
Wei-Liang Chow and Bartel~L. van~der Waerden.
\newblock Zur algebraischen {G}eometrie. {IX}.
\newblock {\em Math. Ann.}, 113(1):692--704, 1937.

\bibitem{Cygan98intersectionTheoryAndSperationExponent}
Ewa Cygan.
\newblock Intersection theory and separation exponent in complex analytic
  geometry.
\newblock {\em Ann. Polon. Math.}, 69(3):287--299, 1998.

\bibitem{CyganKrasinskiTworzewski99separationOfAlgebraicSetsAndLojasiewiczExponent}
Ewa Cygan, Tadeusz Krasi\'{n}ski, and Piotr Tworzewski.
\newblock Separation of algebraic sets and the {{\L}}ojasiewicz exponent of
  polynomial mappings.
\newblock {\em Invent. Math.}, 136(1):75--87, 1999.

\bibitem{DeMilloLipton78probabilisticRemarkOnAlgebraicProgramTesting}
Richard~A. {DeMillo} and Richard~J. {Lipton}.
\newblock A probabilistic remark on algebraic program testing.
\newblock {\em {Inf. Process. Lett.}}, 7:193--195, 1978.

\bibitem{sage}
The~Sage Developers.
\newblock {\em {S}age{M}ath ({V}ersion 8.8)}, 2019.
\newblock \href{http://www.sagemath.org}{http://www.sagemath.org}.

\bibitem{eisenbudHunekeVasconcelos92directMethodsPrimDecomp}
David Eisenbud, Craig Huneke, and Wolmer Vasconcelos.
\newblock Direct methods for primary decomposition.
\newblock {\em Invent. Math.}, 110(1):207--235, 1992.

\bibitem{Fukshansky10algebraicPointsSmallHeightMissingUnionOfVarieties}
Lenny Fukshansky.
\newblock Algebraic points of small height missing a union of varieties.
\newblock {\em J. Number Theory}, 130(10):2099--2118, 2010.

\bibitem{Fulton98IntersectionTheory}
William Fulton.
\newblock {\em Intersection theory}, volume~2.
\newblock Springer-Verlag, Berlin, second edition, 1998.

\bibitem{GelfandKapranovZelvinski94greenBible}
Israel~M. Gelfand, Mikhail~M. Kapranov, and Andrei~V. Zelevinsky.
\newblock {\em Discriminants, resultants, and multidimensional determinants}.
\newblock Mathematics: Theory \& Applications. Birkh\"{a}user Boston, Inc.,
  Boston, MA, 1994.

\bibitem{Hartshorne77algebraicGeometry}
Robin Hartshorne.
\newblock {\em Algebraic geometry}.
\newblock Springer-Verlag, New York-Heidelberg, 1977.
\newblock Graduate Texts in Mathematics, No. 52.

\bibitem{Hatcher02AT}
Allen Hatcher.
\newblock {\em Algebraic topology}.
\newblock Cambridge University Press, Cambridge, 2002.

\bibitem{Hickel01solutionDuneConjectureDeBerensteinYger}
Michel Hickel.
\newblock Solution d'une conjecture de {C}. {B}erenstein--{A}. {Y}ger et
  invariants de contact \`a l'infini.
\newblock {\em Ann. Inst. Fourier (Grenoble)}, 51(3):707--744, 2001.

\bibitem{geogebra5}
M.~Hohenwarter, M.~Borcherds, G.~Ancsin, B.~Bencze, M.~Blossier, J.~\'Eli\'as,
  K.~Frank, L.~G\'al, A.~Hofst\"atter, F.~Jordan, Z.~Kone\v{c}n\'y,
  Z.~Kov\'acs, E.~Lettner, S.~Lizelfelner, B.~Parisse, C.~Solyom-Gecse,
  C.~Stadlbauer, and M.~Tomaschko.
\newblock {G}eo{G}ebra 5.0, 2018.
\newblock \href{http://www.geogebra.org}{www.geogebra.org}.

\bibitem{ImpagliazzoWigderson97PeqBPPIfErequiresExpCircuits}
Russell Impagliazzo and Avi Wigderson.
\newblock {${\rm P}={\rm BPP}$} if {${\rm E}$} requires exponential circuits:
  derandomizing the {XOR} lemma.
\newblock In {\em S{TOC} '97 ({E}l {P}aso, {TX})}, pages 220--229. ACM, New
  York, 1999.

\bibitem{Jabbouri92approximationDiophantiennes}
El~Mostafa Jabbouri.
\newblock Sur un crit\`ere pour l'ind\'{e}pendance alg\'{e}brique de {P}.
  {P}hilippon.
\newblock In {\em Approximations diophantiennes et nombres transcendants
  ({L}uminy, 1990)}, pages 195--202. de Gruyter, Berlin, 1992.

\bibitem{Janivic06GCLC}
Predrag Jani\v{c}i\'{c}.
\newblock G{CLC}---a tool for constructive {E}uclidean geometry and more than
  that.
\newblock In {\em Mathematical software---{ICMS} 2006}, volume 4151 of {\em
  Lecture Notes in Comput. Sci.}, pages 58--73. Springer, Berlin, 2006.

\bibitem{JiKollarShiffman92LojasiewiczInequality}
Shanyu Ji, J\'{a}nos Koll\'{a}r, and Bernard Shiffman.
\newblock A global {{\L}}ojasiewicz inequality for algebraic varieties.
\newblock {\em Trans. Amer. Math. Soc.}, 329(2):813--818, 1992.

\bibitem{KabanetsImpagliazzo04derandomizingPITs}
Valentine Kabanets and Russell Impagliazzo.
\newblock Derandomizing polynomial identity tests means proving circuit lower
  bounds.
\newblock {\em Comput. Complexity}, 13(1-2):1--46, 2004.

\bibitem{Kapur86groebnerBasesForGeometryProblems}
Deepak Kapur.
\newblock Using {G}r\"{o}bner bases to reason about geometry problems.
\newblock {\em J. Symbolic Comput.}, 2(4):399--408, 1986.

\bibitem{kemper02calculationOfRadicalIdealInPosChar}
Gregor Kemper.
\newblock The calculation of radical ideals in positive characteristic.
\newblock {\em J. Symbol. Comput.}, 34(3):229--238, 2002.

\bibitem{Kollar99effectiveNullstellensatzForArbitraryIdeals}
J\'{a}nos Koll\'{a}r.
\newblock Effective {N}ullstellensatz for arbitrary ideals.
\newblock {\em J. Eur. Math. Soc.}, 1(3):313--337, 1999.

\bibitem{krickLogar91algorithmForRadicalIdeal}
Teresa Krick and Alessandro Logar.
\newblock An algorithm for the computation of the radical of an ideal in the
  ring of polynomials.
\newblock In {\em International Symposium on Applied Algebra, Algebraic
  Algorithms, and Error-Correcting Codes}, pages 195--205. Springer, 1991.

\bibitem{KrickPardoSombra01sharpEstimatesForArithmeticNullstellensatz}
Teresa Krick, Luis~Miguel Pardo, and Mart\'{i}n Sombra.
\newblock Sharp estimates for the arithmetic {N}ullstellensatz.
\newblock {\em Duke Math. J.}, 109(3):521--598, 2001.

\bibitem{kronecker1882grundzugeEinerArithmetischenTheorie}
Leopold Kronecker.
\newblock {\em {G}rundz{\"u}ge einer arithmetischen {T}heorie der algebraischen
  {G}r{\"o}ssen... von {L}. {K}ronecker}.
\newblock G.~Reimer, Berlin, 1882.

\bibitem{LaurentRoy01criteriaOfAlgebraicIndependence}
Michel Laurent and Damien Roy.
\newblock Criteria of algebraic independence with multiplicities and
  approximation by hypersurfaces.
\newblock {\em J. Reine Angew. Math.}, 536:65--114, 2001.

\bibitem{Lelong94mesureDeMahler}
Pierre Lelong.
\newblock Mesure de {M}ahler et calcul de constantes universelles pour les
  polyn\^{o}mes de {$n$} variables.
\newblock {\em Math. Ann.}, 299(4):673--695, 1994.

\bibitem{Littlewood14surLaDistributionDesNombresPremiers}
John~E. {Littlewood}.
\newblock {Sur la distribution des nombres premiers.}
\newblock {\em {C. R. Acad. Sci., Paris}}, 158:1869--1872, 1914.

\bibitem{Liu02AlgebraicGeometryAndArithmeticCurves}
Qing Liu.
\newblock {\em Algebraic geometry and arithmetic curves}, volume~6 of {\em
  Oxford Graduate Texts in Mathematics}.
\newblock Oxford University Press, 2002.

\bibitem{Lojasiewicz58divisionDuneDistribution}
Stanis\l{}aw \L{}ojasiewicz.
\newblock Division d'une distribution par une fonction analytique de variables
  r\'{e}elles.
\newblock {\em C. R. Acad. Sci. Paris}, 246:683--686, 1958.

\bibitem{Lojasiewicz59problemeDeLaDivision}
Stanis\l{}aw \L{}ojasiewicz.
\newblock Sur le probl\`eme de la division.
\newblock {\em Studia Math.}, 18:87--136, 1959.

\bibitem{Lojasiewicz71ensemblesSemiAnalytiques}
Stanis\l{}aw \L{}ojasiewicz.
\newblock Sur les ensembles semi-analytiques.
\newblock In {\em Actes du {C}ongr\`es {I}nternational des {M}ath\'{e}maticiens
  ({N}ice, 1970), {T}ome 2}, pages 237--241. Gauthier-Villars, Paris, 1971.

\bibitem{Mahler62someInequalitiesForPolynomialsInSeveralVariables}
Kurt Mahler.
\newblock On some inequalities for polynomials in several variables.
\newblock {\em J. London Math. Soc.}, 37:341--344, 1962.

\bibitem{MasserWustholz83fieldsOfLargeTranscendenceDegree}
David~W. Masser and Gisbert W\"{u}stholz.
\newblock Fields of large transcendence degree generated by values of elliptic
  functions.
\newblock {\em Invent. Math.}, 72(3):407--464, 1983.

\bibitem{MayrMeyer82complexityOfWordProblemsAndPolynomialIdeals}
Ernst~W. Mayr and Albert~R. Meyer.
\newblock The complexity of the word problems for commutative semigroups and
  polynomial ideals.
\newblock {\em Adv. in Math.}, 46(3):305--329, 1982.

\bibitem{MayrRitscher13dimDependentBoundsForGroebnerBases}
Ernst~W. Mayr and Stephan Ritscher.
\newblock Dimension-dependent bounds for {G}r\"{o}bner bases of polynomial
  ideals.
\newblock {\em J. Symbolic Comput.}, 49:78--94, 2013.

\bibitem{MayrToman17complexityOfMembershipProblemsOfPolynomialIdeals}
Ernst~W. Mayr and Stefan Toman.
\newblock Complexity of membership problems of different types of polynomial
  ideals.
\newblock In {\em Algorithmic and experimental methods in algebra, geometry,
  and number theory}, pages 481--493. Springer, Cham, 2017.

\bibitem{NesterenkoPhilippon01introductionToAlgebraicIndependenceTheory}
Yuri~V. Nesterenko and Patrice Philippon, editors.
\newblock {\em Introduction to algebraic independence theory}, volume 1752 of
  {\em Lecture Notes in Mathematics}.
\newblock Springer-Verlag, Berlin, 2001.
\newblock With contributions from F. Amoroso, D. Bertrand, W. D. Brownawell, G.
  Diaz, M. Laurent, Yuri V. Nesterenko, K. Nishioka, Patrice Philippon, G.
  R\'{e}mond, D. Roy and M. Waldschmidt.

\bibitem{Euclid02elements}
Euclid of~Alexandria.
\newblock {\em Euclid's {\it {E}lements}}.
\newblock Green Lion Press, Santa Fe, NM, 2002.
\newblock All thirteen books complete in one volume, The Thomas L. Heath
  translation, Edited by Dana Densmore.

\bibitem{Ore22ueberHoehereKongruenzen}
\O{}ystein Ore.
\newblock {\"Uber h\"ohere Kongruenzen.}
\newblock {Norsk matem. Forenings Skrifter 1, Nr. 7, 15 S}, 1922.

\bibitem{philippon1986criteresIndependanceAlgebrique}
Patrice Philippon.
\newblock Crit\`eres pour l'ind\'{e}pendance alg\'{e}brique.
\newblock {\em Inst. Hautes \'{E}tudes Sci. Publ. Math.}, 64:5--52, 1986.

\bibitem{philippon1991hauteursAlternativesI}
Patrice Philippon.
\newblock Sur des hauteurs alternatives. {I}.
\newblock {\em Math. Ann.}, 289(2):255--283, 1991.

\bibitem{philippon1994hauteursAlternativesII}
Patrice Philippon.
\newblock Sur des hauteurs alternatives. {II}.
\newblock {\em Ann. Inst. Fourier (Grenoble)}, 44(4):1043--1065, 1994.

\bibitem{philippon1995hauteursAlternativesIII}
Patrice Philippon.
\newblock Sur des hauteurs alternatives. {III}.
\newblock {\em J. Math. Pures Appl. (9)}, 74(4):345--365, 1995.

\bibitem{Ritt50differentialAlgebra}
Joseph~Fels Ritt.
\newblock {\em Differential {A}lgebra}.
\newblock American Mathematical Society Colloquium Publications, Vol. XXXIII.
  American Mathematical Society, New York, N. Y., 1950.

\bibitem{Schwartz80fastProbAlgoForVerificationOfPolyIdentities}
Jack~T. Schwartz.
\newblock Fast probabilistic algorithms for verification of polynomial
  identities.
\newblock {\em J. Assoc. Comput. Mach.}, 27(4):701--717, 1980.

\bibitem{ShpilkaYehudayoff10arithmeticCircuits}
Amir Shpilka and Amir Yehudayoff.
\newblock Arithmetic circuits: a survey of recent results and open questions.
\newblock {\em Found. Trends Theor. Comput. Sci.}, 5(3-4):207--388 (2010),
  2009.

\bibitem{Skewes33piMinusLiI}
Stanley Skewes.
\newblock On the difference {$\pi(x)-{\rm li}(x)$}.
\newblock {\em J. London Math. Soc.}, 8(4):277--283, 1933.

\bibitem{Skewes33piMinusLiII}
Stanley Skewes.
\newblock On the difference {$\pi(x)-{\rm li}(x)$}. {II}.
\newblock {\em Proc. London Math. Soc. (3)}, 5:48--70, 1955.

\bibitem{SommeseVerscheldeWampler05intoToNumericalAlgebraicGeometry}
Andrew~J. Sommese, Jan Verschelde, and Charles~W. Wampler.
\newblock Introduction to numerical algebraic geometry.
\newblock In {\em Solving polynomial equations}, volume~14 of {\em Algorithms
  Comput. Math.}, pages 301--335. Springer, Berlin, 2005.

\bibitem{Valiant79completenessClassesInAlgebra}
Leslie~G. Valiant.
\newblock Completeness classes in algebra.
\newblock In {\em Conference {R}ecord of the {E}leventh {A}nnual {ACM}
  {S}ymposium on {T}heory of {C}omputing}, pages 249--261. ACM, New York, 1979.

\bibitem{wikipediaProofByExample}
Wikipedia.
\newblock Proof by example.
\newblock
  {\href{https://en.wikipedia.org/wiki/Proof_by_example}{https://en.wikipedia.org/wiki/Proof\_{}by\_{}example}},
  July 2019.

\bibitem{Wu78decisionProblemAndTheoremProving}
Wen~Ts\"{u}n Wu.
\newblock On the decision problem and the mechanization of theorem-proving in
  elementary geometry.
\newblock {\em Sci. Sinica}, 21(2):159--172, 1978.

\bibitem{Wu84mechanicalTheoremProvingInElementaryGeometries}
Wen~Ts\"{u}n Wu.
\newblock Basic principles of mechanical theorem proving in elementary
  geometries.
\newblock {\em J. Systems Sci. Math. Sci.}, 4(3):207--235, 1984.

\bibitem{Zippel79probabilisticAlgoForSparsePolynomials}
Richard Zippel.
\newblock Probabilistic algorithms for sparse polynomials.
\newblock In {\em Symbolic and algebraic computation ({EUROSAM} '79,
  {I}nternat. {S}ympos., {M}arseille, 1979)}, volume~72 of {\em Lecture Notes
  in Comput. Sci.}, pages 216--226. Springer, Berlin-New York, 1979.

\end{thebibliography}

\end{document}